\input ./style/arxiv-general.cfg
\documentclass[aop,MSNbibl,seceqn,dvips]{arximspdf}
\makeatletter
   \@ifpackageloaded{graphicx}{}{\usepackage{graphicx}}
\makeatother
\doi{10.1214/15-AOP1027}% Updated by VTEXPTS2LaTeX.exe, 20.05.2015
\volume{44}
\issue{4}
\pubyear{2016}
\firstpage{2507}
\lastpage{2553}
\docsubty{FLA}

\makeatletter
\def\cal{\mathcal}
\newcommand{\textscc}[1]{\fontsize{8.36}{10.36}{\selectfont\textsc{#1}}}
\newproclaim{definition}{Definition}

\newproclaim{remark}{Remark}

\newproclaim{assumption}{Assumption}

\newcommand{\eps}{\varepsilon}

\def\dbE{\mathbb{E}}
\def\dbF{\mathbb{F}}
\def\dbL{\mathbb{L}}
\def\dbP{\mathbb{P}}
\def\dbR{\mathbb{R}}
\def\dbS{\mathbb{S}}
\def\dbQ{\mathbb{Q}}

\def\a{\alpha}
\def\b{\beta}
\def\g{\gamma}
\def\e{\varepsilon}
\def\si{\sigma}
\def\t{\tau}
\def\f{\varphi}
\def\th{\theta}
\def\D{\Delta}
\def\L{\Lambda}
\def\O{\Omega}
\def\cA{{\cal A}}
\def\cB{{\cal B}}
\def\cC{{\cal C}}
\def\cD{{\cal D}}
\def\cE{{\cal E}}
\def\cF{{\cal F}}
\def\cG{{\cal G}}
\def\cH{{\cal H}}
\def\cL{{\cal L}}
\def\cP{{\cal P}}
\def\cS{{\cal S}}
\def\cT{{\cal T}}
\def\cU{{\cal U}}

\def\ol{\overline}
\def\ul{\underline}

\def\pa{\partial}
\def\cds{\cdots}

\def\tr{\operatorname{tr}}

\def\limsup{\mathop{\overline{\operatorname{lim}}}}
\def\esup{\mathop{\operatorname{ess\mbox{-}sup}}}
\def\pa{\partial}
\def\wh{\widehat}
\def\cds{\cdots}

\def\1{\mathbf{1}}

\def\Usup{\overline{\cU}}
\def\Usub{\underline{\cU}}

\newtheorem{thmm}{Theorem}[section]
\newtheorem{lem}[thmm]{Lemma}
\newtheorem{cor}[thmm]{Corollary}
\newtheorem{prop}[thmm]{Proposition}
\newproclaim{rem}[thmm]{Remark}
\newproclaim{eg}[thmm]{Example}
\newproclaim{defn}[thmm]{Definition}
\newproclaim{assum}[thmm]{Assumption}

\makeatother

\begin{document}
\begin{frontmatter}

\title{Viscosity solutions of fully nonlinear parabolic path dependent
PDEs: Part II}
\runtitle{Viscosity solutions of fully nonlinear PDEs II}

\begin{aug}
% Corresponding author: Nizar Touzi - nizar.touzi@polytechnique.edu% Updated by VTEXPTS2LaTeX.exe, 22.05.2015 08:05
%Updated by VTEXPTS2LaTeX.exe, 20.05.2015 15:38
\author[A]{\fnms{Ibrahim}~\snm{Ekren}\ead[label=e1]{ibrahim.ekren@math.ethz.ch}},
\author[B]{\fnms{Nizar}~\snm{Touzi}\corref{}\thanksref{T1}\ead[label=e2]{nizar.touzi@polytechnique.edu}}
\and
\author[C]{\fnms{Jianfeng}~\snm{Zhang}\thanksref{T2}\ead[label=e3]{jianfenz@usc.edu}}
\thankstext{T1}{Supported by the Chair \textit{Financial Risks} of the \textit{Risk Foundation} sponsored by Soci\'et\'e
G\'en\'erale, and
the Chair \textit{Finance and Sustainable Development} sponsored by EDF and Calyon.}
\thankstext{T2}{Supported in part by NSF Grant DMS-10-08873.}
\runauthor{I. Ekren, N. Touzi and J. Zhang}
\affiliation{ETH Zurich, Ecole Polytechnique Paris and University of
Southern California}
%\dedicated{}
\address[A]{I. Ekren\\
Department of Mathematics\\
ETH Zurich\\
HG E 66.1 Ramistrasse 101 8092 Zurich\\
Switzerland\\
\printead{e1}}
\address[B]{N. Touzi\\
CMAP\\
Ecole Polytechnique Paris\\
Route de Saclay 91128\\
Palaiseau Cedex\\
France\\
\printead{e2}}
\address[C]{J. Zhang\\
Department of Mathematics\\
University of Southern California\\
3620 S. Vermont Ave, KAP 108\\
Los Angeles, California 90089-1113\\
USA\\
\printead{e3}}
\end{aug}

% HISTORY:
%
\received{\smonth{5} \syear{2013}}% Updated by VTEXPTS2LaTeX.exe,
%20.05.2015 15:38
%
\revised{\smonth{9} \syear{2014}}% Updated by VTEXPTS2LaTeX.exe,
%20.05.2015 15:38

% ABSTRACT
%
\begin{abstract}
In our previous paper [Ekren, Touzi and Zhang (2015)], we introduced a notion of viscosity
solutions for fully nonlinear path-dependent PDEs, extending the
semilinear case of Ekren et al. [\textit{Ann. Probab.} \textbf{42} (2014) 204--236],
which satisfies a partial comparison result under standard
Lipshitz-type assumptions. The main result of this paper provides a
full, well-posedness result under an additional assumption, formulated
on some partial differential equation, defined locally by freezing the
path. Namely, assuming further that such path-frozen standard PDEs
satisfy the comparison principle and the Perron approach for existence,
we prove that the nonlinear path-dependent PDE has a unique viscosity
solution. Uniqueness is implied by a comparison result.
\end{abstract}

% KEYWORDS
\end{frontmatter}

%s1 #&#
\section{Introduction}
\label{sect-Introduction}

This paper is the continuation of our accompanying papers \cite
{ETZ0,ETZ1}. The main objective of this series of three papers is the
following, fully nonlinear parabolic path-dependent partial
differential equation:
%
%e1.1 #&#
\begin{equation}
\label{PPDE-Introduction} \bigl\{-\pa_t u -G\bigl(\cdot,u,\pa_\omega u,
\pa^2_{\omega\omega} u\bigr) \bigr\}(t,\omega) =0,\qquad (t,\omega)\in[0,
T)\times\O.
\end{equation}
Here $\O$ consists of continuous paths $\omega$ on $[0, T]$ starting from
the origin, $G$ is a progressively measurable generator and the path
derivatives $\pa_t u, \pa_\omega u, \pa^2_{\omega\omega} u$ are
defined through a
functional It\^{o} formula, initiated by Dupire \cite{Dupire}; see also
Cont and Fournie \cite{CF}. Such equations were first proposed by Peng
\cite{Peng-ICM}, and they provide a convenient language for many
problems arising in non-Markovian, or say path dependent framework,
with typical examples, including martingales, backward stochastic
differential equations, second-order BSDEs and backward stochastic
PDEs. In particular, the value functions of stochastic controls and
stochastic differential games with both drift and diffusion controls
can be characterized as the solution of the corresponding path
dependent PDEs. This extends the classical results in Markovian
framework to non-Markovian ones. We refer to \cite{ETZ1} and \cite{PZ}
for these connections.

A path dependent PDE can rarely have a classical solution. We thus turn
to the notion of viscosity solutions, which had great success in the
finite dimensional case. There have been numerous publications on
viscosity solutions of PDEs, both in theory and in applications, and we
refer to the classical references \cite{CIL} and \cite{FS}. In our
infinite dimensional case, the major difficulty is that the underlying
state space $\O$ is not locally compact, and thus many tools from the
standard PDE viscosity theory do not apply to the present context. In
our earlier paper \cite{EKTZ}, which studies semilinear path-dependent
PDEs, we replace the pointwise extremality in the standard definition
of viscosity solution in PDE literature with the corresponding
extremality in the context of an optimal stopping problem under a
nonlinear expectation $\cE$. More precisely, we introduce a set of
smooth test processes $\varphi$, which are tangent from above or from
below, to the processes of interest $u$ in the sense of the following
nonlinear optimal stopping problems:
%
%e1.2 #&#
\begin{eqnarray}
\label{optimal-Introduction} \sup_{\tau}\overline{\cE}\bigl[(
\varphi-u)_{\tau}\bigr], \qquad\inf_{\tau}\underline{\cE}
\bigl[(\varphi-u)_{\tau}\bigr]
\nonumber
\\[-8pt]
\\[-8pt]
\eqntext{\mbox{where } \displaystyle\overline{\cE}:=\sup
_{\dbP\in\cP}\dbE^\dbP, \displaystyle\underline{\cE}:=\inf
_{\dbP\in\cP}\dbE^\dbP.}
\end{eqnarray}
Here $\tau$ ranges over a convenient set of stopping times, and $\cP$
is an appropriate set of probability measures. The replacement of the
pointwise tangency by the tangency in the sense of the last optimal
stopping problem is the key ingredient needed to bypass the local
compactness of the underlying space in the standard viscosity solution
theory (or the Hilbert structure, which allows us to access local
compactness by finite realization approximation of the space). Indeed,
the Snell envelope characterization of the solution of the optimal
stopping problem allows us to find a ``point of tangency.''
Interestingly, the structure of the underlying space does not play any
role, and the standard first and second-order conditions of maximality
in the standard optimization theory has the following beautiful
counterpart in the optimal stopping problem: the supermartingale
property (negative drift; notice that drift is related to the second
derivative) of the Snell envelope and the martingale property (zero
drift) up the optimal stopping time (first hitting of the
obstacle/reward process).

In \cite{EKTZ}, we proved existence and uniqueness of viscosity
solutions for semilinear path-dependent PDEs. In particular, the unique
viscosity solution is consistent with the solution to the corresponding
backward SDE.

In \cite{EKTZ}, all probability measures in the class $\cP$ are
equivalent, and consequently $\cP$ is dominated by one measure. In our
fully nonlinear context, the class $\cP$ becomes nondominated,
consisting of mutually singular measures induced by certain
linearization of the nonlinear generator $G$. This causes another major
difficulty of the project: the dominated convergence theorem fails
under $\ol\cE^\cP$. To overcome this, we need some strong regularity
for the involved processes, and thus we require some rather
sophisticated estimates. In particular, the corresponding optimal
stopping problem becomes very technical and is established in a
separate paper \cite{ETZ0}. We remark that the weak compactness of the
class $\cP$ plays a very important role in these arguments.

In \cite{ETZ1}, we introduced the appropriate class $\cP$ for fully
nonlinear path dependent PDEs (\ref{PPDE-Introduction}) and the
corresponding notion of viscosity solutions. We investigated the
connection between our new notion and many other equations in the
existing literature of stochastic analysis, for example, backward SDEs,
second-order BSDEs and backward SPDEs. Moreover, we proved some basic
properties of viscosity solutions, including the partial comparison
principle; that is, for a viscosity subsolution $u^1$ and a classical
supersolution $u^2$, if $u^1_T \le u^2_T$, then $u^1(t,\omega) \le
u^2(t,\omega
)$ for all $(t,\omega) \in[0, T]\times\O$.

In this paper we prove our main result, the comparison principle of
viscosity solutions. That is, for a viscosity subsolution $u^1$ and a
viscosity supersolution $u^2$, if $u^1_T \le u^2_T$, then
$u^1(t,\omega)
\le u^2(t,\omega)$ for all $(t,\omega) \in[0, T]\times\O$.
Again, due to the
lack of local compactness and now also due to our new definition of
viscosity solutions, the standard approach in PDE literature, namely
the doubling variable technique combined with Ishii's lemma, does not
seem to work here. Our strategy is as follows: We start from the above
partial comparison established in \cite{ETZ1}, but we slightly weaken
the smooth requirement of the classical (semi-)solutions. Let $\ol u$
denote the infimum of the classical supersolution and $\ul u$, the
supremum of classical subrsolutions, satisfying appropriate terminal
conditions. Then the partial comparison implies $u^1 \le\ol u$ and
$\ul u \le u^2$. Thus the comparison will be a direct consequence of
the following claim:
%
%e1.3 #&#
\begin{equation}
\label{Perron-Introduction} \ol u = \ul u.
\end{equation}

Then clearly our focus is (\ref{Perron-Introduction}). We first remark
that due to the failure of the dominated convergence theorem under our
new $\ol\cE^\cP$, the approach in \cite{EKTZ} does not work here. In
this paper, we shall follow the alternative approach proposed in \cite
{ETZ1}, Section~7, which is also devoted to semilinear path-dependent
PDEs. However, as explained in \cite{ETZ1}, Remark~7.7, there are
several major difficulties in the fully nonlinear context, and novel
ideas are needed.

Note that (\ref{Perron-Introduction}) is more or less equivalent to
constructing some classical supersolution $\ol u^\e$ and classical
subsolution $\ul u^\e$, for any $\e>0$, such that $\lim_{\e\to
0}[\ol
u^\e- \ul u^\e] = 0$. Our main tool is the following local path-frozen
PDE: for any $(t, \omega)\in[0, T)\times\O$,
%
%e1.4 #&#
\begin{eqnarray}
\label{Frozen-Introduction}  -\pa_t v (s,x) - g^{t,\omega}\bigl(s, v, D v,
D^2 v\bigr) =0,\nonumber\\
\\[-20pt]
\eqntext{s\in[t, t+\e], x\in\dbR^d \mbox{ such
that } |x|\le\e,}
\\
\eqntext{\mbox{where } g^{t,\omega}(s, y,z,\g):= G(s, \omega_{\cdot
\wedge t}, y, z,
\g). }
\end{eqnarray}
Here $D$ and $D^2$ denote the gradient and Hessian of $v$ with respect
to $x$, respectively, and we emphasize that $g^{t,\omega}$ is a
deterministic function, and thus (\ref{Frozen-Introduction}) is a
standard PDE. We shall assume that the above PDE has a unique viscosity
solution (in standard sense), which can be approximated by classical
subsolutions and classical supersolutions. One sufficient condition is
that after certain smooth mollification of $g^{t,\omega}$, the above local
PDE with smooth boundary condition has a classical solution. We then
use this classical solution to construct the desired $\ol u^\e$ and
$\ul u^\e$.

We remark that this approach is very much like Perron's approach in
standard PDE viscosity theory. However, there are two major
differences: First, in the standard Perron approach, $\ol u$ and $\ul
u$ are the extremality of viscosity semi-solutions, while here they are
the extremality of classical semi-solutions. This requires the
smoothness of the above $\ol u^\e$ and $\ul u^\e$ and thus makes their
construction harder. More importantly, the standard Perron approach
assumes the comparison principle and uses it to obtain the existence of
viscosity solutions, while we use (\ref{Perron-Introduction}) to prove
the comparison principle. Thus the required techniques are quite different.

Once we have the comparison principle, then following the idea of the
standard Perron approach, we see $\ol u= \ul u$ is indeed the unique
viscosity solution of the path-dependent PDE, so we have both existence
and uniqueness. Our result covers quite general classes of
path-dependent PDEs, including those not accessible in the existing
literature of stochastic analysis. One particular application is the
existence of the game value for a path-dependent zero sum stochastic
differential game, due to our well-posedness result of the
path-dependent Bellman--Isaacs equation; see Pham and Zhang \cite{PZ}.
We also refer to Henry-Labordere, Tan and Touzi \cite{HTT} and Zhang
and Zhuo \cite{ZZ} for applications of our results to numerical methods
for path-dependent PDEs.

We also note that there is potentially an alternative way to prove the
comparison principle. Roughly speaking, given a viscosity subsolution
$u^1$ and a viscosity supersolution $u^2$, if we could find certain
smooth approximations $u^{i, \e}$, close to $u^i$, such that $u^{1,\e}$
is a classical subsolution and $u^{2,\e}$ is a classical supersolution,
then it follows from partial comparison (actually classical comparison)
that $u^{1,\e}\le u^{2,\e}$, which leads to the desired comparison
immediately by passing $\e$ to~$0$. Indeed, in PDE literature the
convex/concave convolution plays this role. However, in the
path-dependent setting, we did not succeed in finding appropriate
approximations $u^{i,\e}$ which satisfy the desired semi-solution
property. In our current approach, instead of approximating the
(semi-)solution directly, we approximate the path-frozen PDE by
mollifying its generator $g^{t,\omega}$. The advantage of our
approach is
that provided the mollified path-frozen PDE has a classical solution,
it will be straightforward to check that the constructed $\ol u^\e$ and
$\ul u^\e$ are classical semi-solutions.

The price of our approach, however, is that we need classical solutions
of fully nonlinear PDEs. Partially for this purpose, in the present
paper we assume that $G$ is uniformly nondegenerate, which is
undesirable in viscosity theory, and for path-dependent Bellman--Isaacs
equations, we can only deal with the lower dimensional ($d=1$ or $2$)
problems. We shall investigate these important problems and explore
further possible direct approximations of $u^i$ as mentioned above, in
our future research.

The rest of the paper is organized as follows. Section~\ref
{sect-preliminary} introduces the general framework and recalls the
definition of viscosity solutions introduced in our accompanying paper
\cite{ETZ1}. Section~\ref{sect:assumptions} collects all assumptions
needed throughout the paper. The main results are stated in
Section~\ref{sect-wellposedness}, where we also outline strategy of
proof. In
particular, the existence and comparison results follow from the
partial comparison principle, the consistency of the Perron approach
and the viscosity property of the postulated solution of the PPDE.
These crucial results are proved in Sections~\ref
{sect-partialcomparison}, \ref{sect-baru} and \ref{sect:viscosity},
respectively. Finally, Section~\ref{sect-Perron} provides some
sufficient conditions for our main assumption, under which our
well-posedness result is established, together with some concluding remarks.

%s2 #&#
\section{Preliminaries}
\label{sect-preliminary}

In this section, we recall the setup and the notation of \cite{ETZ1}.

%s2.1 #&#
\subsection{The canonical spaces}

Let $\O:=  \{\omega\in C([0,T], \dbR^d)\dvtx \omega_0={\mathbf0}
\}$ be the set
of continuous paths starting from the origin, $B$, the canonical
process, $\dbF$, the natural filtration generated by $B$, $\dbP_0$, the
Wiener measure and $\L:= [0,T]\times\O$. Here and in the sequel, for
notational simplicity, we use ${\mathbf0}$ to denote vectors, matrices or
paths with appropriate dimensions whose components are all equal to
$0$. Let $\dbS^d$ denote the set of $d\times d$ symmetric matrices, and
\begin{eqnarray*}
x \cdot x'& :=& \sum_{i=1}^d
x_i x'_i\qquad \mbox{for any } x, x'
\in\dbR^d,
\nonumber
\\[-8pt]
\\[-8pt]
\nonumber
 \g\dvtx \g'& :=& \tr\bigl[\g\g'
\bigr] \qquad\mbox{for any } \g, \g'\in\dbS^d.
\end{eqnarray*}
We define a semi-norm on $\O$ and a pseudometric on $\L$ as follows:
for any $(t, \omega), ( t', \omega') \in\L$,
%
%e2.1 #&#
\begin{equation}\qquad
\label{rho} \|\omega\|_{t} := \sup_{0\le s\le t} |
\omega_s|,\qquad  {\mathbf{d}}_\infty \bigl((t, \omega),\bigl(
t', \omega'\bigr) \bigr) := \bigl|t-t'\bigr| + \bigl\|
\omega_{\cdot \wedge t} - \omega'_{\cdot \wedge t'} \bigr\|_T.
\end{equation}
Then $(\O, \|\cdot\|_{T})$ is a Banach space, and $(\L, {\mathbf
{d}}_\infty)$ is
a complete pseudometric space. We shall denote by $\dbL^0(\cF_T)$ and
$\dbL^0(\L)$ the collection of all $\cF_T$-measurable random variables
and $\dbF$-progressively measurable processes, respectively. Let
$C^0(\L
)$ [resp., $\mathrm{UC}(\L)$] be the subset of $\dbL^0(\L)$ whose elements are
continuous (resp., uniformly continuous) in $(t,\omega)$ under
${\mathbf{d}}_\infty
$. The corresponding subsets of bounded processes are denoted by
$C^0_b(\L)$ and $\mathrm{UC}_b(\L)$. Finally, $\dbL^0(\L, \dbR^d)$ denote the
space of $\dbR^d$-valued processes with entries in $\dbL^0(\L)$, and we
define similar notation for the spaces $C^0$, $C^0_b$, $\mathrm{UC}$ and $\mathrm{UC}_b$.

We next introduce the shifted spaces. Let $0\le s\le t\le T$.
\begin{itemize}[$-$]
\item[$-$] Let $\O^t:=  \{\omega\in C([t,T], \dbR^d)\dvtx \omega_t ={\mathbf
0} \}$ be the
shifted canonical space; $B^{t}$ the shifted canonical process on
$\O^t$; $\dbF^{t}$ the shifted filtration generated by $B^{t}$, $\dbP
^t_0$ the Wiener measure on $\O^t$, and $\L^t := [t,T]\times\O^t$.

\item[$-$] Define $\|\cdot\|^s_t$ on $\O^s$ and ${\mathbf{d}}^s_\infty$ on
$\L^s$ in the
spirit of (\ref{rho}), and the sets $\dbL^0(\L^t)$ etc. in an obvious
way. %, $C_b^0(\L^t)$, \mathrm{UC}$(\L^t)$, and \mathrm{UC}$_b(\L^t)$
%and Definition~\ref{defn-spaceC0}.

\item[$-$] For $\omega\in\O^s$ and $\omega'\in\O^t$, define the
concatenation path $\omega
\otimes_{t} \omega'\in\O^s$ by
\[
\bigl(\omega\otimes_t \omega'\bigr) (r) :=
\omega_r\1_{[s,t)}(r) + \bigl(\omega_{t} +
\omega'_r\bigr)\1_{[t, T]}(r)\qquad\mbox{for all }
r\in[s,T].
\]

\item[$-$] Let $\xi\in\dbL^0(\cF^s_T)$ and $X\in\dbL^0(\L^s)$. For $(t,
\omega)
\in\L^s$, define $\xi^{t,\omega} \in\dbL^0(\cF^t_T)$ and
$X^{t,\omega}\in\dbL
^0(\L^t)$ by
% an $\cF^{s}_{T}$-measurable random variable $\xi$, an $
%\dbF^{s}$-progressively measurable process $X$ on $\O^s$, and $t
%\in(s,T]$, define the shifted $\cF^{t}_{T}$-measurable random variable
%$\xi^{t,\o}$ and $\dbF^{t}$-progressively measurable process $X^{t,
%\o}$ on $\O^t$ by:
%
\[
\xi^{t, \omega}\bigl(\omega'\bigr) :=\xi\bigl(\omega
\otimes_t \omega'\bigr), \qquad X^{t, \omega}\bigl(
\omega'\bigr) := X\bigl(\omega \otimes_t
\omega'\bigr),\qquad \mbox{for all } \omega'\in
\O^t.
\]
\end{itemize}

It is clear that for any $(t,\omega) \in\L$ and any $u\in C^0(\L
)$, we
have $u^{t,\omega} \in C^0(\L^t)$. The other spaces introduced
before enjoy
the same property.

We denote by $\cT$ the set of $\dbF$-stopping times, and by $\cH
\subset
\cT$, the subset of those hitting times $\textsc{h}$ of the form
%
%e2.2 #&#
\begin{equation}
\label{cT} \textsc{h}:= \inf\{t\dvtx B_t \notin O\} \wedge
t_0, % = \inf\{ t: d(
%\o_t, O^c) =
%0\} \wedge t_0,
\end{equation}
for some $0< t_0\le T$, and some open and convex set $O \subset\dbR^d$
containing ${\mathbf0}$.
The set $\cH$ will be important for our optimal stopping problem, which
is crucial for the comparison and the stability results. We note that
$\textsc{h}= t_0$ when $O=\dbR^d$, and for any $\textsc{h}\in\cH$,
%
%e2.3 #&#
\begin{equation}
\label{che} 0<\textsc{h}_\e\le\textsc{h} \mbox{ for $\e$ small
enough, where } \textsc{h}_\e:= \inf\bigl\{ t\ge0\dvtx |B_t|=
\e\bigr\} \wedge\e.\hspace*{-25pt}
\end{equation}
%
%Moreover,
%\beaa
%\ch>0, \mbox{$\ch$ is lower semicontinuous, and $\ch_1 \wedge\ch_2
%\in\cH$ for any $\ch_1, \ch_2\in\cH$}.
%\eeaa
Define $\cT^t$ and $\cH^t$ in the same spirit. For any $\t\in\cT$
(resp., $\textsc{h}\in\cH$) and any $(t,\omega)\in\L$ such that
$t<\t(\omega)$ [resp.,
$t<\textsc{h}(\omega)$], it is clear that $\t^{t,\omega}\in\cT
^t$ (resp., $\textsc{h}^{t,\omega
}\in\cH^t$).

Finally, the following types of regularity will be important in our framework:
%
%de2.1 #&#
\begin{defn}
\label{defn-spaceC0} Let $u \in\dbL^0 (\L)$.

\begin{longlist}[(iii)]
\item[(i)] We say $u$ is right continuous in $(t,\omega)$ under ${\mathbf
{d}}_\infty$ if
for any $(t,\omega) \in\L$ and any $\e>0$, there exists $\delta
>0$ such that
for any $(s, \tilde\omega)\in\L^t$ satisfying ${\mathbf
{d}}_\infty((s, \tilde\omega),
(t, {\mathbf0})) \le\delta$, we have $|u^{t,\omega}(s,\tilde\omega
) - u(t,\omega)|\le\e$.

\item[(ii)] We say $u\in\Usub$ if $u$ is bounded from above, right continuous
in $(t,\omega)$ under ${\mathbf{d}}_\infty$ and there exists a
modulus of
continuity function $\rho$ such that for any $(t,\omega),
(t',\omega')\in\L$,
%
%e2.4 #&#
\begin{equation}
\label{USC} u(t,\omega) - u\bigl(t', \omega'\bigr)
\le \rho \bigl({\mathbf{d}}_\infty \bigl((t,\omega),
\bigl(t',\omega'\bigr) \bigr) \bigr)\qquad \mbox{whenever } t
\le t'.
\end{equation}
\item[(iii)] We say $u\in\Usup$ if $-u \in\Usub$.
\end{longlist}
\end{defn}

The progressive measurability of $u$ implies that $u(t,\omega) =
u(t,\omega_{\cdot
\wedge t})$, and it is clear that $\Usub \cap \Usup= \mathrm{UC}_b(\L)$.
We also recall from \cite{ETZ0} Remark~3.2 that condition~(\ref{USC})
implies that $u$ has left-limits and positive jumps.

%s2.2 #&#
\subsection{Capacity and nonlinear expectation}

For every constant $L>0$, we denote by $\cP_L$ the collection of all
continuous semimartingale measures $\dbP$ on $\O$ whose drift and
diffusion characteristics are bounded by $L$ and $\sqrt{2L}$,
respectively. To be precise, let $\tilde\O:=\O^3$ be an enlarged
canonical space, $\tilde B:= (B, A, M)$ be the canonical processes and
$\tilde\omega= (\omega, a, m)\in\tilde\O$ be the paths. A
probability measure
$\dbP\in\cP_L$ means that there exists an extension $\dbQ$ of $\dbP$
on $\tilde\O$ such that
%
%e2.5 #&#
\begin{eqnarray}\qquad
\label{cPL} B&=&A+M \qquad\mbox{$A$ is
absolutely continuous, $M$ is a martingale},
\nonumber\\
\bigl|\alpha^\dbP\bigr|&\le& L,\qquad\displaystyle \frac{1}{2}\tr\bigl(\bigl(
\beta^\dbP\bigr)^2\bigr)\le L\qquad \mbox {where }
\a^\dbP_t:= {d A_t\over dt},
\b^\dbP_t
:= \sqrt{d\langle
M\rangle
_t\over dt}\\
\eqntext{\dbQ\mbox{-a.s.}}
\end{eqnarray}
Similarly, for any $t\in[0, T)$, we may define $\cP_L^t$ on $\O^t$ and
$\cP^t_\infty:= \bigcup_{L>0} \cP^t_L$.

The set $\cP^t_L$ induces the following capacity:
%
%e2.6 #&#
\begin{equation}
\cC^L_t[A] := \sup_{\dbP\in\cP_L^t}\dbP[A],\qquad
\mbox{for all }A\in\cF^t_T.
\end{equation}
We denote by $\dbL^1(\cF^t_T,\cP^t_L)$ the set of all $\cF
^t_T$-measurable r.v. $\xi$ with\break $\sup_{\dbP\in\cP^t_L}\dbE
^{\dbP}[|\xi
|]<\infty$. The following nonlinear expectation will play a crucial role:
%
%e2.7 #&#
\begin{eqnarray}
\label{cE} \overline{\cE}^L_t[\xi] = \sup
_{\dbP\in\cP^t_L}\dbE^{\dbP}[\xi]\quad \mbox{and}\quad \underline{
\cE}^L_t[\xi] = \inf_{\dbP\in\cP^t_L}
\dbE^{\dbP}[\xi] = -\overline{\cE}^L_t[-\xi]
\nonumber
\\[-8pt]
\\[-8pt]
\eqntext{\mbox{for all }\xi\in\dbL^1\bigl(\cF^t_T,
\cP^t_L\bigr).}
\end{eqnarray}

%de2.2 #&#
\begin{defn}
Let $X\in\dbL^0(\L)$ satisfy $X_t\in\dbL^1(\cF_t,\cP_L)$ for all
$0\le
t\le T$. We say that $X$ is an $\overline{\cE}^L$-supermartingale
(resp., submartingale, martingale) if, for any $(t,\omega)\in\L$ and
any $\t
\in\cT^t$, $\overline{\cE}^L_t[X^{t,\omega}_\t]\le$ (resp., $\ge
,=$) $X_t(\omega)$.
\end{defn}

We now state the Snell envelope characterization of optimal stopping
under the above nonlinear expectation operators. Given a bounded
process $X\in\dbL^0(\L)$, consider the nonlinear optimal stopping problem
%
%e2.8 #&#
\begin{eqnarray}
\label{cS} \overline{\cS}^L_t[X](\omega) := \sup
_{\t\in\cT^t} \overline\cE^L_t
\bigl[X^{t,\omega}_{\t} \bigr]\quad \mbox{and}\quad \underline{
\cS}^L_t[X](\omega) := \inf_{\t\in\cT^t}
\underline\cE^L_t \bigl[X^{t,\omega}_{\t}
\bigr],
\nonumber
\\[-8pt]
\\[-8pt]
 \eqntext{(t,\omega)\in\L.}
\end{eqnarray}
By definition, we have $\overline{\cS}^L[X]\ge X$ and $\overline{\cS
}^L_T[X] = X_T$.

%th2.3 #&#
\begin{thmm}[(\cite{ETZ0})]\label{thmm-optimal}
Let $X\in\Usub$ be bounded, $\textsc{h}\in\cH$ and set $\wh
X_t:=X_t\1_{\{
t<\textscc{h}\}}+X_{\textscc{h}-}\1_{\{t\ge\textscc{h}\}}$. Define
\[
Y := \overline\cS^L [\wh X ] \quad\mbox{and}\quad \t^*:=\inf\{t\ge0\dvtx
Y_t=\wh X_t\}.
\]
Then $Y_{\t^*} = \wh X_{\t^*}$, $Y$ is an $\overline\cE
^L$-supermartingale on $[0, \textsc{h}]$, and an $\overline\cE^L$-martingale
on $[0, \t^*]$. Consequently, $\t^*$ is an optimal stopping time.
\end{thmm}

We remark that the nonlinear Snell envelope $Y$ is continuous in $[0,
\textsc{h})$ and has left limit at $\textsc{h}$. However, in general
$Y$ may have a
negative left jump at $\textsc{h}$.

%s2.3 #&#
\subsection{The path derivatives}

We define the path derivatives via the functional It\^{o} formula,
initiated by Dupire \cite{Dupire}.

%de2.4 #&#
\begin{defn}
\label{defn-spaceC12} We say $u\in C^{1,2}(\L)$ if $u\in C^0(\L)$, and
there exist $\pa_t u \in C^0(\L)$, $\pa_\omega u \in C^0(\L, \dbR
^d)$, $\pa
^2_{\omega\omega} u\in C^0(\L, \dbS^d)$ such that for any %$(t,\o)
%\in[0,T)
%\times\O$ and any
$\dbP\in\cP^0_\infty$, $u$ is a $\dbP$-semimartingale satisfying
%
%e2.9 #&#
\begin{equation}
\label{Ito} d u = \pa_t u \,dt+ \pa_\omega u \cdot d
B_t + \tfrac{1}2 \pa ^2_{\omega\omega} u \dvtx d
\langle B\rangle_t,\qquad  0\le t\le T, \dbP\mbox{-a.s.} %d u^{t,\o}(s, B^t) = (\pa_t u)^{t,\o} ds+ (\pa_\o u)^{t,\o} \cd d
%B^t_s + \frac12(\pa^2_{\o\o} u)^{t,\o} : d \la B^t\ra_s, t\le s\le
%T, \dbP\mbox{-a.s.}
\end{equation}
%
%Moreover, we say $u\in C^{1,2}_b(\L)$ if $u$, $\pa_t u$, $\pa_\o u$,
%and $\pa^2_{\o\o} u$ are all bounded.
\end{defn}
We remark that the above $\pa_t u$, $\pa_\omega u$ and $\pa
^2_{\omega\omega} u$, if
they exist, are unique, and thus are called the time derivative,
first-order and second-order space derivatives of $u$, respectively. In
particular, it holds that
%
%e2.10 #&#
\begin{equation}
\label{pat} \pa_t u(t,\omega) = \lim_{h\downarrow0}
{1\over h}\bigl[u (t + h, \omega_{\cdot
\wedge t} ) - u (t, \omega )
\bigr].
\end{equation}
We refer to \cite{ETZ1}, Remark~2.9, and \cite{BMZ}, Remarks 2.3, 2.4,
for various discussions on these path derivatives, especially on their
comparison with Dupire's original definition. See also Remark~\ref
{rem-derivative} below. We define $C^{1,2}(\L^t)$ similarly. It is
clear that, for any $(t,\omega)$ and $u\in C^{1,2}(\L)$, we have
$u^{t,\omega
}\in C^{1,2}(\L^t)$, and $\pa_\omega(u^{t,\omega}) = (\pa
_\omega u)^{t,\omega}$, $\pa
^2_{\omega\omega} (u^{t,\omega}) = (\pa^2_{\omega\omega}
u)^{t,\omega}$.

For technical reasons, we shall extend the space $C^{1,2}(\L)$ slightly
as follows.

%de2.5 #&#
\begin{defn}
\label{defn-barC}
Let $t\in[0,T]$, $u\dvtx \L^t \to\dbR$. We say $u\in\ol C^{1,2}(\L^t)$
if there exist an increasing sequence of $\{\textsc{h}_i, i\ge1\}
\subset\cT
^t$, a partition $\{E^i_j, j\ge1\}\subset\cF_{\textscc{h}_i}^t$ of
$\O^t$ for
each $i$, a constant $n_{i}\ge1$ for each $i$, and $\f^i_{jk}\in\mathrm{
UC}_b(\L)$, $\psi^i_{jk} \in C^{1,2}(\L) \cap \mathrm{UC}_b(\L)$ for each
$(i,j)$ and $1\le k \le n_{i}$, such that, denoting $\textsc{h}_0 := t$,
$E^0_1:=\O^t$:
\begin{longlist}[(iii)]
\item[{(i)}] for each $i$ and $\omega$, $\textsc{h}^{\textscc
{h}_i,\omega}_{i+1} \in\cH^{\textscc{h}_i(\omega
)}$ whenever $\textsc{h}_i(\omega) < T$, the set $\{i\dvtx \textsc
{h}_i(\omega) < T\}$ is finite
for each $\omega$ and $\lim_{i\to\infty} \cC^L_s[\textsc
{h}_i^{s,\omega} < T] =0$ for
any $(s,\omega)\in\L^t$ and $L>0$;

\item[{(ii)}] for each $(i,j)$, $\omega, \omega'\in E^i_j$ such that
$\textsc{h}_i(\omega)
\le\textsc{h}_i(\omega')$, it holds for all $\tilde\omega\in\O$,
%
%e2.11 #&#
\begin{equation}
\label{chO} 0\le\textsc{h}_{i+1}\bigl(\omega'
\otimes_{\textscc{h}_i(\omega')} \tilde\omega\bigr)- \textsc{h}_{i+1}(\omega
\otimes _{\textscc{h}_i(\omega)} \tilde\omega) \le\textsc{h}_i\bigl(\omega
'\bigr) - \textsc{h}_i(\omega);
\end{equation}
here we abuse the notation that $(\omega\otimes_s \tilde\omega)_r
:= \omega_r \1
_{[t, s)}(r) + (\omega_s + \tilde\omega_{r-s}) \1_{[s, T]}(r)$;

\item[{(iii)}] for each $i$, $\f^i_{jk}$, $\psi^i_{jk}$, $\pa_t \psi
^i_{jk}$, $\pa_\omega\psi^i_{jk}$, $\pa_{\omega\omega}^2\psi
^i_{jk}$ are uniformly
bounded, and $\f^i_{jk}$, $\psi^i_{jk}$ are uniformly continuous,
uniformly in $(j, k)$ (but may depend on $i$);
%for each $i$, $\f^i_{jk}$ and $\psi^i_{jk}$ have a modulus of
%continuity function $\rho^i$, and $\pa_t \psi^i_{jk}$, $\pa_\o
%\psi^i_{jk}$, $\pa_{\o\o}^2\psi^i_{jk}$ have a bound $C_i$, where $
%\rho^i$ and $C_i$ may depend on $i$ but do not depend on $(j,k)$.

%For each $(i, j)$, there exist $\f^i_j\in{ UC}_b(\L)$ and $\psi^i_j
%\in C^{1,2}(\L) \cap UC_b(\L)$ such that such that $\f^i_j(t, B_{\cd
%\wedge\ch_i})$ is differentiable in $t$ on $[\ch_i, T)$, with a
%modulus of continuity function $\rho^i$, which may depend on $i$ but
%does not depend on $j$, such that
\item[{(iv)}] $u$ is continuous in $t$ on $[0,T]$, and for each $i$,
$\omega\in
\O$ and $ \textsc{h}_i(\omega) \le s\le\textsc{h}_{i+1}(\omega)$,
%
%e2.12 #&#
\begin{eqnarray}
\label{barC}&& u(s,\omega)
\nonumber
\\[-8pt]
\\[-8pt]
\nonumber
&&\qquad= \sum_{j\ge1}\sum
_{k=1}^{n_{i} } \bigl[\f ^i_{jk}
\bigl(\textsc{h}_i(\omega), \omega\bigr) \psi^i_{jk}
\bigl(s-\textsc{h}_i(\omega), \omega_{\textscc{h}_i(\omega
)+s} -
\omega_{\textscc{h}_i(\omega)} \bigr) \bigr] \1_{E^i_j}.
\end{eqnarray}
\end{longlist}
\end{defn}
The main idea of the above space is that the processes in $ \ol
C^{1,2}(\L^t)$ are piecewise smooth. However, purely for technical
reasons, we require rather technical conditions. For example, (\ref
{chO}) and (\ref{barC}) are mainly needed for Proposition~\ref
{prop-comparison} below. We remark that these technical requirements
may vary from time to time. In particular, the space here requires a
more specific structure than the corresponding space in \cite{EKTZ} and
that in \cite{ETZ1} Section~7, both dealing with semilinear PPDEs.
Nevertheless, by abusing the notation slightly, we still denote it as $
\ol C^{1,2}(\L^t)$.

%We remark that the derivatives in (iv) above exist due to \reff{barC},
%and
Let $u\in\overline C^{1,2}(\L^t)$. One may easily check that
$u^{s,\omega}
\in\overline C^{1,2}(\L^s)$ for any $(s,\omega)\in\L^t$. For any
$\dbP\in
\cP^t_\infty$, it is clear that the process $u$ is a local $\dbP
$-semimartingale on $[t,T]$ and a $\dbP$-semimartingale on $[t,\textsc{h}_i]$
for all $i$,
and
%
%e2.13 #&#
\begin{equation}\qquad
\label{Ito2} d u_s= \pa_t u_s \,ds+
\tfrac{1}2\pa^2_{\omega\omega} u_s \dvtx d\bigl
\langle B^t\bigr\rangle_s+ \pa_\omega
u_s \cdot \,dB^t_s,\qquad  t\le s< T, \dbP\mbox{-a.s.}
\end{equation}
By setting $\textsc{h}_1 := T$, $n_{0}:=1$, $\f^0_{11} := 1$ and
$\psi
^0_{11}:= u$, we see that $C^{1,2}(\L^t)\subset\ol C^{1,2}(\L^t)$.

%s2.4 #&#
\subsection{Fully nonlinear path dependent PDEs}
%\label{sect-PPDE}
%%\setcounter{equation}{0}

Following the accompanying paper \cite{ETZ1}, we continue our study of
the following fully nonlinear parabolic path-dependent partial
differential equation (PPDE, for short):
%
%e2.14 #&#
\begin{equation}\qquad
\label{PPDE} \cL u (t,\omega) := \bigl\{-\pa_t u - G \bigl(\cdot,
u, \pa_\omega u, \pa^2_{\omega\omega} u \bigr)\bigr\} (t,
\omega) = 0,\qquad (t,\omega)\in\L,
\end{equation}
where the generator $G \dvtx \L\times\dbR\times\dbR^d \times\dbS^d
\rightarrow\dbR$ satisfies the conditions reported in Section~\ref
{sect:assumptions}.

For any $u\in\dbL^0(\L)$, $(t,\omega) \in[0, T)\times\O$ and
$L>0$, define
%
%e2.15 #&#
\begin{eqnarray}
\label{cA} %
 \underline\cA^{ L} u(t,
\omega) &:=& \bigl\{\f\in C^{1,2}\bigl(\L^{ t}\bigr) \dvtx \bigl(
\f-u^{t,\omega}\bigr)_t = 0 = \underline\cS^L_t
\bigl[\bigl(\f-u^{t,\omega}\bigr)_{\cdot\wedge\textscc{h}} \bigr]\nonumber\\
 &&\hspace*{160pt}\mbox{for some }
\textsc{h}\in\cH^t \bigr\},
\nonumber
\\[-8pt]
\\[-8pt]
\nonumber
\overline\cA^{ L} u(t,\omega) &:=& \bigl\{\f\in C^{1,2}\bigl(
\L^{ t}\bigr) \dvtx \bigl(\f-u^{t,\omega}\bigr)_t =0=
\overline\cS^L_t \bigl[\bigl(\f-u^{t,\omega}
\bigr)_{\cdot\wedge\textscc{h}} \bigr] \\
&&\hspace*{160pt}\mbox{for some } \textsc{h}\in\cH^t \bigr
\},\nonumber
\end{eqnarray}
where $\overline{\cS}^L$ and $\underline{\cS}^L$ are the nonlinear
Snell envelopes defined in (\ref{cS}).

%de2.6 #&#
\begin{defn}
\label{defn-viscosity}
{(i)} Let $L>0$. We say $u\in\Usub$ (resp., $\Usup$) is a viscosity
$L$-subsolution (resp., $L$-supersolution) of PPDE (\ref{PPDE}) if, for
any $(t,\omega)\in[0, T)\times\O$ and any $\f\in\underline\cA
^{L}u(t,\omega
)$ [resp., $\f\in\overline\cA^{L}u(t,\omega)$],
\[
\bigl\{-\pa_t \f-G^{t,\omega}\bigl(\cdot, \f,\pa_\omega
\f,\pa ^2_{\omega\omega}\f\bigr) \bigr\}(t,{\mathbf0}) \le (
\mbox{resp.,} \ge)\ 0.
\]

 {(ii)} We say $u\in\Usub$ (resp., $\Usup$) is a viscosity
subsolution (resp., supersolution) of PPDE (\ref{PPDE}) if $u$ is
viscosity $L$-subsolution (resp., $L$-supersolution) of PPDE~(\ref
{PPDE}) for some $L>0$.

 {(iii)} We say $u\in\mathrm{ UC}_b(\L)$ is a viscosity
solution of
PPDE (\ref{PPDE}) if it is both a viscosity subsolution and a viscosity
supersolution.
\end{defn}
As pointed out in \cite{ETZ1} Remark~3.11(i), without loss of
generality in (\ref{cA}), we may always set $\textsc{h}= \textsc
{h}^t_\e$ for some
small $\e>0$,
%
%e2.16 #&#
\begin{equation}
\label{chet} \textsc{h}^t_\e:= \inf\bigl\{s> t\dvtx
\bigl|B^t_s\bigr|\ge\e\bigr\} \wedge(t+\e).
\end{equation}

%s3 #&#
\section{Assumptions}
\label{sect:assumptions}

This section collects all of our assumptions on the nonlinearity $G$,
the terminal condition $\xi$ and the underlying path-frozen PDE.

%s3.1 #&#
\subsection{Assumptions on the nonlinearity and terminal conditions}

We first need the conditions on the nonlinearity $G$ as assumed in
\cite{ETZ1}.

%as3.1 #&#
\begin{assum}\label{assum-G}
The nonlinearity $G$ satisfies:
\begin{longlist}[(iii)]
\item[{(i)}] for fixed $(y,z,\g)$, $G(\cdot, y,z,\g)\in\dbL^0(\L)$
%%is $
%\dbF$-progressively measurable,
and $|G(\cdot, 0, {\mathbf0}, {\mathbf0})|\le C_0$;

\item[{(ii)}] $G$ is uniformly Lipschitz continuous in $(y,z,\g)$, with a
Lipschitz constant~$L_0$;

\item[{(iii)}] for any $(y, z,\g)$, $G(\cdot, y,z,\g)$ is right continuous
in $(t,\omega)$ under ${\mathbf{d}}_\infty$;

\item[{(iv)}] $G$ is elliptic, that is, nondecreasing in $\g$. % with a
%modulus of continuity function $\rho_0$.
\end{longlist}
\end{assum}

Our main well-posedness result requires the following strengthening of
(iii) and (iv) above:

%\begin{assum}\label{assum-Guniform}
%$G$ is uniformly continuous in $(t,\o)$ under $\dbf_\infty$ with a
%modulus of continuity function $\rho_0$.
%\end{assum}

%as3.2 #&#
\begin{assum}
\label{assum-Guniform}
(i) $G$ is uniformly continuous in $(t,\omega)$ under ${\mathbf
{d}}_\infty$ with a
modulus of continuity function $\rho_0$.

(ii) For each $\omega$, $G$ is uniformly elliptic. That is, there
exits a
constant $c_0 >0$ such that $G(\cdot,\g_1) - G(\cdot,\g_2) \ge
c_0\tr(\g
_1-\g_2)$ for any $\g_1 \ge\g_2$.
\end{assum}

Condition (i) is needed for our uniform approximation of $G$ below; in
particular it is used (only) in the proof of Lemma~\ref{lem-phie}. We
should point out though, for the semi-linear PPDE and the
path-dependent HJB considered in \cite{ETZ1}, Section~4, this condition
is violated when $\si$ depends on $(t,\omega)$. However, this is a
technical condition due to our current approach for uniqueness.
Condition (ii) is used to ensure the existence of the viscosity
solution for the path-frozen PDE (\ref{PDEe}) below. See also Example~\ref{eg-rep}.

Our first condition on the terminal condition $\xi$ is the following:

%as3.3 #&#
\begin{assum}
\label{assum-xi0}
$\xi\in\dbL^0(\cF_T)$ is bounded and uniformly continuous in
$\omega$
under $\|\cdot\|_T$, with the same modulus of continuity function
$\rho
_0$ as in Assumption~\ref{assum-Guniform}(i).
\end{assum}

%re3.4 #&#
\begin{rem}
\label{rem-xi}
The continuity of a random variable in terms of $\omega$ seems less
natural in stochastic analysis literature. However, since by nature we
are in the weak formulation setting, such continuity is in fact natural
in many applications. This is emphasized in the two following examples:

\begin{itemize}
\item[$-$] Let $V_0:=\dbE^{\dbP_0}[g(X^\si_\cdot)]$, for some bounded function
$g\dvtx \O\longrightarrow\dbR$, and some bounded progressively measurable
process $\sigma$, with
\[
dX^\sigma_t = \si_t \,dB_t,\qquad
\dbP_0\mbox{-a.s.}
\]
In the weak formulation, we define $\dbP^\sigma$ as the probability
measure induced by the process $X^\sigma$, and we re-write $V_0:=\dbE
^{\dbP^\sigma}[g(B_\cdot)]$. Thus the uniform continuity requirement
reduces to that of the function $g$.

\item[$-$] Similarly, the stochastic control problem in strong formulation
$V_0:= \sup_{\underline\si\le\sigma\le\overline\si}\dbE^{\dbP
_0}[g(X^\si_\cdot)]$ for some constants $0\le\underline\si\le
\overline
\si$, may be expressed in the weak formulation as $V_0=\sup_{\underline
\si\le\sigma\le\overline\si}\dbE^{\dbP^\sigma}[g(B_\cdot)]$, thus
reducing the uniform continuity requirement of the terminal data to
that of the function $g$.
\end{itemize}
\end{rem}

Our next assumption is a purely technical condition needed in our proof
of uniqueness. To be precise, it will be used only in the proof of
Lemma~\ref{lem-the} below to ensure the function $\th^\e_n$ constructed
there is continuous in its parameter $\pi_n$. When we have a
representation for the viscosity solution, for example, in the
semilinear case in \cite{ETZ1}, Section~7, we may construct the $\th
^\e
_n$ explicitly and thus avoid the following assumption:

For all $\e>0$, $n\ge0$ and $0\le T_0 < T_1\le T$, denote
%
%e3.1 #&#
\begin{eqnarray}
\label{One2} &&\Pi^\e_n(T_0,
T_1) := \bigl\{\pi_n = (t_i,
x_i)_{1\le i\le n}\dvtx T_0<t_1<
\cds<t_n<T_1,
\nonumber
\\[-8pt]
\\[-8pt]
\nonumber
&&\hspace*{155pt}|x_i|\le\e\mbox{ for all } 1\le i
\le n \bigr\}.
\end{eqnarray}
For all $\pi_n\in\Pi^\e_n(T_0, T_1) $, we denote by $\omega^{\pi
_n}\in\O
^{T_0}$ the linear interpolation of $(T_0, {\mathbf0})$, $(t_i, \sum_{j=1}^i x_j)_{1\le i\le n}$, and $(T, \sum_{j=1}^n x_j)$.

%as3.5 #&#
\begin{assum}
\label{assum-xi2}
There exist $0=T_0<\cds<T_N = T$ such that for each $i = 0,\ldots, N-1$,
for any $\e$ small, any $n$ and any $\omega\in\O$, $\tilde\omega
\in\O
^{T_{i+1}}$, the functions $\pi_n \mapsto\xi(\omega\otimes
_{T_i}\omega^{\pi
_n}\otimes_{T_{i+1}} \tilde\omega)$ and $\pi_n \mapsto G(t,
\omega\otimes
_{T_i}\omega^{\pi_n}\otimes_{T_{i+1}} \tilde\omega, y, z, \g)$
are uniformly
continuous in $\Pi^\e_n(T_i, T_{i+1})$, uniformly on $t\ge T_{i+1}$,
$(y,z,\g)\in\dbR\times\dbR^d\times\dbS^d$ and $\tilde\omega
\in\O^{T_{i+1}}$.
\end{assum}

We note that the uniform continuity of $\xi$ and $G$ implies that the
above mappings are continuous in $\pi_n\in\Pi^\e_n(T_i, T_{i+1})$, but
not necessarily uniformly continuous. In particular, they may not have
limits on the boundary of $\Pi^\e_n(T_i, T_{i+1})$, namely when $t_i =
t_{i+1}$ but $x_i \neq x_{i+1}$.
We conclude this subsection with a sufficient condition for Assumption~\ref{assum-xi2}, where for $\omega\in\O$, we use the notation
$\overline{\omega
}_t:=\max_{s\le t}\omega_s$ and $\underline{\omega}_t:=\min_{s\le t}\omega_s$,
defined componentwise.

%le3.6 #&#
\begin{lem}
\label{lem-xi}
Let $\xi(\omega) = g(\omega_{T_1}, \ldots, \omega_{T_N},
\overline\omega_{T_1},\ldots,
\overline\omega_{T_N}, \underline\omega_{T_1},\ldots, \underline
\omega_{T_N}, \omega)$
for some $0=T_0<T_1<\cds<T_N=T$, and some bounded uniformly continuous
function $(\theta,\omega)\in\dbR^{3dN}\times\O\longmapsto
g(\theta
,\omega)\in\dbR$. Assume further that for all $\theta$, $i$ and
$\omega\in
\O$, there exists a modulus of continuity function $\rho$ and $p>0$
(which may depend on the above parameters), such that
\begin{eqnarray*}
\bigl| g \bigl(\theta, \omega\otimes_{T_i} \omega^1 \otimes
_{T_{i+1}} \tilde\omega \bigr) - g \bigl(\theta, \omega\otimes_{T_i}
\omega^2 \otimes _{T_{i+1}} \tilde\omega \bigr) \bigr| \le\rho
\biggl(\int_{T_i}^{T_{i+1}} \bigl|\omega^1_t-
\omega ^2_t\bigr |^p\,dt \biggr),
\end{eqnarray*}
for all $\omega^1, \omega^2\in\O^{T_i}, \tilde\omega\in\O^{T_{i+1}}$.
Then $\xi$ satisfies Assumptions \ref{assum-xi0} and \ref{assum-xi2}.
\end{lem}
\begin{pf}
Clearly $\xi$ satisfies Assumption~\ref
{assum-xi0}. For
$\omega\in\O
$, $\tilde\omega\in\O^{T_{i+1}}$ and $\pi_n, \tilde\pi_n\in\Pi
^\e
_n(T_i, T_{i+1})$, denote $\hat\omega^{\pi_n} := \omega\otimes
_{T_i}\omega^{\pi
_n}\otimes_{T_{i+1}}\tilde\omega$ and $\hat\omega^{\tilde\pi
_n} := \omega\otimes
_{T_i}\omega^{\tilde\pi_n}\otimes_{T_{i+1}}\tilde\omega$. Then
\[
\bigl|\xi\bigl(\hat\omega^{\tilde\pi_n} \bigr) - \xi\bigl(\hat\omega^{\tilde
\pi_n}
\bigr) \bigr| \le \rho_0 \Biggl(\sum_{k=1}^n
|\theta_k -\tilde\theta_k| \Biggr) + \rho \biggl(\int
_{T_i}^{T_{i+1}} \bigl|\omega^{\pi_n}_t -
\omega ^{\tilde\pi
_n}_t\bigr |^p\,dt \biggr),
\]
where $\rho_0$ is the modulus of continuity function of $g$ with
respect to $(\theta,\omega)$.
Then one can easily check that the $\pi_n\longmapsto\xi(\hat\omega
^{\pi
_n})$ is uniformly continuous in $\Pi^\eps_n(T_i, T_{i+1})$.
\end{pf}

%s3.2 #&#
\subsection{Path-frozen PDEs}
Our main tool for proving the comparison principle for viscosity
solutions, or, more precisely, for constructing the $\ol u^\e$ and
$\ul
u^\e$, mentioned in the \hyperref[sect-Introduction]{Introduction}, so as to prove (\ref
{Perron-Introduction}), is some path-frozen PDE.
Define the following deterministic function on $[t,\infty)\times\dbR
\times\dbR^d\times\dbS^d$:
\[
g^{t,\omega}(s,y,z,\g) := G(s\wedge T, \omega_{\cdot\wedge t}, y,z,\g),\qquad (t,
\omega)\in\L.
\]
For any $\e>0$ and $\eta\ge0$, we denote $T_\eta:=(1+\eta)T$, $\e
_\eta:= (1+\eta)\e$ and
%
%e3.2 #&#
\begin{eqnarray}
\label{Oet}  O_\e&:=&\bigl\{x
\in\dbR^d\dvtx |x|<\e\bigr\},\qquad \overline O_\e:=\bigl\{x\in
\dbR^d\dvtx |x|\le\e \bigr\},\nonumber\\
 \pa O_\e&:=&\bigl\{x\in
\dbR^d\dvtx |x|=\e\bigr\},
\nonumber
\\[-8pt]
\\[-8pt]
\nonumber
Q^{\e,\eta}_t&:= &[t,T_\eta)\times O_{\e_\eta},\qquad
\overline Q^{\e,\eta}_t := [t, T_\eta]\times\overline
O_{\e_\eta},\\
 \pa Q^{\e,\eta}_t &:=& \bigl([t,T_\eta]
\times\partial O_{\e_\eta} \bigr) \cup \bigl(\{T_\eta\}\times
O_{\e_\eta} \bigr), \nonumber
\end{eqnarray}
and we further simplify the notation for $\eta=0$ as
\[
Q^{\e}_t := Q^{\e,0}_t, \qquad \overline
Q^{\e}_t := \overline Q^{\e,0}_t,\qquad  \pa
Q^{\e}_t := \pa Q^{\e,0}_t.
\]
Our additional assumption is formulated on the following localized and
path-frozen PDE defined for every $(t,\omega)\in\L$:
%
%e3.3 #&#
\begin{equation}
\label{PDEe} \mbox{(E)}^{t,\omega}_{\eps,\eta} \qquad \mathbf{L}^{t,\omega}v
:= -\pa_t v - g^{t,\omega} \bigl(s, v, D v, D^2 v
\bigr) = 0 \qquad\mbox{on } Q^{\e,\eta}_t.
\end{equation}
Notice that for fixed $(t,\omega)$, this is a standard deterministic
partial differential equation.

%le3.7 #&#
\begin{lem}
\label{lem-PDEcomparison}
Under Assumptions \ref{assum-G} and \ref{assum-Guniform}\textup{(ii)}, PDE
(\ref{PDEe}) satisfies the comparison principle for bounded viscosity
solutions (in standard sense, as in \cite{CIL}). Moreover, for any
$h\in C^0(\pa Q^{\e,\eta}_t)$, PDE (\ref{PDEe}) with the boundary
condition $h$ has a (unique) bounded viscosity solution $v$.
\end{lem}
\begin{pf}
The comparison principle follows from standard
theory; see, for
example, \cite{CIL}. Moreover, as we will see later, the $\overline v$
and $\underline v$ defined in (\ref{rep-boundingPDE}) are viscosity
supersolution and subsolution, respectively, of the PDE (\ref{PDEe}) and
satisfy $\overline v = \underline v = h$ on $\pa Q^{\e,\eta}_t$. Then
the existence follows from the standard Perron approach in the spirit
of \cite{CIL}, Theorem~4.1.
\end{pf}

We will use the following additional assumption:

%\begin{assum}\label{assum-existence} One of the following conditions
%holds:

%(i) $G$ is convex or concave in $ \g$;
% $ g^{t,\o}_\d(\cdot,\g):=\inf_{A\in\dbS^d_+}\big\{ g^{t,\o}(\cdot,
%\g+A)- \d
%I_d:A\big\}>-\infty$ for $0\le\d\le c_0$, for some $c_0>0$, and $g^{t,
%\o}_\d\longrightarrow g$ as $\d\searrow0$,

%(ii) For any $\e>0$, $\eta\ge0$, $(t,\o)\in\L$, and any $h\in C^0(
%\pa Q^{\e,\eta}_t)$, the PDE \reff{PDEe} with boundary condition $h$
%has a (unique) viscosity solution $v$.
%\end{assum}
%
%as3.8 #&#
\begin{assum}\label{assum-comparison}
For any $\e>0,\eta\ge0$, $(t,\omega)\in\L$ and $h\in C^0(\pa
Q^{\e,\eta
}_t)$, we have $\overline v = v=\underline v$, where $v$ is the unique
viscosity solution of PDE (\ref{PDEe}) with boundary condition $h$, and
%
%e3.4 #&#
\begin{eqnarray}
\label{barv}
\overline v(s,x) &:=& \inf
\bigl\{w(s,x)\dvtx w \mbox{ classical supersolution of $\mathrm{(E)}$}^{t,\omega}_{\eps,\eta}
\nonumber\\
&&\hspace*{119pt}\mbox{and } w \ge h \mbox{ on } \pa Q^{\e,\eta}_t \bigr\},
\nonumber
\\[-8pt]
\\[-8pt]
\nonumber
\underline v(s,x) &:=& \sup \bigl\{w(s,x)\dvtx w \mbox{ classical subsolution of
$\mathrm{(E)}$}^{t,\omega}_{\eps,\eta} \\
&&\hspace*{111pt}\mbox{and } w \le h \mbox{ on } \pa
Q^{\e,\eta}_t \bigr\}.\nonumber
\end{eqnarray}
\end{assum}

We first note that the above sets of $w$ are not empty. Indeed, one can
check straightforwardly that for any $\delta>0$ and denoting
$\lambda_\delta:= {C_0
+ L_0\|h\|_\infty\over\delta} + L_0$,
%
%e3.5 #&#
\begin{equation}
\label{barw} \overline w (t, x) := \|h\|_\infty+ \delta
e^{\lambda_\delta
(T_\eta-t)},\qquad  \underline w(t,x) := -\|h\|_\infty- \delta
e^{\lambda_\delta
(T_\eta-t)}
\end{equation}
satisfy the requirement for $\overline v(s,x)$ and $\underline v(s,x)$,
respectively.
We also observe that our definition (\ref{barv}) of $\overline v$ and
$\underline v$ is different from the corresponding definition in the
standard Perron approach \cite{Ishii}, in which the $w$ is a viscosity
supersolution or subsolution. It is also different from the recent
development of Bayraktar and Sirbu~\cite{BS}, in which the $w$ is a so
called stochastic supersolution or subsolution. Loosely speaking,
%%assuming that the comparison result holds for the equation (E)$^{t,
%\o}_{\e,\eta}$,
our Assumption~\ref{assum-comparison} requires that the viscosity
solution of (E)$^{t,\omega}_{\e,\eta}$ can be approximated by a
sequence of
classical supersolutions and a sequence of classical subsolutions. We
shall discuss further this issue in Section~\ref{sect-Perron} below. In
particular, we will provide some sufficient conditions for Assumption~\ref{assum-comparison} to hold.

\section{Main results}
\label{sect-wellposedness}
%\setcounter{equation}{0}

%Let $\mbox\mathrm{ UC}_b(\cF_T,\dbR)$ denote the collection of all scalar
%bounded $\cF_T$-measurable r.v.
The following theorem is the main result of this paper:

%th4.1 #&#
\begin{thmm}\label{thmm-wellposedness}
Let Assumptions \ref{assum-G}, \ref{assum-Guniform}, \ref{assum-xi0},
\ref{assum-xi2} and \ref{assum-comparison} hold true:

\begin{longlist}[(ii)]
\item[{(i)}] Let $u^1\in\Usub$ be a viscosity subsolution and $u^2\in
\Usup$ a viscosity supersolution of PPDE (\ref{PPDE}) with
$u^1(T,\cdot)
\le\xi\le u^2(T,\cdot)$. Then $u^1\le u^2$ on $\L$.

\item[{(ii)}] PPDE (\ref{PPDE}) with terminal condition $\xi$ has a unique
viscosity solution $u\in \mathrm{UC}_b(\L)$.
\end{longlist}
\end{thmm}

%s4.1 #&#
\subsection{Strategy of the proof}
There are two key ingredients for the proof of this main result. The
first is the following partial comparison, proved in Section~\ref
{sect-partialcomparison}, which extends the corresponding result in
Proposition~5.3 of \cite{ETZ1} to the set $\overline{C}^{1,2}(\L)$. The
reason for extending $C^{1,2}(\L)$ to $\ol C^{1,2}(\L)$ is that
typically we can construct the approximations $\ol u^\e$ and $\ul u^\e
$, mentioned in the \hyperref[sect-Introduction]{Introduction}, only in the space $\ol C^{1,2}(\L)$,
and not in $C^{1,2}(\L)$.

%pr4.2 #&#
\begin{prop}
\label{prop-comparison}
Assume Assumption~\ref{assum-G} holds true. Let $ u^2\in\Usup$ be a
viscosity supersolution of PPDE (\ref{PPDE}) and $ u^1\in\overline
C^{1,2}(\L)$ bounded from above satisfying $\cL u^1(t,\omega) \le0$ for
all $(t,\omega) \in\L$ with $t<T$. If $ u^1(T,\cdot) \le
u^2(T,\cdot)$, then $
u^1\le u^2$ on $\L$.

A similar result holds if we switch the roles of $u^1$ and $u^2$.
\end{prop}

The second key ingredient follows the spirit of the Perron approach as
in \cite{EKTZ}. Let
%
%e4.1 #&#
\begin{eqnarray}
\label{baru} \overline u(t,\omega) &:= &\inf \bigl\{\psi_t\dvtx \psi
\in\overline\cD_T^\xi(t,\omega) \bigr\},
\nonumber
\\[-8pt]
\\[-8pt]
\nonumber
 \underline u(t,
\omega)& := &\sup \bigl\{\psi_t\dvtx \psi\in\underline
\cD_T^\xi(t,\omega) \bigr\},
\end{eqnarray}
where
%
%e4.2 #&#
\begin{eqnarray}
\label{cD}  \overline
\cD_T^\xi(t,\omega)& := &\bigl\{ \psi\in\overline
C^{1,2}\bigl(\L^t\bigr)\dvtx \psi^- \mbox{ bounded},\nonumber\\
&&\hspace*{6pt} (\cL
\psi)^{t,\omega} \ge0 \mbox{ on } [t,T)\times\O^t,
\psi_T \ge\xi^{t,\omega} \bigr\},
\nonumber
\\[-8pt]
\\[-8pt]
\nonumber
\underline\cD_T^\xi(t,\omega) &:=& \bigl\{\psi\in
\overline C^{1,2}\bigl(\L^t\bigr)\dvtx \psi^+ \mbox{ bounded},\\
&&\hspace*{6pt}(\cL\psi)^{t,\omega} \le0 \mbox{ on } [t,T)\times\O^t,
\psi_T \le\xi ^{t,\omega} \bigr\}.\nonumber
\end{eqnarray}
By using the functional It\^o formula (\ref{Ito2}), and following the
arguments in \cite{ETZ1}, Theorem~3.16, we obtain a similar result as
the partial comparison of Proposition~\ref{prop-comparison}, implying that
%
%e4.3 #&#
\begin{equation}
\label{u<u} \underline u \le \overline u.
\end{equation}
Moreover, these processes satisfy naturally a partial dynamic
programming principle which implies the following viscosity properties.

%pr4.3 #&#
\begin{prop}\label{prop-viscosity}
Let Assumptions \ref{assum-G}, \ref{assum-Guniform} and \ref{assum-xi0}
hold true.
Then the processes $\overline u$ and $ \underline u$ are bounded,
uniformly continuous viscosity super solutions and subsolutions,
respectively, of PPDE (\ref{PPDE}).
\end{prop}

This result will be proved in Section~\ref{sect:viscosity}. A crucial
step for our proof is to show the consistency of the Perron approach in
the sense that equality holds in the last inequality, under our
additional assumptions.

%pr4.4 #&#
\begin{prop}\label{prop-perron}
Under the conditions of Theorem~\ref{thmm-wellposedness}, with $N=1$ in
Assumption~\ref{assum-xi2}, we have $\overline u = \underline u$.
\end{prop}

The proof of this proposition is reported in Section~\ref{sect-baru}.
Given Propositions \ref{prop-comparison}, \ref{prop-viscosity} and
\ref
{prop-perron}, Theorem~\ref{thmm-wellposedness} follows immediately.

\begin{pf*}{Proof of Theorem~\ref{thmm-wellposedness}} We
prove the
theorem in three steps:

\textit{Step} 1. We first consider the case $N=1$ in Assumption~\ref
{assum-xi2}. By Proposition~\ref{prop-comparison}, we have $u^1 \le
\overline u$ and $\underline u \le u^2$. Then Proposition~\ref
{prop-perron} implies $u^1\le u^2$ immediately, which implies (i) and
the uniqueness of the viscosity solution. Finally, by Propositions \ref
{prop-perron} and \ref{prop-viscosity}, $u:=\overline{u}=\underline{u}$
is a viscosity solution of (\ref{PPDE}).

\textit{Step} 2. For general $N$, it follows from step 1 that the
comparison, existence and uniqueness of the viscosity solution holds on
$[T_{N-1}, T_N]$. Let $u$ denote the unique viscosity solution on
$[T_{N-1}, T_N]$ with terminal condition $\xi$, constructed by the
Perron approach. Now consider PPDE (\ref{PPDE}) on $[T_{N-2}, T_{N-1}]$
with terminal condition $u(T_{N-1}, \cdot)$. We shall prove below that
$u(T_{N-1}, \cdot)$ satisfies the requirement of step 1. Then we may
extend the comparison, existence and uniqueness of the viscosity
solution to the interval $[T_{N-2}, T_N]$. By repeating the arguments
backwardly, we complete the proof of Theorem~\ref{thmm-wellposedness}.

\textit{Step} 3. It remains to verify Assumptions \ref{assum-xi0} and
\ref
{assum-xi2} with $N=1$ for $u(T_{N-1}, \cdot)$ on $[T_{N-2}, T_{N-1}]$.
First, by Proposition~\ref{prop-viscosity} it is clear that $u(T_{N-1},
\cdot)$ is bounded. Given $\omega\in\O$, note that PPDE (\ref
{PPDE}) on
$[T_{N-1}, T_N]$ can be viewed as a PPDE with generator $G^{T_{N-1},
\omega
}$ and terminal condition $\xi^{T_{N-1},\omega}$. Then, following the
arguments in Lemma~\ref{lem-uUC}(i) below, one can easily show that
$u(T_{N-1}, \omega)$ is uniformly continuous in $\omega$, and it
follows from
Assumption~\ref{assum-xi2} that $u(T_{N-1}, \omega\otimes_{T_{N-2}}
\omega^{\pi
_n})$ is uniformly continuous in $\pi_n\in\break \Pi^\e_n(T_{N-2}, T_{N-1})$.
\end{pf*}

%s4.2 #&#
\subsection{Heuristic analysis on Proposition \texorpdfstring{\protect\ref{prop-perron}}{4.4}}

\label{sect-Heuristic}
While highly technical, Proposition~\ref{prop-comparison} follows along
the same lines as the partial comparison of \cite{ETZ1}, Proposition~5.3. Proposition~\ref{prop-viscosity} has a corresponding result in PDE
literature, and is proved in the spirit of the stability result of
\cite
{ETZ1}, Theorem~5.1. In this subsection, we provide some heuristic
discussions on Proposition~\ref{prop-perron}, focusing on the case
$\ol
u_0 = \ul u_0$, and the rigorous arguments will be carried out in
Section~\ref{sect-baru} below.

We shall follow \cite{ETZ1}, Section~7, where Proposition~\ref
{prop-perron} is proved in a much simpler, semi-linear setting. The
idea is to construct $\ol u^\e\in\overline\cD_T^\xi(0, 0)$ and
$\ul
u^\e\in\underline\cD_T^\xi(0, 0)$ such that $\lim_{\e\to0}[\ol
u^\e
_0 - \ul u^\e_0] = 0$. To be precise, modulus some technical
properties, the approximations $\ol u^\e, \ul u^\e$ should satisfy:

\begin{itemize}
\item they are piecewise smooth and $\cL\ol u^\e\ge0\ge\cL\ul
u^\e$;

\item they are continuous in $t$;

\item $\ol u^\e_T$ and $\ul u^\e_T$ are close to $\xi$.
\end{itemize}

To achieve this, we shall discretize the path $\omega$ so that we can
utilize the path-frozen PDE (\ref{PDEe}). We note that such
discretization of $\omega$ will not induce big errors, thanks to the
uniform continuity of the involved processes. Fix $\e>0$, and set
$\textsc{h}
_0:= 0$,
\[
\textsc{h}_{i+1} := \bigl\{t\ge\textsc{h}_i\dvtx
|B_t -B_{t_i}| = \e\bigr\} \wedge T.
\]
Denote $\hat\pi_n := \{(\textsc{h}_i, B_{\textscc{h}_i}), 0\le i\le
n\}$. Let $\pi_n =
\{(t_i, x_i), 0\le i\le n\}$ be a typical value of $\hat\pi_n(\omega
)$, and
$\omega^{\pi_n}\in\O$ be the linear interpolation of $\pi_n$. The main
idea is to construct a sequence of deterministic functions $v^\e_n(\pi
_n; t, x)$ so that we may construct the desired $\ol u^\e$ and $\ul
u^\e
$ from a common process $u^\e_t:= v^\e_n(\hat\pi_n; t, B_t -
B_{\textscc{h}
_n})$, $\textsc{h}_n \le t <\textsc{h}_{n+1}$. For this purpose, we
require $v^\e
_n$, and hence $u^\e$, satisfying the following three corresponding properties:

$\bullet$ For each $\pi_n$, the function $v^\e_n(\pi_n; \cdot)$
is in
$C^{1,2}(Q^\e_{t_n})$ and is a classical solution of a certain
mollified path-frozen PDE,
%
%e4.4 #&#
\begin{equation}
\label{PDEven} -\pa_t v^\e_n -
g^{\pi_n}_\e\bigl(t, v^\e_n, D
v^\e_n, D^2 v^\e_n
\bigr) = 0,
\end{equation}
where $g^{\pi_n}_\eps= g^{t_n, \omega^{\pi_n}}_\eps$.
Consequently, the
process $u^\e$ is approximately a classical solution of PPDE (\ref
{PPDE}) on $[\textsc{h}_n, \textsc{h}_{n+1}]$, thanks to the fact that
$g^{\hat\pi_n(\omega
)}_\e(t,\cdot)$ is a good approximation of $G(t,\omega,\cdot)$.

$\bullet$ $v^\e_n(\hat\pi_n; \textsc{h}_{n+1}, B_{\textscc
{h}_{n+1}}- B_{\textscc{h}_n}) = v^\e
_{n+1}(\hat\pi_{n+1}; \textsc{h}_{n+1}, {\mathbf0})$ so that $u^\e$ is
continuous
in $t$ and is more or less in $\ol C^{1,2}(\L)$.

$\bullet$ $v^\e_n(\pi_n; T, x)$ is constructed from $\xi$, so that
$u^\e
_T$ is close to $\xi$.

 Now by the uniform continuity of $\xi$ and $G$, we will see that
$\ol u^\e:= u^\e+ \rho_0(2\e)$ and $\ul u^\e:= u^\e- \rho_0(2\e)$
satisfy the desired classical semi-solution property. Clearly $\ol u^\e
- \ul u^\e\le2\rho_0(2\e)$, implying the result.

In \cite{ETZ1}, Section~7, the functions $v^\e_n$ can be constructed
explicitly via approximating backward SDEs. In the present setting,
since we do not have a representation for the candidate solution, we
cannot construct $v^\e_n$ directly. By some limiting procedure, in
Lemma~\ref{lem-the} below, we shall construct certain deterministic
functions $\th^\e_n$ which satisfy all the above three properties,
except that $\th^\e_n$ is only a viscosity solution of PDE (\ref
{PDEven}). Now to construct smooth $v^\e_n$ from $\th^\e_n$, we apply
Assumption~\ref{assum-comparison}. In fact, given the viscosity
solution $\th^\e_n$, Assumption~\ref{assum-comparison} allows us to
construct the classical supersolution $\ol v^\e_n$ and the classical
subsolution $\ul v^\e_n$, rather than one single smooth function $v^\e
_n$, such that $\ul v^\e_n \le\th^\e_n \le\ol v^\e_n$, and $\ol
v^\e
_n - \ul v^\e_n$ is small. This procedure is carried out in Lemma~\ref
{lem-phie} below, and the construction is done piece by piece,
forwardly on each random interval $[\textsc{h}_n, \textsc{h}_{n+1}]$.

%re4.5 #&#
\begin{rem}
\label{rem-derivative}
As we see in the above discussion, the processes we will use to
prove the comparison takes the form $v(\Pi_n; t, B_t - B_{\textscc
{h}_n})$, $\textsc{h}
_n \le t<\textsc{h}_{n+1}$, for some deterministic function $v$, which is
smooth in $(t,x)$. Then it suffices to apply the standard It\^{o}
formula on $v$, rather than the functional It\^{o} formula. Indeed,
under our assumptions, we can prove rigorously the well-posedness of
viscosity solutions, including existence, stability and comparison and
uniqueness, without using the functional It\^{o} formula. In other
words, technically speaking, we can establish our theory without
involving the path derivatives. However, we do feel that the path
derivatives and the functional It\^{o} formula are the natural and
convenient language in this path-dependent framework. In particular, it
is much more natural to talk about classical solutions of PPDEs by
using the path derivatives. Moreover, the current proof relies heavily
on the discretization of the underlying path $\omega$, with the help
of the
path-frozen PDEs. This discretization induces the above piecewise
Markovian structure. The functional It\^{o} formula allows us to
explore in future research other approaches without using such discretization.
\end{rem}

%s4.3 #&#
\subsection{The bounding equations}

The proof of Proposition~\ref{prop-perron} requires some estimates,
which involve the following particular example analyzed in \cite
{ETZ1}. Recall the constants $L_0$ $C_0$, and $c_0$ from Assumptions
\ref{assum-G} and \ref{assum-Guniform}, and consider the operators
%
%e4.5 #&#
\begin{eqnarray}
\label{overlineg} \overline
g_0(z, \g)& :=& \sup_{|\a|\le L_0, \sqrt{2c_0} \le|\b| \le\sqrt{2L_0} } \biggl[ \a\cdot z +
{1\over2} \b^2 \dvtx \g\biggr],\nonumber\\
\overline g(y,z,\g)&:=&
\overline g_0(z,\g)+L_0|y|+C_0,
\nonumber
\\[-8pt]
\\[-8pt]
\nonumber
\underline g_0(z, \g)& := &\inf_{|\a|\le L_0, \sqrt{2c_0} \le|\b|\le\sqrt{2L_0} } \biggl[ \a
\cdot z +{1\over2} \b^2 \dvtx \g\biggr],\\
  \underline
g(y,z,\g)&:=&\underline g_0(z,\g)- L_0 |y|-C_0,\nonumber
%\label{underlineg}
\end{eqnarray}
which clearly satisfy Assumptions \ref{assum-G} and \ref
{assum-Guniform}, and
%
%e4.6 #&#
\begin{equation}
\label{Goverlineg} \underline g \le G\le\overline g.
\end{equation}
These operators induce the PPDEs
%
%e4.7 #&#
\begin{eqnarray}
\label{cLbar} \overline{\cL}u&:=&-\partial_tu-\overline g(u,
\partial_\omega u,\partial_{\omega
\omega}u)=0 \quad\mbox{and}
\nonumber
\\[-8pt]
\\[-8pt]
\nonumber
 \underline{
\cL}u&:=&-\partial_tu-\underline g(u,\partial_\omega u,\partial
_{\omega\omega}u)=0.
\end{eqnarray}
Let $\cB^t_{L_0} :=  \{ b \in\dbL^0(\L^t)\dvtx |b|\le L_0 \}$ and
%
%e4.8 #&#
\begin{eqnarray}
\label{cEc0} \cP^t_{L_0,c_0} &:=& \bigl\{
\cP^t_{L_0}\dvtx \bigl|\b^\dbP\bigr|\ge
\sqrt{2c_0}\bigr\},\qquad \overline\cE^{L_0, c_0}_t := \sup
_{\dbP\in\cP^t_{L_0, c_0}} \dbE ^\dbP,
\nonumber
\\[-8pt]
\\[-8pt]
\nonumber
  \underline\cE^{L_0, c_0}_t
&:=& \inf_{\dbP\in\cP^t_{L_0,c_0}} \dbE ^\dbP.
\end{eqnarray}
Following the arguments in our accompanying paper (\cite{ETZ1},
Proposition~4), we see that for a bounded, uniformly continuous $\cF
_T$-measurable r.v. $\xi$,
%
%e4.9 #&#
\begin{eqnarray}
\label{rep-boundingPPDE} %
\overline{w}(t,\omega) &:= &\sup
_{b\in\cB^t_{L_0}} \overline\cE^{L_0, c_0}_{t} \biggl[
\xi^{t,\omega}e^{\int_t^T b_r\,dr} + C_0 \int_t^T
e^{\int_t^s b_r\,dr}\,ds \biggr],
\nonumber
\\[-8pt]
\\[-8pt]
\nonumber
\underline w(t,\omega) &:=& \inf_{b\in\cB^t_{L_0}} \underline
\cE^{L_0,c_0}_{t} \biggl[\xi^{t,\omega}e^{\int_t^T b_r\,dr} -
C_0 \int_t^T e^{\int_t^s b_r\,dr}\,ds
\biggr]
\end{eqnarray}
are viscosity solutions of the PPDE $\overline{\cL}\overline{w}:=0$ and
$\underline{\cL}\underline{w}:=0$, respectively.

By Lemma~\ref{lem-PDEcomparison}, the PDE version of (\ref{cLbar}),
%
%e4.10 #&#
\begin{eqnarray}
\label{Lbar-PDE} \overline{\mathbf{L}}v&:=&-\partial_tv-\overline g
\bigl(v,Dv,D^2v\bigr)=0 \quad \mbox{and}
\nonumber
\\[-8pt]
\\[-8pt]
\nonumber
 \underline{\mathbf{L}}v&:=&-
\partial_tv-\underline g\bigl(v,Dv,D^2v\bigr)=0 \qquad\mbox {in }
Q^{\e,\eta}_t,
\end{eqnarray}
satisfies the comparison principle. Moreover, we have the following:
%
%le4.6 #&#
\begin{lem}
\label{lem-rep-boundingPDE}
Under Assumptions \ref{assum-G} and \ref{assum-Guniform}\textup{(ii)}, for any
$h\in C^0 (\partial Q^{\e,\eta}_t )$, the following functions are
the unique viscosity solutions of PDEs (\ref{Lbar-PDE}) with boundary
condition $h$:
%
%e4.11 #&#
\begin{eqnarray}
\label{rep-boundingPDE} %
 \overline{v}(t, x) &:=& \sup
_{b\in\cB^t_{L_0}} \overline\cE^{L_0,c_0}_{t}
\biggl[e^{\int_t^{\textscc{h}} b_r\,dr}h\bigl(\textsc{h},x+B^t_{\textscc{h}}
\bigr) + C_0 \int_t^{\textscc{h}}
e^{\int_t^s b_r\,dr}\,ds \biggr],
\nonumber
\\[-8pt]
\\[-8pt]
\nonumber
\underline{v}(t, x) &:=& \inf_{b\in\cB^t_{L_0}} \underline
\cE^{L_0,c_0}_{t} \biggl[e^{\int_t^{\textscc{h}} b_r\,dr}h\bigl(
\textsc{h},x+B^t_{\textscc{h}}\bigr) - C_0 \int
_t^{\textscc{h}} e^{\int_t^s b_r\,dr}\,ds \biggr],
\end{eqnarray}
where $\textsc{h}:=\textsc{h}^{t,x}:=\{s>t\dvtx (s,x+B^t_s)\notin Q^{\e
,\eta}_t\}$.
\end{lem}
\begin{pf}First, by the arguments in \cite{ETZ0}, one may
easily check
that $\overline v$ and $\underline v$ are continuous and satisfy
dynamic programing principle for $t<\textsc{h}$, which implies the viscosity
property immediately. Then it remains to check the boundary conditions.
For $x\in\overline O_{\e_\eta}$, since $t\le H^{t,x} \le T$ and $h$ is
uniformly continuous with certain modulus of continuity function $\rho
_h$, it is clear that
%
%e4.12 #&#
\begin{eqnarray}
\label{boundaryest}
&&\bigl|\overline{v}(t, x) - h(T,x)\bigr|\nonumber\\
&&\qquad\le\sup_{b\in\cB^t_{L_0}}
\overline\cE^{L_0}_{t} \biggl[ \biggl|e^{\int_t^{\textscc{h}} b_r\,dr}h\bigl(
\textsc{h},x+B^t_{\textscc{h}}\bigr) + C_0 \int
_t^{\textscc{h}} e^{\int_t^s b_r\,dr}\,ds - h(T,x) \biggr| \biggr]
\nonumber
\\
&&\qquad=\sup_{b\in\cB^t_{L_0}} \overline\cE^{L_0}_{t}
\biggl[ \biggl|\bigl[e^{\int_t^{\textscc{h}}
b_r\,dr}-1\bigr]h\bigl(\textsc{h},x+B^t_{\textscc{h}}
\bigr) + \bigl[h\bigl(\textsc{h},x+B^t_{\textscc{h}}\bigr) -
h(T,x)\bigr]\nonumber\\
&&\hspace*{236pt}{} + C_0 \int_t^{\textscc{h}}
e^{\int_t^s b_r\,dr}\,ds \biggr| \biggr]
\\
&&\qquad\le C\overline\cE^{L_0}_{t} \bigl[ H-t +
\rho_h\bigl(T-\textsc{h}+ \bigl|B^t_{\textscc{h}}\bigr|\bigr)
\bigr] \nonumber\\
&&\qquad\le C\overline\cE^{L_0}_{t} \bigl[ T-t +
\rho_h\bigl(T-t +\bigl \|B^t\bigr\|_T\bigr) \bigr]\nonumber\\
&&\qquad\to0,\nonumber
\end{eqnarray}
as $t\uparrow T$. Furthermore, let $t<T$ and ${\mathbf0}\neq x\in O_{\e
_\eta}$. Note that for any $a>0$ and $\dbP\in\cP^t_{L_0, c_0}$,
\begin{eqnarray*}
\dbP \bigl(\textsc{h}^{t,x} - t \ge a \bigr) &\le& \dbP \biggl(\sup
_{t\le s\le t+a} {x\over|x|} \cdot B^t_s
\le\e_\eta- |x| \biggr)
\\
&\le& \dbP \biggl(\sup_{t\le s\le t+a} \int_t^s
{x\over|x|} \cdot \b^\dbP_r \,d
W^\dbP_r \le\e_\eta- |x| + L_0 a
\biggr).
\end{eqnarray*}
Let $A_s := \int_t^s {x^T\over|x|} (\b^\dbP_r)^2 {x\over|x|}\,dr$ and
$\t_s := \inf\{r \ge t\dvtx A_r \ge s\}$. Then $M_s :=\break  \int_t^{\t_s}
{x^T\over|x|} \b^\dbP_r \,d W^\dbP_r$ is a $\dbP$-Brownian motion, and
$A_s \ge2c_0 (s-t)$. % which implies $\t_s -t \le{s-t \over2c_0}$.
Thus
\begin{eqnarray*}
\dbP \bigl(\textsc{h}^{t,x} - t \ge a \bigr) &\le& \dbP \Bigl(\sup
_{t\le s\le t+2
c_0 a} M_s \le\e_\eta- |x| +
L_0 a \Bigr) \\
&=& \dbP_0\bigl(\|B\|_{ 2c_0 a} \le
\e_\eta- |x| + L_0 a \bigr)
\\
&=&\dbP_0\bigl(|B_{ 2c_0a}| \le\e_\eta- |x| +
L_0 a \bigr) \\
&=& \dbP_0\biggl(|B_1| \le
{1\over\sqrt{2c_0a}}\bigl[\e_\eta- |x| + L_0 a\bigr] \biggr)
\\
&\le& {C\over\sqrt{a}}\bigl[\e_\eta- |x| + L_0 a\bigr].
\end{eqnarray*}
Set $a := \e_\eta- |x| $, and we get
\[
\dbP \bigl(\textsc{h}^{t,x} - t \ge\e_\eta- |x| \bigr) \le C
\sqrt {\e_\eta- |x|}.
\]
Following similar arguments to those in (\ref{boundaryest}), one can
easily show that for some modulus of continuity function $\rho$,
\begin{eqnarray*}
\bigl|\overline{v}(t, x) - h(t, \tilde x)\bigr| \le\rho\bigl( \e_\eta- |x|\bigr)\qquad \mbox
{where } \tilde x:= {|x|\over\e_\eta} x \in\pa O_{\e_\eta}.
\end{eqnarray*}
Then, for $t_0<T$, $x_0 \in\pa O_{\e_\eta}$, $t<T$ and $x\in O_{\e
_\eta
}$, noting that
\[
|x-\tilde x| \le\e_\eta-|x| = |x_0|-|x|
\le|x-x_0|,
\]
we have, as $(t,x) \to(t_0, x_0)$,
\begin{eqnarray*}
\bigl|\overline{v}(t, x) - h(t_0, x_0)\bigr| &\le& \bigl|
\overline{v}(t, x) - h(t, \tilde x)\bigr| + \bigl|h(t, \tilde x) - h(t_0,
x_0)\bigr|
\\
&\le& \rho\bigl(\e_\eta-|x|\bigr) + \rho_h\bigl(|t-t_0| +
|x_0 - \tilde x|\bigr)
\\
&\le& \rho\bigl(|x_0-x|\bigr) + \rho_h\bigl(|t-t_0| +2
|x_0 - x| \bigr) \to0.
\end{eqnarray*}
This implies that $\overline v$ is continuous on $\overline Q^{\e,\eta
}$. Similarly one can prove the result for $\underline v$.
\end{pf}

We remark that (\ref{rep-boundingPPDE}) provides representation for
viscosity solutions of PPDEs~(\ref{cLbar}), even in the degenerate case
$c_0=0$. However, this is not true for the PDEs (\ref{Lbar-PDE}), due to
the boundedness of the domain $Q^{\e,\eta}_t$, which induces the
hitting time $\textsc{h}$ and ruins the required regularity, as we
will see in
next example.
%
%ex4.7 #&#
\begin{eg}
\label{eg-rep}
Assume Assumption~\ref{assum-G} holds, but $G$ is degenerate, and thus
$c_0=0$. Let $d=1$, and set $h(s,x) := s$ on $\pa Q^{\e,\eta}_t$. Then
the $\overline v$ defined by (\ref{rep-boundingPDE}) is discontinuous in
$[0, T_\eta) \times\pa O_{\e_\eta} \subset\pa Q^{\e,\eta}_0$ and thus
is not a viscosity solution of the PDE (\ref{Lbar-PDE}).
\end{eg}
\begin{pf}
It is clear that
\[
\overline v(t,x) = \overline\cE^{L_0}_{t}
\biggl[e^{L_0 (\textsc{h}-
t) } \textsc{h} + C_0 \int_t^{\textscc{h}}
e^{L_0(s-t)}\,ds \biggr],
\]
where the integrand is increasing in $\textsc{h}$ which takes values
on $[t,
T_\eta]$. Then, by taking the $\dbP$ corresponding to $\a= \b= 0$, we
have $\textsc{h}= T_\eta$, $\dbP$-a.s. and thus
\[
\overline v(t,x) = e^{L_0 (T_\eta- t) } T_\eta+ C_0 \int
_t^{T_\eta} e^{L_0(s-t)}\,ds,\qquad (t,x)\in
Q^{\e,\eta}_0.
\]
However, we have $\overline v(t,x) = t$ on $\pa Q^{\e,\eta}_0$, so
$\overline v$ is discontinuous in $[0, T_\eta) \times\pa O_{\e_\eta}$.
\end{pf}

%s4.4 #&#
\subsection{A change of variables formula}
We conclude this section with a change of variables formula, which is
interesting in its own right.
We have previously observed in \cite{ETZ1}, Remark~3.15, that the
classical change of variables formula is not known to hold true for our
notion of viscosity solutions under Assumption~\ref{assum-G}. We now
show that it holds true under the additional Assumption~\ref{assum-comparison}.

Let $u \in C^{1,2}_b(\L)$ and $\Phi\in C^{1,2}([0,T]\times\dbR)$.
Assume $\Phi$ is strictly increasing in~$x$, and let $\Psi$ denote its
inverse function. Note that $\Psi$ is increasing in $x$ and $\Psi_x
>0$. Define
%
%e4.13 #&#
\begin{equation}
\label{change} \tilde u(t,\omega) := \Phi\bigl(t, u(t,\omega)\bigr) \qquad\mbox{and
thus } u(t,\omega) = \Psi\bigl(t, \tilde u(t,\omega)\bigr).
\end{equation}
Then direct calculation shows that
%
%e4.14 #&#
\begin{eqnarray}
\label{PPDE-change} \cL u (t,\omega) &=& \Psi_x \bigl(t,\tilde u(t,\omega)
\bigr)\tilde \cL\tilde u(t,\omega)\quad \mbox{and}
\nonumber
\\[-8pt]
\\[-8pt]
\nonumber
 \tilde\cL\tilde u &:=& -
\pa_t\tilde u - \tilde G \bigl(t,\omega, \tilde u, \pa_\omega
\tilde u, \pa^2_{\omega\omega}\tilde u \bigr),
\end{eqnarray}
where
\begin{eqnarray*}
&&\tilde G(t,\omega, y, z, \g) \\
&&\qquad:= \frac{\Psi_t(t,y)
+G (t,\omega, \Psi(t,y), \Psi_x(t,y) z,
\Psi_{xx}(t,y) z^2 + \Psi_x(t,y) \g )
}{\Psi_x(t, y)}.
\end{eqnarray*}
Then the following result is obvious:

%pr4.8 #&#
\begin{prop}
\label{prop-change}
Under the above assumptions on $\Psi$, $u$ is classical solution (resp.,
supersolution, subsolution) of $\cL u =0$ if and only if $\tilde u:=
\Phi(t,u)$ is a classical solution (resp., supersolution, subsolution)
of $\tilde\cL\tilde u=0$.
\end{prop}

Moreover, we have the following:

%th4.9 #&#
\begin{thmm}
\label{thmm-change}
Assume both $(G, \xi)$ and $ (\tilde G, \Phi(T,\xi) )$ satisfy
the conditions of Theorem~\ref{thmm-wellposedness}. Then $u$ is the
viscosity solution of PPDE (\ref{PPDE}) with terminal condition $\xi$ if
and only if $\tilde u := \Phi(t, u)$ is the viscosity solution of PPDE
(\ref{PPDE-change}) with terminal condition $\tilde\xi:= \Phi(T,
\xi)$.
\end{thmm}

\begin{pf}One may easily check that $\overline{w} = \Phi
(t, \overline u),
\underline{w} = \Phi(t, \underline u)$, where
\begin{eqnarray*}
\overline{w}(t,\omega) &:=& \inf \bigl\{\psi_t\dvtx \psi\in\overline
C^{1,2}\bigl(\L^t\bigr), \psi^- \mbox { bounded}, \tilde\cL
\psi\ge0, \psi_T \ge\Phi\bigl(T, \xi ^{t,\omega}\bigr) \bigr\} ;
\\
\underline{w}(t,\omega) &:=& \sup \bigl\{\psi_t\dvtx \psi\in
\overline C^{1,2}\bigl(\L ^t\bigr), \psi^+\mbox{ bounded},
\tilde\cL\psi\le0, \psi_T \le\Phi\bigl(T, \xi^{t,\omega}\bigr)
\bigr\}.
\end{eqnarray*}
Then the result follows immediately from Proposition~\ref{prop-perron}
and the arguments in the proof of Theorem~\ref{thmm-wellposedness}.
\end{pf}

We observe that the above operator $\tilde G$ is quadratic in the
$z$-variable, so we need somewhat stronger conditions to ensure the
well-posedness.

%s5 #&#
\section{Partial comparison of viscosity solutions}
\label{sect-partialcomparison}

In this section, we prove Proposition~\ref{prop-comparison}. The proof
is crucially based on the optimal stopping problem reported in Theorem~\ref{thmm-optimal}.

We first prove a lemma. Recall the partition $\{E^i_j, j\ge1\} \subset
\cF_{\textscc{h}_i}$, the constant $n_i$ and the uniform continuous
mappings $\f
^i_{jk}$ and $\psi^i_{jk}$ in (\ref{barC}) corresponding to $ u^1\in
\overline C^{1,2}(\L)$. For $\delta>0$, let $0=t_0<t_1<\cds<t_N=T$ such
that $t_{k+1}-t_k \le\delta$ for $k=0,\ldots, N-1$, and define
$t_{N+1}:=T+\delta$.

%le5.1 #&#
\begin{lem}\label{lem:partial-claim2-prep}
For all $i,j\ge1$, there is a partition $(\tilde E^i_{j,k})_{k\ge
1}\subset\cF_{\textscc{h}_i}$ of $E^i_j$ and a sequence $(p_k)_{k\ge
1}$ taking
values $0,\ldots, N$, such that
\begin{eqnarray*}
&&\textsc{h}_i\in[t_{p_k},t_{p_k+1}) \qquad\mbox{on }
\tilde E^i_{j,k},\qquad \sup_{\omega,\omega'\in\tilde E^i_{j,k}} \bigl\|
\omega_{\cdot \wedge\textscc{h}_i(\omega)}-\omega'_{\cdot \wedge
\textscc{h}_i(\omega')} \bigr\| \le\delta\quad\mbox{and}
\\
&& \min_{\omega\in\tilde E^i_{j,k}}\textscc{h}_i(\omega) =
\textsc{h}_i\bigl(\omega^i_{j,k}\bigr)=:\tilde
t^{i}_{j,k} \qquad\mbox{for some } \omega^i_{j,k}
\in\tilde E^i_{j,k}.
\end{eqnarray*}
\end{lem}

\begin{pf}Since $i,j$ are fixed, we simply denote
$E:=E^i_j$ and $\textsc
{h}:=\textsc{h}_i$.
Denote $E_k := E \cap\{t_k \le\textsc{h}< t_{k+1}\}$, $k\le n$.
Then $\{E_k\}_k\subset\cF_{\textscc{h}}$ forms a partition of $E$.
Since $\O$
is separable, there exists a finer partition $\{E_{k,l}\}_{k,l} \subset
\cF_{\textscc{h}}$ such that, for any $\omega, \omega' \in
E_{k,l}$, $\|\omega_{\cdot\wedge
\textscc{h}(\omega)} - \omega'_{\cdot\wedge\textscc{h}(\omega
')}\| \le\delta$.
%Similarly, we may form a partition $\{E^{i-1}_{k, j,l}: k, j,l \ge1\}
%\subset\cF_{\ch_{i-1}}$. Consider the finer partition $E^i_{k, j,l}
%\cap E^{i-1}_{k', j',l'}\in\cF_{\ch_i}$ over all possible $(k,j,l)$
%and $(k', j', l')$, and renumerate them as $\{\tilde E^i_j, j\ge1\}$.
%For notational simplicity we renumerate $\{E_{k,l}: k,l \ge1\}$ as $\{
%\tilde E_{k}\}_{k\ge1}$.

Next, for each $E_{k,l}$,
%if there exists $\o^{i,j} \in\tilde E^i_j$ such that $\ch_i(\o^{i,j})
%= \min_{\o\in\tilde E^i_j} \ch_i(\o)$, we set $\o^{i,j,m} :=
%\o^{i,j}$ for all $m\ge1$. Otherwise
there is a sequence $\omega^{k,l,m}\in E_{k,l}$ such that
$t_{k,l,m}:=\textsc{h}(\omega
^{k,l,m}) \downarrow\inf_{\omega\in E_{k,l}}\textsc{h}(\omega)$.
Denote $t_{k,l,0}:=
t_{k+1}$. %, where $t_{k_p} \le\ch< t_{k_p+1}$ on $\tilde E_k$.
Define $E_{k,l,m} :=\break E_{k,l} \cap \{t_{k,l,m+1} \le\textsc{h}<
t_{k,l,m}\}
\in\cF_{\textscc{h}_i}$, and renumerate them as $(\tilde E_k)_{k\ge1}$.
We then verify directly that $(\tilde E_k)_{k\ge1}$
defines a partition of $E$ satisfying the required conditions.
\end{pf}

\begin{pf*}{Proof of Proposition~\ref{prop-comparison}} We
only prove
$u^1_0 \le u^2_0$. The inequality for general $t$ can be proved
similarly. %Without loss of generality, we a
Assume $ u^2$ is a viscosity $L$-supersolution and $ u^1\in\ol
C^{1,2}(\L)$ with corresponding hitting times $\textsc{h}_i$, $i\ge0$.
By Proposition~3.14 of~\cite{ETZ1}, we may assume without loss of
generality that
%
%e5.1 #&#
\begin{equation}
\label{Gmonotone} G(t,\omega, y_1, z, \g) - G(t,\omega, y_2,
z, \g) \ge y_2-y_1 \qquad\mbox{for all } y_1 \le
y_2.
\end{equation}
We now prove the proposition in three steps. Throughout the proof, denote
\[
\hat u := u^1 - u^2.
\]
Since $u^1$ is bounded from above and $u^2$ bounded from below, we see
that $\hat u^+$ is bounded.

\textit{Step} 1. We first show that for all $i\ge0$ and
$\omega\in
\O$,
%
%e5.2 #&#
\begin{equation}
\label{partial-claim} \hat u^+_{\textscc{h}_i}(\omega) \le \overline
\cE^L_{\textscc{h}_i(\omega)} \bigl[ \bigl(\hat u^+_{\textscc{h}_{i+1}-}
\bigr)^{\textscc{h}_i,\omega
} \bigr].
\end{equation}
Since $(u^1)^{t,\omega} \in\overline C^{1,2}(\L^t)$, clearly it suffices
to consider $i=0$. Assume on the contrary that
%
%e5.3 #&#
\begin{equation}
\label{partial-c} 2Tc := \hat u^+_0({\mathbf0}) - \overline
\cE^L_0 \bigl[\hat u_{\textscc{h}_1-}^+ \bigr] > 0.
\end{equation}
Recall (\ref{barC}). Notice that $E^0_1 = \O$ and that $\f^0_{1k}(0,
{\mathbf0})$ are constants, and we may assume without loss of generality
that $n_0=1$ and
\[
u^1_t = \psi(t, B), \qquad 0\le t\le\textsc{h}_1,
\]
where $\psi\in C^{1,2}(\L)\cap \mathrm{UC}_b(\L)$ with bounded derivatives.
%%By
%Definition~\ref{defn-barC} (iv), $u^1$ is bounded on $[0, \ch_1]$, and
%denote by $C$ its upper bound.
Denote
\[
X_t := \bigl( \psi_t - u^2_t
\bigr)^+ + ct,\qquad  0\le t\le T.
\]
Since $u^2$ is bounded from below, by the definition of $\Usub$, one
may easily check that
\begin{eqnarray*}
\mbox{$X$ is a bounded process in $\Usub$, and } X_t := \hat
u^+_t+ c t, 0\le t\le\textsc{h}_1.
\end{eqnarray*}
Define
\begin{eqnarray*}
\wh X := X \1_{[0, \textscc{h}_1)} + X_{\textscc{h}_1-}\1_{[\textscc
{h}_1, T]};\qquad Y:=
\overline {\cS}^L[\wh X],\qquad \tau^*:=\inf\{t\ge0\dvtx Y_t=\wh
X_t\}.
\end{eqnarray*}
Applying Theorem~\ref{thmm-optimal} and by (\ref{partial-c}), we have
\begin{eqnarray*}
\overline\cE^L_0[\wh X_{\t^*}] =
Y_0 \ge X_0 = \hat u^+_0({\mathbf0}) = 2Tc +
\overline\cE^L_0\bigl[\hat u_{\textscc{h}_1-}^+\bigr] \ge
Tc+ \overline\cE^L_0[\wh X_{\textscc{h}_1}].
\end{eqnarray*}
Then there exists $\omega^*\in\O$ such that $t^*:=\tau^*(\omega
^*)<\textsc{h}_1(\omega
^*)$. Next, by the $\overline\cE^L$-supermartingale property of $Y$ of
Theorem~\ref{thmm-optimal}, we have
\begin{eqnarray*}
\hat u^+\bigl(t^*,\omega^*\bigr) + c t^* = X_{t^*}\bigl(\omega^*\bigr)
= Y_{t^*}\bigl(\omega^*\bigr) \ge \overline\cE^L_{t^*}
\bigl[X_{\textscc{h}_1^{t^*,\omega
^*}-}^{t^*,\omega^*} \bigr] \ge\overline\cE^L_{t^*}
\bigl[ c\textsc{h}_1^{t^*,\omega^*}\bigr] > ct^*,
\end{eqnarray*}
implying that $0<\hat u^+(t^*, \omega^*)=\hat u(t^*, \omega^*)$.
Since $u^2\in
\Usup$, by (\ref{che})\vspace*{1pt} there exists $\textsc{h}\in\cH^{t^*}$ such that
%
%e5.4 #&#
\begin{equation}
\label{u1-u2>0} \textsc{h}< \textsc{h}_1^{t^*,\omega^*}\quad\mbox{and}\quad
\hat u^{t^*, \omega^*}_t>0\qquad \mbox{for all } t\in\bigl[t^*, \textsc{h}
\bigr].
\end{equation}
Then $X^{t^*, \omega^*}_t = \varphi_t-(u^2)^{t^*, \omega^*}_t$ for
all $t\in
[t^*, \textsc{h}]$, where $\f(t,\omega):= \psi^{t^*,\omega
^*}(t,\omega) + c t$. Observe
that $\f\in C^{1,2}(\L^{t^*})$.
%, a consequence of our assumption $ u^1\in C^{1,2}(\L(\ch_1))$.
Using again the $\overline\cE^L$-supermartingale property of $Y$ of
Theorem~\ref{thmm-optimal}, we see that for all $\t\in\cT^{t^*}$,
\begin{eqnarray*}
\bigl(\f-\bigl(u^2\bigr)^{t^*, \omega^*} \bigr)_{t^*}&=&
X_{t^*}\bigl(\omega^*\bigr)= Y_{t^*}\bigl(\omega^*\bigr) \ge
\overline\cE^L_{t^*} \bigl[Y^{t^*,\omega^*}_{\t\wedge\textscc
{h}}
\bigr] \ge \overline\cE^L_{t^*} \bigl[X^{t^*,\omega^*}_{\t\wedge\textscc{h}}
\bigr] \\
&=& \overline\cE^L_{t^*} \bigl[ \bigl(\f-
\bigl(u^2\bigr)^{t^*, \omega^*} \bigr)_{\t\wedge\textscc
{h}} \bigr].
\end{eqnarray*}
That is, $\f\in\overline\cA^L u^2(t^*, \omega^*)$, and by the viscosity
$L$-supersolution property of $u^2$,
\begin{eqnarray*}
 0 &\le& \bigl\{-\pa_t \f - G\bigl(\cdot, u^2,
\pa_\omega\f, \pa^2_{\omega\omega}\f\bigr) \bigr\}
\bigl(t^*,\omega^*\bigr)\\
& = &-c- \bigl\{\pa_t u^1 + G
\bigl(\cdot, u^2,\pa_\omega u^1,
\pa^2_{\omega\omega} u^1\bigr) \bigr\} \bigl(t^*,\omega^*
\bigr)
\\
& \le& -c- \bigl\{\pa_t u^1 + G\bigl(\cdot,
u^1,\pa_\omega u^1, \pa^2_{\omega\omega}
u^1\bigr) \bigr\} \bigl(t^*,\omega^*\bigr),
\end{eqnarray*}
where the last inequality follows from (\ref{u1-u2>0}) and (\ref
{Gmonotone}). Since $c>0$, this is in contradiction with the subsolution
property of $ u^1$ and thus completes the proof of (\ref{partial-claim}).

%re5.2 #&#
\begin{rem}
The rest of the proof is only needed in the case where $ u^1\in
\overline
{C}^{1,2}(\L)\setminus C^{1,2}(\L)$. Indeed, if $ u^1\in C^{1,2}(\L)$,
then $H_1=T$, and it follows from step 1 that $\hat u^+_0\le
\overline{\cE}^L_0 [\hat u_{T-}^+ ]\le\overline{\cE
}^L_0
[\hat u_T^+ ]=0$, and then $ u^1_0\le u^2_0$. In fact, this is the
partial comparison principle proved in \cite{ETZ1}, Proposition~5.3.
%By a direct change of the time origin, we also prove similarly that $
%u^1\le\tild u^2$ on $\L$.
\end{rem}

\textit{Step} 2. We continue by using the following result
which will
be proved in step~3:
%
%e5.5 #&#
\begin{eqnarray}
\label{lem:partial-claim2}
&& \mbox{For $i\ge1$, $\dbP\in
\cP_L$ and $\cP_L(\dbP, \textsc {h}_i) :=
\bigl\{\dbP '\in\cP_L\dvtx \dbP' = \dbP
\mbox{ on } \cF_{\textscc{h}_i}\bigr\}$, we have}
\nonumber
\\[-8pt]
\\[-8pt]
\nonumber
&&\D_i := \hat u^+_{\textscc{h}_i-} - \esup^{ \dbP}_{\dbP'\in\cP_L(\dbP, \textscc{h}_i)}
\dbE^{\dbP'} \bigl[\hat u^+_{\textscc{h}_{i+1}-} |\cF_{\textscc{h}_i} \bigr]
\le0,\qquad \dbP\mbox{-a.s.}
\end{eqnarray}
Then by standard arguments, we have
\[
\dbE^\dbP \bigl[ \hat u^+_{\textscc{h}_i-} \bigr] \le \sup
_{\dbP'\in\cP_L(\dbP, \textscc{h}_i)} \dbE^{\dbP'} \bigl[\hat u^+_{\textscc{h}_{i+1}-}
\bigr] \le \overline\cE^L_0 \bigl[\hat u^+_{\textscc{h}_{i+1}-}
\bigr].
\]
Since $\dbP\in\cP_L$ is arbitrary, this leads to
$\overline\cE^L_0
[ \hat u^+_{\textscc{h}_i-}]
\le
\overline\cE^L_0
 [\hat u^+_{\textscc{h}_{i+1}-} ]$, and by induction, $
\hat u^+_0
\le
\overline\cE^L_0
 [\hat u^+_{\textscc{h}_i-} ]$, for all $i$.
Notice that $\hat u^+$ is bounded, $\lim_{i\to\infty}\cC
^L_0[\textsc{h}_i < T]
= 0$ by Definition~\ref{defn-barC}(i) and $u^2_{T-}\ge u^2_T$ by the
definition of $\overline{\cU}$. Then, sending $i\to\infty$, we obtain
$\hat u^+_0
\le
\overline\cE^L_0 [\hat u^+_{T-} ] \le
\overline\cE^L_0 [\hat u^+_T ]
= 0$,
which completes the proof of $ u^1_0\le u^2_0$.

\textit{Step} 3. It remains to prove (\ref{lem:partial-claim2}).
Clearly it suffices to prove it on each $E^i_j$. As in the proof of
Lemma~\ref{lem:partial-claim2-prep}, we omit the dependence on the
fixed pair $(i,j)$, thus writing $E:=E^i_j$, $n=n_i$, $\textsc
{h}:=\textsc{h}_i$, $\textsc{h}
_1:=\textsc{h}_{i+1}$, $\f_k := \f^i_{j,k}$, $\psi_k:= \psi
^i_{j,k}$, $\D:= \D
_i$, and let $C$ denote the common bound of $\f_k, \psi_k$ and $\rho$,
the common modulus of continuity function of $\f_k, \psi_k$,
$1\le
k\le n$. We also denote $\tilde E_k := \tilde E^i_{j,k}$, $\omega
^k:=\omega
^i_{j,k}$ and $\tilde t_k:=\tilde t^{i}_{j,k}$, as defined in Lemma~\ref
{lem:partial-claim2-prep}.

Fix an arbitrary $\dbP\in\cP_L$ and $\e>0$.
Since $u^2\in\Usup$, we have $u^2_{\textscc{h}-}\ge u^2_{\textscc
{h}}$. Then, for each
$k$, it follows from (\ref{partial-claim}) that
\[
\hat u^+_{\textscc{h}-}\bigl(\omega^k\bigr) \le\hat
u^+_{\textscc{h}}\bigl(\omega^k\bigr) \le \dbE^{\dbP_k}
\bigl[\bigl(\hat u^+_{\textscc{h}_1-}\bigr)^{\tilde t_k,\omega^k} \bigr] + \e\qquad\mbox{for
some } \dbP_k \in\cP_L^{\tilde t_k}.
\]
Define $\hat\dbP\in\cP_L(\dbP, \textsc{h})$ such that for $\dbP
$-a.e. $\omega\in
\tilde E_k$, the $\hat\dbP^{\textscc{h}(\omega),\omega
}$-distribution of $B^{\textscc{h}(\omega)}$
is equal to the $\dbP_k$-distribution of $B^{\tilde t_k}$, where $\hat
\dbP^{\textscc{h}(\omega),\omega}$ denotes the r.c.p.d. Then $\dbP
$-a.s. on $\tilde E_k$,
\begin{eqnarray*}
&&\dbE^{\hat\dbP} \bigl[\hat u^+_{\textscc{h}_1-} |\cF_{\textscc{h}} \bigr](
\omega) \\
&&\qquad= \dbE^{\hat\dbP^{\textscc{h}(\omega),\omega}} \bigl[\hat u^+ \bigl(\textsc{h}_1
\bigl(\omega\otimes_{\textscc{h}(\omega)} B^{\textscc
{h}(\omega)}_\cdot\bigr)-, \omega
\otimes_{\textscc{h}(\omega)} B^{\textscc{h}(\omega
)}_\cdot \bigr) \bigr]
\\
&&\qquad= \dbE^{\dbP_k} \bigl[\hat u^+ \bigl(\textsc{h}_1\bigl(
\omega\otimes_{\textscc{h}(\omega)} \tilde B^{\tilde t_k}_\cdot\bigr)-,
\omega\otimes_{\textscc{h}(\omega)} \tilde B^{\tilde t_k}_. \bigr) \bigr],
\end{eqnarray*}
where $\tilde B^{\tilde t_k}_s := B^{\tilde t_k}_{s - \textscc
{h}(\omega) + \tilde
t_k}$, $s\ge\textsc{h}(\omega)$. Recalling that $\hat u^+$ is
bounded, $\dbP
$-a.s. this provides
%
%e5.6 #&#
\begin{eqnarray}
\label{Diest} \Delta(\omega) &\le& \hat u^+_{\textscc{h}-}(\omega) -
\dbE^{\hat\dbP} \bigl[\hat u^+_{\textscc{h}_1-} |\cF_{\textscc{h}} \bigr](
\omega)
\nonumber
\\
&\le& \e+ \sum_{k\ge1}\1_{\tilde E_k}(\omega)
\bigl(\hat u^+_{\textscc{h}-}(\omega)-\hat u^+_{\textscc{h}-}\bigl(\omega
^k\bigr) \bigr)
\nonumber
\\
&&{} +\sum_{k\ge1}\1_{\tilde E_k}(\omega)
\dbE^{\dbP_k} \bigl[\bigl(\hat u^+_{\textscc{h}_1-}\bigr)^{\tilde t_k,\omega^k} -
\hat u^+ \bigl(\textsc{h}_1\bigl(\omega\otimes_{\textscc{h}(\omega
)}\tilde
B^{\tilde t_k}_\cdot\bigr)-, \omega\otimes_{\textscc{h}(\omega)}\tilde
B^{\tilde t_k}_\cdot \bigr) \bigr]
\nonumber
\\[-8pt]
\\[-8pt]
\nonumber
&\le& \e + \sum_{k\ge1}\1_{\tilde E_k}(\omega)
\bigl(\hat u_{\textscc{h}-}(\omega)-\hat u_{\textscc{h}-}\bigl(
\omega^k\bigr) \bigr)^+
\\
&&{} +\sum_{k\ge1}\1_{\tilde E_k}(\omega)
\dbE^{\dbP_k} \bigl[ \bigl((\hat u_{\textscc{h}_1-})^{\tilde t_k,\omega^k} -\hat u
\bigl(\textsc{h}_1\bigl(\omega\otimes_{\textscc{h}(\omega
)}\tilde
B^{\tilde t_k}_\cdot\bigr)-, \omega\otimes_{\textscc{h}(\omega)} \tilde
B^{\tilde t_k}_\cdot \bigr) \bigr)^+\nonumber\\
&&\hspace*{292pt}{}\wedge C \bigr].
\nonumber
\end{eqnarray}
We now estimate the above error for fixed $\omega\in\tilde E_k$:

(1) To estimate the terms of the first sum, we recall that
${\mathbf{d}}
_\infty ((\textsc{h}(\omega),\omega),\break   (\tilde t_k, \omega
^k) ) \le2\delta$ on $\tilde
E_k$, by Lemma~\ref{lem:partial-claim2-prep}. Then since $u^1$ is
continuous, it follows from (\ref{barC}) that on $\tilde E_k$,
\begin{eqnarray*}
u^1_{\textscc{h}_i-}(\omega) - u^1_{\textscc{h}_i-}\bigl(
\omega^j\bigr) &=& u^1_{\textscc{h}_i}(\omega) -
u^1_{\textscc{h}_i}\bigl(\omega^j\bigr) \\
&= &\sum
_{l=1}^{n} \bigl[\f_l\bigl(\textsc{h}(
\omega), \omega\bigr)-\f _l\bigl(\tilde t_k,
\omega^k\bigr) \bigr] \psi_l(0,{\mathbf0})\\
& \le &Cn\rho(2
\delta).
\end{eqnarray*}
Moreover, denoting by $\rho_2$ the modulus of continuity of $- u^2\in
\Usub$ in (\ref{USC}), we see that
\begin{eqnarray*}
&&u^2_{\textscc{h}-}\bigl(\omega^k\bigr) -
u^2_{\textscc{h}-}(\omega) \\
&&\qquad= u^2\bigl(\tilde
t_k-, \omega^k\bigr) - u^2(\tilde
t_k-, \omega) + u^2(\tilde t_k-, \omega) -
u^2\bigl(\textsc{h}(\omega)-, \omega\bigr)
\\
&&\qquad\le \rho_2(\delta) + \sup_{\textscc{h}(\omega) - \delta\le t\le
\textsc{h}(\omega)}\bigl[
u^2(t-, \omega) - u^2\bigl(\textsc{h} (\omega)-, \omega
\bigr)\bigr].
\end{eqnarray*}
By the last two estimates, we see that the first sum in (\ref{Diest})
%
%e5.7 #&#
\begin{equation}
\sum_{k\ge1}\1_{\tilde E_k}(\omega) \bigl(\hat
u_{\textscc{h}-}(\omega)-\hat u_{\textscc{h}-}\bigl(\omega^k\bigr)
\bigr)^+ \longrightarrow 0 \qquad\mbox{as } \delta\searrow0.
\end{equation}

(2) Recall from Lemma~\ref{lem:partial-claim2-prep} that $0\le
\textsc{h}(\omega)-\tilde t_k \le\delta$. Then (\ref{chO}) leads to
\begin{eqnarray*}
0 \le\bigl[\textsc{h}_1\bigl(\omega^k
\otimes_{\tilde t_k} \tilde B^{\tilde
t_k}_\cdot\bigr) - \tilde
t_k\bigr] - \bigl[\textsc{h}_1\bigl(\omega
\otimes_{\textscc{h}(\omega)} \tilde B^{\tilde t_k}_\cdot\bigr) - \textsc{h}(
\omega)\bigr] \le \textsc{h}(\omega) - \tilde t_k \le\delta,
\end{eqnarray*}
and therefore, denoting $\eta_\delta(\omega):=\delta+\sup\{
|\omega_s-\omega_t|: 0\le t\le
T, t\le s\le(t+\delta)\wedge T\}$,
%
%e5.8 #&#
\begin{eqnarray}
\label{Diest0} &&{\mathbf{d}}_\infty \bigl( \bigl(\textsc{h}_1
\bigl(\omega^k\otimes _{\tilde t_k} \tilde B^{\tilde
t_k}_\cdot
\bigr) - \tilde t_k,\tilde B^{\tilde t_k} \bigr), \bigl(\textsc
{h}_1\bigl(\omega\otimes _{\textscc{h}(\omega)}\tilde B^{\tilde t_k}_\cdot
\bigr) - \textsc {h}(\omega), \tilde B^{\tilde t_k} \bigr) \bigr)
\nonumber
\\[-8pt]
\\[-8pt]
\nonumber
&&\qquad
\le
\eta_\delta\bigl(\tilde B^{\tilde t_k}\bigr) \le \eta_\delta
\bigl(B^{\tilde t_k}\bigr).
\end{eqnarray}
Then, by using (\ref{barC}) again, we see that
%
%e5.9 #&#
\begin{eqnarray}
\label{Diest2} && \bigl(u^1\bigr)_{\textscc{h}^{\tilde t_k,\omega^k}_1-}^{{\tilde t_k,\omega^k}} -
u^1 \bigl(\textsc{h}_1\bigl(\omega\otimes_{\textscc{h}(\omega)}
\tilde B^{\tilde t_k}_\cdot\bigr)-, \omega\otimes_{\textscc{h}(\omega)}
\tilde B^{\tilde t_k}_\cdot \bigr)
\nonumber
\\
&&\qquad= u^1 \bigl(\textsc{h}_1\bigl(\omega^k
\otimes_{\tilde t_k} B^{\tilde
t_k}_\cdot\bigr),
\omega^k\otimes_{\tilde t_k} \tilde B^{\tilde t_k}_\cdot
\bigr) - u^1 \bigl(\textsc{h}_1\bigl(\omega
\otimes_{\textscc{h}_1(\omega)} \tilde B^{\tilde t_1}_\cdot\bigr), \omega
\otimes_{\textscc{h}_1(\omega)} \tilde B^{\tilde t_k}_\cdot \bigr)
\nonumber
\\
&&\qquad= \sum_{l=1}^{n} \bigl[
\f_l\bigl(\tilde t_k, \omega^k\bigr)
\psi_l \bigl(\textsc{h}_1\bigl(\omega
\otimes_{\textscc{h}(\omega)} \tilde B^{\tilde t_k}_\cdot\bigr) - \textsc{h}(
\omega),\tilde B^{\tilde t_k} \bigr) \\
&&\hspace*{50pt}{}- \f_l\bigl(\textsc{h}(
\omega), \omega\bigr) \psi_l \bigl(\textsc{h}_1\bigl(
\omega^k \otimes_{\tilde t_k} \tilde B^{\tilde t_k}_\cdot
\bigr) - \tilde t_k,\tilde B^{\tilde t_k} \bigr) \bigr]
\nonumber\\
&&\qquad \le C n \bigl[\rho(2\delta) + \rho \bigl(\eta_\delta
\bigl(B^{\tilde
t_k}\bigr) \bigr) \bigr].\nonumber %\le Cn\rho\Big(\eta_\d(B^{\tilde t_k})\Big).
\end{eqnarray}
We now similarly estimate the corresponding term with $u^2$. Since
$\tilde t_k \le\textsc{h}(\omega)$, by~(\ref{USC}) and (\ref
{Diest2}) we have
\begin{eqnarray*}
&& u^2 \bigl(\textsc{h}_1\bigl(\omega
\otimes_{\textscc{h}(\omega)} \tilde B^{\tilde t_k}_\cdot\bigr)-, \omega
\otimes_{\textscc{h}(\omega)} \tilde B^{\tilde t_k}_\cdot \bigr) -
\bigl(u^2_{\textscc{h}_1-}\bigr)^{\tilde t_k,\omega^k}
\nonumber
\\
&&\qquad= \bigl(-u^2\bigr) \bigl(\textsc{h}_1\bigl(
\omega^k\otimes_{\tilde t_k} B^{\tilde
t_k}_\cdot
\bigr)-,\omega^k\otimes_{\tilde t_k} \tilde B^{\tilde t_k}_\cdot
\bigr) - \bigl(-u^2\bigr) \bigl(\textsc{h}_1\bigl(\omega
\otimes_{\textscc{h}(\omega)} \tilde B^{\tilde t_k}_\cdot\bigr)-,\\
&&\hspace*{272pt}{}\omega
\otimes_{\textscc{h}(\omega)} \tilde B^{\tilde t_k}_\cdot \bigr)
\nonumber
\\
&&\qquad\le \rho \bigl({\mathbf{d}}_\infty \bigl(\bigl(\textsc{h}_1
\bigl(\omega\otimes_{\textscc{h}(\omega)} \tilde B^{\tilde t_k}\bigr),\omega
\otimes_{\textscc{h}(\omega)} \tilde B^{\tilde t_k}\bigr), \bigl(\textsc{h}_1
\bigl(\omega^k\otimes_{\tilde t_k} \tilde B^{\tilde t_k}\bigr),
\omega^k\otimes_{\tilde t_k}\tilde B^{\tilde t_j}\bigr) \bigr)
\bigr)
\nonumber
\\
&&\qquad\le \rho \bigl( {\mathbf{d}}_\infty \bigl(\bigl(\textsc{h}(\omega),
\omega\bigr),\bigl(\tilde t_k,\omega^k\bigr) \bigr)
\\
&&\hspace*{8pt}\qquad\quad{} +{\mathbf{d}}_\infty \bigl( \bigl(\textsc{h}_1\bigl(
\omega^k\otimes _{\tilde t_k} \tilde B^{\tilde t_k}\bigr) -
\tilde t_k,\tilde B^{\tilde t_k} \bigr), \bigl(
\textsc{h}_1\bigl(\omega\otimes_{\textscc{h}(\omega)} \tilde B^{\tilde t_k}
\bigr) - \textsc{h}(\omega),\tilde B^{\tilde t_k} \bigr) \bigr) \bigr)
\\
&&\qquad\le \rho \bigl( 2\delta+ \eta_\delta\bigl(B^{\tilde t_k}\bigr)
\bigr).
\end{eqnarray*}
%
%by \reff{Diest0}.
Combining the above with (\ref{Diest2}), this implies that the second
summation in (\ref{Diest}) satisfies
\begin{eqnarray*}
&&\sum_{k\ge1}\1_{\tilde E_k}(\omega)
\dbE^{\dbP_k} \bigl[ \bigl((\hat u_{\textscc{h}_1-})^{\tilde t_k,\omega^k} -\hat u
\bigl(\textsc{h}_1\bigl(\omega\otimes_{\textscc{h}(\omega
)}\tilde
B^{\tilde t_k}_\cdot\bigr)-, \omega\otimes_{\textscc{h}(\omega)} \tilde
B^{\tilde t_k}_\cdot \bigr) \bigr)^+\wedge C \bigr]
\\
&&\qquad \le \sum_{k\ge1}\dbE^{\dbP_k} \bigl[ \bigl(C
n (\rho+\rho_2) \bigl(2\delta+\eta_\delta
\bigl(B^{\tilde t_k}\bigr) \bigr) \bigr)\wedge C \bigr]\1_{\tilde E_k}(
\omega)
\\
&& \qquad\le Cn\overline\cE^L_0 \bigl[ (\rho+
\rho_2) \bigl(2\delta+\eta_\delta (B) \bigr)\wedge C \bigr].
\end{eqnarray*}
One can easily check that $\lim_{\delta\to0} \overline\cE^L_0
[(\rho+
\rho_2) (2\delta+ \eta_\delta(B) )\wedge C ] =0$.
Then by sending $\delta
\to0$ and $\e\to0$ in (\ref{Diest}), we complete the proof of (\ref
{lem:partial-claim2}).
\end{pf*}

%s6 #&#
\section{Consistency of the Perron approach}
\label{sect-baru}

This section is dedicated to the proof of Proposition~\ref
{prop-perron}. We follow the strategy outlined in Section~\ref
{sect-Heuristic}, which is based on the idea in \cite{ETZ1},
Proposition~7.5. However, as pointed out in \cite{ETZ1}, Remark~7.7,
due to fully nonlinearity, the arguments here are much more
involved. We shall divide the proof into several lemmas. As in the
previous section, we may assume without loss of generality that $G$
satisfies the monotonity (\ref{Gmonotone}).

We start with some estimates for viscosity solutions of PDE (\ref{PDEe}).

%le6.1 #&#
\begin{lem}
\label{lem-stabilitye}
Let Assumptions \ref{assum-G} and \ref{assum-Guniform}\textup{(ii)} hold true.
Let $h^i\dvtx \pa Q^\e_t \to\dbR$ be continuous and $v^i$ be the viscosity
solution of the PDE $\mathrm{(E)}^{t,\omega}_{\eps,0}$ with boundary condition
$h^i$, $i=1,2$. Then, denoting $\delta v := v^1-v^2$, $\delta h :=
h^1-h^2$, on
$Q^\e_t$ we have
%
%e6.1 #&#
\begin{eqnarray}
\label{Deltav} \delta v (s,x) \le\overline\cE^{L_0,c_0}_s
\bigl[(\delta h)^+\bigl(\textsc{h}, x+ B^s_{\textscc{h}}\bigr)
\bigr],
\nonumber
\\[-8pt]
\\[-8pt]
\eqntext{\mbox{where } \textsc{h}:= T\wedge\inf\bigl\{r\ge s\dvtx \bigl|x+
B_r^s\bigr| = \e\bigr\}.}
\end{eqnarray}
\end{lem}
\begin{pf} Let $w$ denote the right-hand side of (\ref
{Deltav}). Following
the arguments in Lemma~\ref{lem-rep-boundingPDE}, it is clear that $w$
is the unique viscosity solution of PDE with boundary condition
$(\delta h)^+$,
%
%e6.2 #&#
\begin{equation}
\label{DPDErep} -\partial_tw -\overline g_0
\bigl(Dw,D^2 w\bigr) =0\qquad \mbox{on } {\cal O}_t^\eps.
\end{equation}
Let $K$ be a smooth nonnegative kernel with unit total mass. For all
$\eta>0$, we define the mollification $w^\eta:=w*K^\eta$ of $w$. Then
$w^\eta$ is smooth, and it follows from a convexity argument of Krylov
\cite{krylov00} that $w^\eta$ is a classical supersolution of
%
%e6.3 #&#
\begin{eqnarray}
\label{wetasupersol} -\partial_tw^\eta -\overline
g_0\bigl(Dw^\eta,D^2 w^\eta\bigr)
&\ge& 0 \qquad\mbox{on } {\cal O}_t^\eps,
\nonumber
\\[-8pt]
\\[-8pt]
\nonumber
 w^\eta&=&(\delta
h)^+*K^\eta\qquad\mbox{on } \pa{\cal O}_t^\eps.
\end{eqnarray}
We claim that
%
%e6.4 #&#
\begin{equation}\label{weta}
\begin{tabular}{p{190pt}@{}}
 $\tilde w^\eta+ v^2$ \mbox{supersolution of
the PDE} $\mbox{(E)}^{t,\omega
}_{\eps,0}$, \mbox{where} $\tilde
w^\eta:= w^\eta+ \bigl\|w^\eta-(\delta h)^+
\bigr\|_{\dbL
^\infty(\pa Q_t^\e)}$.
\end{tabular}
\end{equation}
Then, noting that $\tilde w^\eta+ v^2 = w^\eta+h^2 + \|w^\eta
-(\delta
h)^+\|_{\dbL^\infty(\pa Q_t^\e)} \ge h^1 = v^1$ on $\pa Q^\e_t$, we
deduce from the comparison result of Lemma~\ref{lem-PDEcomparison} that
$\tilde w^\eta+ v^2\ge v^1$ on $\overline Q^\e_t$. Sending $\eta
\searrow0$, this implies
that $ w + v^2\ge v^1$, which is the required result.

It remains to prove that $\tilde w^\eta+ v^2$ is a supersolution of
the PDE $\mbox{(E)}^{t,\omega}_{\eps,0}$. Let $(t_0,x_0)\in{\cal
O}_t^\eps
$, $\phi\in C^{1,2}({\cal O}_t^\eps)$ be such that $0=(\phi-\tilde
w^\eta- v^2)(t_0,x_0)=\max(\phi-\tilde w^\eta-v^2)$. Then it follows
from the viscosity supersolution property of $v^2$ that ${\mathbf
L}^{t,\omega
}(\phi-\tilde w^\eta)(t_0,x_0)\ge0$. Hence, at the point $(t_0,x_0)$,
by (\ref{Gmonotone}) and (\ref{wetasupersol}), we have
\begin{eqnarray*}
\mathbf{L}^{t,\omega}\phi &\ge& \mathbf{L}^{t,\omega}\phi-
\mathbf{L}^{t,\omega}\bigl(\phi-\tilde w^\eta\bigr)
\\
&=& -\pa_t w^\eta -g^{t,\omega} \bigl(\cdot,\phi, D
\phi,D^2 \phi\bigr)\\
&&{} +g^{t,\omega} \bigl(\cdot,\phi-\tilde
w^\eta,D \bigl(\phi-w^\eta\bigr),D^2 \bigl(
\phi-w^\eta \bigr) \bigr)
\\
&\ge& -\pa_t w^\eta -g^{t,\omega} \bigl(\cdot,\phi,D
\phi,D^2 \phi\bigr) +g^{t,\omega} \bigl(\cdot,\phi,D \bigl(
\phi-w^\eta\bigr),D^2 \bigl(\phi-w^\eta\bigr)
\bigr)
\\
&\ge& \overline g_0\bigl(Dw^\eta,D^2
w^\eta\bigr) -\a\cdot D w^\eta -\g\dvtx D^2
w^\eta \ge 0,
\end{eqnarray*}
where $|\a|\le L_0$ and $ |\g| \le L_0$, thanks to Assumption~\ref
{assum-G}. This proves (\ref{weta}).
\end{pf}

%s6.1 #&#
\subsection{Viscosity solutions of a discretized path-frozen PDE}

Denote $\Pi_n^\eps:=\Pi_n^\eps(0,T)$ in (\ref{One2}), and by
$\overline
\Pi^\e_n$ its closure. Under Assumption~\ref{assum-xi2} (with $N=1$),
clearly one may extend the mapping $\pi_n \in\Pi^\e_n \longmapsto
\xi(\omega
^{\pi_n})$ continuously to the compact set $\overline\Pi^\e_n$, and we
shall still denote it as $\xi(\omega^{\pi_n})$ for all $\pi_n \in
\overline
\Pi^\e_n$.

We first construct some stopping times, in light of Definition~\ref
{defn-barC}. For $\pi_n\in\Pi^\e_n$ and $(t,x) \in Q^\e_{t_n}$, define
the sequence $\textsc{h}_m^{\e,\pi_n, t, x}:=\textsc{h}_m$ as
follows: First, $\textsc{h}_0
:= t$, and
%
%e6.5 #&#
\begin{eqnarray}
\label{thech}  \textsc{h}_1 &:=&
\inf\bigl\{s\ge t\dvtx | x+B^{ t}_s| = \e\bigr\} \wedge T,\nonumber\\
\textsc{h}_{m+1} &:=& \bigl\{ s> \textsc{h}_m\dvtx
\bigl|B^t_s - B^t_{\textscc{h}_m}\bigr|=\e\bigr\}
\wedge T,\qquad m\ge 1;
\\
\pi^m_n\bigl(t,x,B^t\bigr) &:= &\bigl(
\pi_n, \bigl(\textsc{h}_1, x+ B^t_{\textscc
{h}_1}
\bigr), \bigl(\textsc{h}_{2}, B^t_{\textscc{h}_2} -
B^t_{\textscc{h}_1}\bigr), \ldots,\nonumber\\
&&\hspace*{97pt} \bigl(\textsc{h}_{m},
B^t_{\textscc{h}_m} - B^t_{\textscc{h}_{m-1}}\bigr) \bigr).\nonumber
\end{eqnarray}
It is clear that $\pi^m_n(t,x,B^t) \in\Pi^\e_{n+m}$ whenever
$\textsc{h}_m < T$.

%le6.2 #&#
\begin{lem}
\label{lem-chen}
$\{\textsc{h}^{\e,\pi_n, t, x}_m, m\ge0\}$ satisfies the
requirements of
Definition \ref{defn-barC}\textup{(i)--(ii)}, with $E^m_j = \O^t$ in \textup{(ii)}.
\end{lem}

\begin{pf}
For notational simplicity, we omit the
superscripts $^{\e,\pi_n,
t, x}$. It is clear that $\textsc{h}^{\textscc{h}_m,\omega}_{m+1}
\in\cH^{\textscc{h}_m(\omega)}$
whenever $\textsc{h}_m(\omega) < T$. Next, if $\textsc{h}_m(\omega
) < T$ for all $m$, then
$|B^t_{\textscc{h}_{m+1}}-B^t_{\textscc{h}_m}|(\omega) = \e$ for all
$m$. This contradicts
the fact that $\omega$ is (left) continuous at $\lim_{m\to\infty}
\textsc{h}_m(\omega
)$, and thus $\textsc{h}_m(\omega) = T$ when $m$ is large enough.
Moreover, for
each $m$,
\begin{eqnarray*}
\{\textsc{h}_m < T\} &\subset&\bigl\{\bigl|B^t_{\textscc{h}_{i+1}}-B^t_{\textscc
{h}_i}\bigr|
= \e, i=1,\ldots, m-1\bigr\}\\
& \subset&\Biggl\{\sum_{i=1}^{m-1}\bigl|B^t_{\textscc{h}_{i+1}}-B^t_{\textscc
{h}_i}\bigr|^2
\ge (m-1)\e^2\Biggr\}.
\end{eqnarray*}
Then, for any $L>0$,
%
%e6.6 #&#
\begin{eqnarray}
\label{cLhm} \cC^L_t[\textsc{h}_m<T] &\le&
{1\over(m-1)\e^2} \overline\cE ^L_t \Biggl[\sum
_{i=1}^{m-1}\bigl|B^t_{\textscc{h}_{i+1}}-B^t_{\textscc{h}_i}\bigr|^2
\Biggr]
\nonumber
\\[-8pt]
\\[-8pt]
\nonumber
&\le& {CL^2\over(m-1)\e
^2} \to0 \qquad\mbox{as } m\to\infty.
\end{eqnarray}
Similarly one can show that $\lim_{m\to\infty} \cC^L_s[\textsc
{h}^{s,\omega}_m<T]
=0$ for any $(s,\omega)\in\L^t$.
Finally, for $\omega, \tilde\omega\in\O$ and using the notation in
Definition~\ref{defn-barC}(ii), we have
\begin{eqnarray*}
\textsc{h}_{m+1}(\omega\otimes_{\textscc{h}_m(\omega)} \tilde \omega) &=& T
\wedge\inf\bigl\{t\ge\textsc{h}_i(\omega )\dvtx |\tilde
\omega_{t-\textscc{h}_m(\omega)}| = \e\bigr\} \\
&=& T\wedge \bigl[\textsc{h}_m(
\omega) + \tilde\textsc{h} (\tilde\omega)\bigr],
\end{eqnarray*}
where $\tilde\textsc{h}(\tilde\omega) := \inf\{t\dvtx |\tilde\omega
_t| = \e\}$ is
independent of $\omega$.
Then, given $\textsc{h}_n(\omega) \le\textsc{h}_n(\omega')$,
(\ref{chO}) follows immediately.
\end{pf}

We next prove the existence of the functions $\th^\e_n$, as mentioned
in Section~\ref{sect-Heuristic}, which allows us to construct classical
super and subsolutions in Lemma~\ref{lem-phie} below.

%le6.3 #&#
\begin{lem}\label{lem-the}
Let Assumptions \ref{assum-G}, \ref{assum-Guniform}\textup{(ii)}, \ref
{assum-xi0} and \ref{assum-xi2} with $N=1$ hold true. Then there exists
a sequence of continuous functions $\th^\e_n\dvtx (\pi_n, (t,x))\in
\overline\Pi^\e_{n+1} \mapsto\dbR$, bounded uniformly in $(\eps,n)$,
such that
%
%e6.7 #&#
\begin{eqnarray}
\label{then}&&\theta^\eps_n(
\pi_n;\cdot) \mbox{ is a viscosity solution of $\mathrm{(E)}$}^{t_n,\omega^{\pi_n}}_{\eps,0};\nonumber
\\
&&\th^\e_n(\pi_n; t,x) = \xi \bigl(
\omega^{\pi_n, (t,x)} \bigr)\qquad \mbox{if } t=T, \\
&& \th^\e _n(
\pi_n; t,x) = \th^\e_{n+1}\bigl(
\pi_n, (t,x); t,0\bigr) \qquad\mbox{if } |x|=\e.\nonumber
\end{eqnarray}
\end{lem}
\begin{pf} \textit{Step} 1.
We first prove the lemma in the cases $G = \overline g$ and
$G=\underline g$, as introduced in (\ref{overlineg}).
%Indeed, as in \cite{ETZ1} Section~7 for semilinear PPDEs, in these
%cases we may have explicit representation for the required functions.
For any $n$, denote
\[
\overline\th^\e_{n,n}(\pi_n;
t_n, {\mathbf0}) := \xi\bigl(\omega^{\pi_n}\bigr),
\]
which is continuous for $\pi_n \in\overline\Pi^\e_n$, thanks to
Assumption~\ref{assum-xi2} (with $N=1$).
For $m:= n-1, \ldots, 0$, let $\th:= \overline\th^\e_{n,m}(\pi
_m;\cdot)$
be the unique viscosity solution of the PDE
%
%e6.8 #&#
\begin{eqnarray}
\label{olthen} \overline{\mathbf{L}}\th&:=& -\pa_t \th- \overline g
\bigl(\th, D \th, D^2 \th\bigr) =0\qquad \mbox{in } Q^\e
_{t_m},
\nonumber
\\[-8pt]
\\[-8pt]
\nonumber
 \th(t, x)& =& \overline\th^\e_{n,m+1}\bigl(
\pi_m, (t,x); t, {\mathbf 0}\bigr) \qquad\mbox{on } \pa Q^\e_{t_m}.
\end{eqnarray}
%
%Applying \reff{rep-boundingPDE} repeatedly, the viscosity solution of
%the above PDE has the representation
%\beaa
%\overline\th^\e_{n,m}(\pi_m, t, x)
%&:=&
%\sup_{b\in\cB^t_{L_0}}
%\overline\cE^{L_0}_{t}
% \Big[e^{\int_t^{\ch_{n-m}} b_rdr}\xi\Big(\o^{\pi^{n-m}_m(t,x,B^t)}
%\Big)
% + C_0 \int_t^{\ch_{n-m}} e^{\int_t^s b_rdr}\,ds
%\Big].
%\eeaa
%In particular, it is easily checked that
Applying Lemma~\ref{lem-stabilitye} repeatedly and recalling Assumption~\ref{assum-xi2} (with $N=1$) again, we see that $ \overline\th^\e
_{n,m}(\pi_m; t, x)$ are uniformly bounded and continuous in all
variables $(\pi_m, t, x)$.
Now for any $\pi_m \in\overline\Pi^\e_m$ and $(t,x) \in\overline
Q^\e
_{t_m}$, define
\begin{eqnarray*}
\overline\th^\e_m(\pi_m, t, x) &:=& \sup
_{b\in\cB^t_{L_0}} \overline\cE^{L_0}_{t}
\biggl[e^{\int_t^T b_r\,dr} \limsup_{n\to\infty} \xi \bigl(\omega
^{\pi
^{n-m}_m(t,x,B^t)} \bigr) + C_0 \int_t^T
e^{\int_t^s b_r\,dr}\,ds \biggr].
\end{eqnarray*}
Then, by (\ref{cLhm}),
\begin{eqnarray}
\bigl| \overline\th^\e_m(\pi_m, t, x) -
\overline\th^\e_{n,m}(\pi_m, t, x)\bigr| \le C
\cC_{t_m}^{L_0} [\textsc{h}_{n-m} < T ] \le
{C\over
(n-m-1)\e^2} \longrightarrow0
\nonumber
\\
\eqntext{\mbox{as } n\to\infty.}
\end{eqnarray}
This implies that $\overline\th^\e_{m}(\pi_m; t, x)$ are uniformly
bounded, uniform in $(\e,m)$ and are continuous in all variables $(\pi
_m, t, x)$. Moreover, by stability of the viscosity solutions, we see that
\begin{eqnarray*}
&& \overline\th^\e_m(\pi_m; \cdot)\qquad
\mbox{is the viscosity solution of PDE (\ref{olthen}) in $Q^\e_{t_m}$}\\
&&\hspace*{65pt}\mbox{with the boundary condition}
\\
&& \overline\th^\e_m(\pi_m; T,x) = \xi
\bigl(\omega^{\pi_m, (T,x)} \bigr), \qquad |x|\le\e,\\
&& \overline\th^\e_m(
\pi_m; t,x) = \overline\th^\e_{m+1}\bigl(
\pi_m, (t,x); t,0\bigr),\qquad  |x|=\e.
\end{eqnarray*}
Similarly we may define from $\underline g$ the following $\underline
\th^\e_m$ satisfying the corresponding properties:
\begin{eqnarray*}
\underline\th^\e_m(\pi_m, t, x) &:=& \inf
_{b\in\cB^t_{L_0}} \underline\cE^{L_0}_{t}
\biggl[e^{\int_t^T b_r\,dr} \limsup_{n\to\infty} \xi \bigl(\omega
^{\pi
^{n-m}_m(t,x,B^t)} \bigr) - C_0 \int_t^T
e^{\int_t^s b_r\,dr}\,ds \biggr].
\end{eqnarray*}
\textit{Step} 2. We now prove the lemma for $G$. Given the
construction of
step 1, define
\begin{eqnarray*}
\overline\th^{\e,m}_m(\pi_m; t,x) :=
\overline\th^\e_m(\pi_m; t,x),\qquad  \underline
\th^{\e,m}_m(\pi_m; t,x) := \underline
\th^\e_m(\pi _m; t,x);\qquad m\ge1.
\end{eqnarray*}
For $i=m-1,\ldots, 0$, by Lemma~\ref{lem-PDEcomparison} we may define
$\overline\th^{\e, m}_i$ and $\underline\th^{\e, m}_i$ as the unique
viscosity solution of the PDE (E)$^{t_i,\omega^{\pi_i}}_{\eps,0}$ with
boundary conditions $\overline\th^{\e, m}_i=\overline\th^{\e,m}_{i+1}$
and $\underline\th^{\e, m}_i=\underline\th^{\e, m}_{i+1}$ on
$\partial
Q^\e_{t_i}$.
Note that for $ (t,x) \in\pa Q^\e_{t_{m}}$,
\begin{eqnarray*}
\overline\th^{\e,m}_m(\pi_m; t,x) =\overline
\th^{\e
,m+1}_{m+1}\bigl(\pi _m^{t,x}; t,0
\bigr),\qquad \underline\th^{\e,m}_m(\pi_m; t,x) =
\underline \th^{\e
,m+1}_{m+1}\bigl(\pi_m^{t,x};
t,0\bigr).
\end{eqnarray*}
Since $\underline g\le g^{t,\omega}\le\overline g$, it follows from the
comparison result of the PDEs defined by the operators $\overline g$
and $\underline g$ that
\[
\overline\th^{\e,m}_m(\pi_m; \cdot) \ge
\overline\th^{\e,m+1}_m(\pi_m; \cdot) \ge
\underline\th^{\e,m+1}_m(\pi_m; \cdot) \ge
\underline\th^{\e,m}_m(\pi_m; \cdot) \qquad
\mbox{in } Q^\e_{t_m}.
\]
Then, by an immediate backward induction, the comparison result of\break
Lemma~\ref{lem-PDEcomparison} implies
%
%e6.9 #&#
\begin{eqnarray}
\label{themmonotone} \overline\th^{\e,m}_i(\pi_i;
\cdot) \ge\overline\th^{\e,m+1}_i(\pi_i; \cdot) \ge
\underline\th^{\e,m+1}_i(\pi_i; \cdot) \ge
\underline\th^{\e,m}_i(\pi_i; \cdot)
\nonumber
\\[-8pt]
\\[-8pt]
\eqntext{\mbox{in }
Q^\e_{t_i},  \mbox{ for all } i\le m.}
\end{eqnarray}
Denote $\delta\th^{\e,m}_i:= \overline\th^{\e,m}_i-\underline
\th^{\e
,m}_i$. For any $\pi_i$ and any $(t,x)\in Q^\e_{t_i}$, recall the
notation in (\ref{thech}). Applying Lemma~\ref{lem-stabilitye}
repeatedly, and following similar but much easier arguments as those in
Lemma~\ref{lem:partial-claim2}, we see that
\begin{eqnarray*}
\bigl|\delta\th^{\e,m}_i(\pi_i; t,x)\bigr| \le
\overline\cE^{L_0}_t \bigl[ \bigl|\delta\th^{\e,m}_m
\bigl(\pi^{m-i}_i\bigl(t,x,B^t\bigr);
\textsc{h}_{m-i},0 \bigr)\bigr | \bigr].
\end{eqnarray*}
Note that $\delta\th^{\e,m}_i(\pi_i; t,x) = 0$ when $t=T$. Then,
by (\ref{cLhm}) again,
\begin{eqnarray*}
\bigl|\delta\th^{\e,m}_i(\pi_i; t,x)\bigr| \le C
\cC^{L_0}_t [\textsc{h}_{m-i} < T ] \le
{C\over(m-i-1)\e
^2}\to0\qquad \mbox{as } m\to\infty.
\end{eqnarray*}
Together with (\ref{themmonotone}), this implies the existence of $\th
^\e
_i$ such that $\overline\th^{\e,m}_i \searrow\th^\e_i$,
$\underline\th
^{\e,m}_i \nearrow\th^\e_i$, as $m\to\infty$. Clearly $\th^\e
_i$ are
uniformly bounded and continuous. Finally, it follows from the
stability of the viscosity solutions that $\th^\e_i$ satisfies (\ref{then}).
\end{pf}

%s6.2 #&#
\subsection{Approximating classical super and subsolutions of the PPDE}

We now apply Assumption~\ref{assum-comparison} to $\th^\e_n$ to
construct smooth approximations of $\ol u$ and $\ul u$, namely the $\ol
u^\e$ and $\ul u^\e$ mentioned in Section~\ref{sect-Heuristic}. Define
$\textsc{h}^\e_i := \textsc{h}^{\e, (0, {\mathbf0}), (0, {\mathbf0})}_i$,
that is,
\begin{eqnarray*}
\textsc{h}^\e_0:=0\quad\mbox{and}\quad\textsc{h}^\e_{n+1}
:= T\wedge\inf \bigl\{t \ge\textsc{h}^\e_{n}\dvtx
|B_t - B_{\textscc{h}^\e_{n}}| = \e \bigr\}\qquad \mbox{for all } n\ge0.
\end{eqnarray*}
Let $\hat\pi_n$ denote the sequence $ (\textsc{h}^\e_i,
B_{\textscc{h}^\e_i}
)_{1\le i\le n}$, and $\omega^\e:= \lim_{n\to\infty} \omega
^{\hat\pi_n}$. It
is clear that
%
%e6.10 #&#
\begin{eqnarray}
\label{oerror} \bigl\|\omega- \omega^\e\bigr\|_T\le2\e
\quad\mbox{and} \quad \bigl\|\omega^{\hat\pi
_n}_{\cdot\wedge\textscc{h}
_n}-\omega\bigr\|_{\textscc{h}_{n+1}} \le2
\e \qquad\mbox{for all } n, \omega.
\end{eqnarray}

Recall the common modulus of continuity function $\rho_0$ of $G$ in
Assumption~\ref{assum-Guniform}, and let $\th^\e_n$ be given as in
Lemma~\ref{lem-the}. We then approximate $\th^\e_0$ by a piecewise
smooth processes in $\overline C^{1,2}(\L)$.
%
%le6.4 #&#
\begin{lem}
\label{lem-phie}Under the conditions of Theorem~\ref
{thmm-wellposedness}, with $N=1$ in Assumption~\ref{assum-xi2}, there
exists $\psi^\e\in\overline C^{1,2}(\L)$ bounded from below with
corresponding stopping times $\textsc{h}^\e_n$ such that
%
%e6.11 #&#
\begin{eqnarray}
\label{psie} % \left.\ba{c}
\psi^\e(0,{\mathbf0}) &=&
\th^\e_0(0,{\mathbf0}) + \e+ T\rho_0(2\e),
\nonumber
\\[-8pt]
\\[-8pt]
\nonumber
\psi^\e(T,\omega) &\ge&\xi\bigl(\omega^\e\bigr), \qquad \cL
\psi^\e\ge 0 \qquad\mbox{on } [0,T). % \\
% \mbox{and}
%-\pa_t \psi^\e
%- g^{\ch^\e_n,\o^{\hat\pi_n}}
% (\cdot,\psi^\e, \pa_\o\psi^\e, \pa^2_{\o\o}\psi^\e)
%\ge0, \mbox{on} \big[\ch^\e_n, \ch^\e_{n+1}\big).
%\ea\right.
\end{eqnarray}
\end{lem}

\begin{pf}For notational simplicity, in this proof we omit the
superscript $\e$ and denote $\th_n := \th_n^\e$, $\psi= \psi^\e$ etc.
Moreover, we extend the domain of $\th_n(\pi_n;\cdot)$ to $[t_n,
\infty
)\times\dbR^d$,
\begin{eqnarray*}
\th_n(\pi_n;t,x) := \th \bigl(\pi_n;t
\wedge T, \mathrm{proj}_{\ol
O_\eps
}(x) \bigr),
\end{eqnarray*}
where $ \mathrm{proj}_{\ol O_\eps}$ is the orthogonal projection on
$\ol
O_\eps$, the closed centered ball with radius $\eps$.
We shall construct $\psi$ on each $[\textsc{h}_n, \textsc{h}_{n+1})$
forwardly, by
induction on $n$.

\textit{Step} 1. First, let $\eta>0$, $\lambda>0$ be
small numbers which
will be decided later. Consider PDEs (\ref{PDEe}) and (\ref
{Lbar-PDE}) on
$Q^{\e,\eta}_0$, and recall the operators $\underline{\mathbf{L}}$ and
$\overline{\mathbf{L}}$ at (\ref{Lbar-PDE}). Thanks to Lemma~\ref
{lem-PDEcomparison}, let $v_0^{\eta,\lambda}$, $\overline v_0^{\eta
,\lambda}$ and
$\underline{v}_0^{\eta,\lambda}$ denote the unique viscosity
solutions of
PDEs (E)$^{0, {\mathbf0}}_{\eps,\eta}$, $\overline{\mathbf{L}}v=0$ and
$\underline{\mathbf{L}}v=0$, respectively, with the same boundary condition
$\th_0 +\lambda$ on $\pa Q^{\e,\eta}_0$. %Here $\overline{
%\mathbf{L}}$ is
%defined in \reff{olthen}. We define $\underline{v}_0^{\eta,\l}$
%similarly by using $\underline g$.

By comparison, we have $\underline v_0^{\eta,\lambda}\le v_0^{\eta
,\lambda}\le
\overline v_0^{\eta,\lambda}$. Then, by using the estimate in Lemma~\ref
{lem-stabilitye}, one can easily show that there exist $\eta
_0(\lambda)$ and
$C_0(\lambda)$, which may depend on $L_0$, $\lambda$ and the
regularity of $\th
_0$, such that, for all $\eta\le\eta_0(\lambda)$,
\begin{eqnarray*}
0\le v^{\eta,\lambda}_0-\th_0 \le C_0(
\lambda)\qquad \mbox{on } \overline Q^{\e,\eta}_0\setminus
Q^\e_0 \mbox{ with }  C_0(\lambda) \searrow0,
\mbox{ as } \lambda\searrow0.
\end{eqnarray*}
In particular, the above inequalities hold on $\pa Q^\e_0$. Then, by
the comparison principle, Lemmas \ref{lem-PDEcomparison} and~\ref
{lem-stabilitye}, we have
\[
0\le v^{\eta,\lambda}_0- \th_0 \le C_0(
\lambda) \qquad\mbox{in } \overline Q^{\e,\eta}_0.
\]
Fix $\lambda_0$ such that $C_0(\lambda_0) < {\e\over4}$, and set
$\eta_0:= \eta
_0(\lambda_0)$. Then $v^{\eta_0,\lambda_0}_0< \th_0+{\e\over4}$.
On the other
hand, by Assumption~\ref{assum-comparison}, there exists $v_0\in
C^{1,2}(Q^{\e,\eta_0}_0)$ satisfying
\begin{eqnarray*}
%\label{psiij}
v_0(0, {\mathbf0}) &\le & v^{\eta_0,\lambda_0}_0(0,{
\mathbf0})+ {\e\over4} <\th_0(0,{\mathbf0})+
{\eps\over2},\\
 \mathbf{L}^{0,{\mathbf0}} v_0 &\ge&0\qquad
\mbox{in } Q^{\e,\eta_0}_0,\qquad v_0\ge v^{\eta_0, \lambda_0}_0
\qquad\mbox{on } \pa Q^{\e,\eta_0}_0.
\end{eqnarray*}
By the comparison principle and Lemma~\ref{lem-PDEcomparison}, the last
inequality on $\pa Q^{\e,\eta_0}_0$ implies that
\[
% \label{psiest1}
v_0\ge v^{\eta_0,\lambda_0}_0 \ge
\th_0 \qquad\mbox{on } \overline Q^{\e,\eta_0}_0.
\]
By modifying $v_0$ outside of $Q^{\e, {\eta_0/2}}_0$ and by the
monotonicity (\ref{Gmonotone}), without loss of generality we may
assume $v_0\in C^{1,2}([0, T], \dbR^d)$ with bounded derivatives such that
\begin{eqnarray*}
% \label{w0}
v_0(0,{\mathbf0}) = \th_0(0, {\mathbf0}) +
{\e\over2}, \qquad\mathbf {L}^{0,{\mathbf
0}} v_0 \ge0\qquad
\mbox{in } Q^{\e}_0,\qquad  v_0\ge\th_0\qquad
\mbox{on } \pa Q^{\e}_0.
\end{eqnarray*}
We now define
%
%e6.12 #&#
\begin{eqnarray}
\label{psi0} \psi(t,\omega) := v_0(t, \omega_t) +
{\e\over2} + \rho_0(2\e) (T-t),\qquad  t\in [0,
\textsc{h}_1].
\end{eqnarray}
Note that $(t, \omega_t) \in Q^\e_0$ for $t<\textsc{h}_1$,
$(\textsc{h}_1, \omega_{\textscc{h}_1})
\in\pa Q^\e_0$, and $\th_0$ is bounded. Then
%
%e6.13 #&#
\begin{eqnarray}
\label{psi1}  \psi(0,{\mathbf0}) &=&
\th_0(0,{\mathbf0}) +\e+ T\rho_0(2\e),
\nonumber
\\[-8pt]
\\[-8pt]
\nonumber
v_0(\textsc{h}_1, \omega) &\ge&\th_0(
\textsc{h}_1, \omega _{\textscc{h}_1}) = \th_1(\hat
\pi_1; \textsc{h}_1, 0), \qquad\psi\ge-C \mbox{ on } [0,
\textsc{h}_1].
\end{eqnarray}
Moreover, by monotonicity (\ref{Gmonotone}) again, and by Assumption~\ref
{assum-Guniform} and (\ref{oerror}),
%
%e6.14 #&#
\begin{eqnarray}
\label{psiest1} \cL\psi(t,\omega) &=& \rho_0(2\e) - \pa_t
v_0(t,\omega_t) - G\bigl(t, \omega, \psi, D
v_0(t,\omega_t), D^2v_0(t,
\omega_t)\bigr)
\nonumber
\\
&\ge&\rho_0(2\e) - \pa_t v_0(t,
\omega_t) - G\bigl(t, \omega, v_0(t,\omega_t),
D v_0(t,\omega _t), D^2v_0(t,
\omega_t)\bigr)
\nonumber
\\[-8pt]
\\[-8pt]
\nonumber
&\ge& - \pa_t v_0(t,\omega_t) -
g^{0,{\mathbf0}}\bigl(t, v_0(t,\omega_t), D
v_0(t,\omega _t), D^2v_0(t,
\omega_t)\bigr)
\nonumber
\\
&=& {\mathbf L}^{0, {\mathbf0}}v_{0} (t,\omega_t) \ge0\qquad
\mbox {for } 0\le t <\textsc{h}_1(\omega).\nonumber
\end{eqnarray}
Here we use the fact that $\pa_\omega[v_0(t, \omega_t)] = (D
v_0)(t, \omega_t)$;
see \cite{ETZ1}, Remark~2.9(i).

\textit{Step} 2. Let $ \eta$, $\lambda$, $\delta$ be small
positive numbers
which will be decided later. Set $s_i := (1-\delta)^i T$, $i\ge0$. Since
$\overline O_\e$ is compact, there exists a partition $D_1, \ldots, D_n$
%%of $\{y\dvtx |y|=\e\}$
such that $|y-\tilde y|\le T\delta$ for any $y, \tilde y\in D_j$,
$j=1,\ldots
, n$. For each $j$, fix a point $y_j\in D_j$. Now for each $(i,j)$, let
$v_{ij}^{\eta,\lambda}$ denote the unique viscosity solution of the PDE
(E)$^{s_i,\omega^{(s_i,y_j)}}_{\eps,\eta}$ with the boundary condition
$v^{\eta,\lambda}_{ij}(t,x) = \th_1(s_i, y_j;t, x) +\lambda$ on
$\pa Q^{\e,\eta
}_{s_i}$. Here $\omega^{(s_i,y_j)}$ denotes the linear interpolation of
$(0, {\mathbf0}), (s_i, y_j), (T, y_j)$. Then, by the same arguments as in
step 1, there exist $\eta_0(\lambda)$ and $C_0(\lambda)$, which
may depend on
$L_0$, $\lambda$ and the regularity of $\th_1$, but independent of
$\delta$ and
$(i, j)$, such that for all $\eta\le\eta_0(\lambda)$,
\begin{eqnarray*}
&&0\le v^{\eta,\lambda}_{ij}(t,x)-\th_1(s_i,
y_j;t, x) \le C_0(\lambda) \qquad\mbox{on } \overline
Q^{\e,\eta}_{s_i}\setminus Q^\e_{s_i}
\mbox{ and}\\
&& C_0(\lambda) \searrow0 \qquad\mbox{as } \lambda\searrow0.
\end{eqnarray*}
Following the arguments in step 1, we may fix $\lambda_0$, $\eta_0$,
independently of $\delta$ and $(i,j)$, and there exists $v_{ij} \in
C^{1,2}([s_i, T], \dbR^d)$ with bounded derivatives such that
\begin{eqnarray*}
v_{ij}(s_i, {\mathbf0}) &= &\th_1(s_i,
y_j; s_i, {\mathbf0}) + {\e\over4},\qquad
\mathbf{L}^{s_i, \omega^{(s_i, y_j)}} v_{ij} \ge0 \qquad\mbox{in } Q^{\e}_{s_i},\\
v_{ij}&\ge&\th_1(s_i, y_j;\cdot)\qquad
\mbox{on } \pa Q^{\e}_{s_i}.
\end{eqnarray*}
Denote
\[
E^1_{ij} := \{s_{i+1} <\textsc{h}_1
\le s_i\} \cap\{B_{\textscc{h}_1} \in D_j\} \in
\cF_{\textscc{h}_1}. %& \f^1_{i,j}(t,\o) := v_0(t,\o_t) +{\e\over8} - v_0(s_i, y_j),\q
%\psi^1_{i,j}(t,\o) := v_{ij}(s_i + t, \o_t),&
\]
Here we are using $(i,j)$ instead of $j$ as the index, and clearly
$E^1_{ij}$ form a partition of $\O$. We then define $\psi$ on
$[\textsc{h}_1,
\textsc{h}_2]$ in the form of (\ref{barC}) with $n_1=2$,
%
%e6.15 #&#
\begin{eqnarray}
\label{psi2} \psi_t &:=& \sum_{i,j}
\biggl[ v_0(\textsc{h}_1,B_{\textscc{h}_1}) +
v_{ij}(s_i+t-\textsc{h}_1, B_t -
B_{\textscc{h}_1}) - v_{ij}(s_i, {\mathbf0}) +
{\e\over2} \biggr] \1 _{E^1_{ij}}
\nonumber
\\[-8pt]
\\[-8pt]
\nonumber
&&{} + \rho_0(2\e)
(T-t),\qquad  t\in[\textsc{h}_1, \textsc{h}_2].
\end{eqnarray}
We show that $\psi$ satisfies all the requirements on $[\textsc{h}_1,
\textsc{h}_2]$
when $\delta$ is small enough.

$\bullet$ First, by (\ref{psi2}), we have
\begin{eqnarray*}
\psi_{\textscc{h}_1}& =& \sum_{i,j} \biggl[
v_0(\textsc{h}_1,B_{\textscc
{h}_1}) +
{\e\over2} \biggr] \1_{E^1_{ij}} + \rho_0(2\e) (T-
\textsc{h}_1)\\
&= &v_0(\textsc {h}_1,B_{\textscc{h}_1})
+ {\e\over
2} +\rho_0(2\e) (T-\textsc{h}_1),
\end{eqnarray*}
which is consistent with (\ref{psi0}), and thus $\psi$ is continuous at
$t=\textsc{h}_1$.

$\bullet$ We next check, similar to (\ref{psiest1}), that
%
%e6.16 #&#
\begin{equation}
\label{cLpsi} \cL\psi(t,\omega) \ge0,\qquad \textsc{h}_1 \le t <
\textsc{h}_2.
\end{equation}
Note that $(\textsc{h}_1, B_{\textscc{h}_1}) \in\pa Q^\e_0$ and
$0\le s_i - \textsc{h}_1 \le
s_i - s_{i+1} = \delta s_i \le\delta T$ on $E^1_{ij}$, then
\begin{eqnarray*}
v_0(\textsc{h}_1,B_{\textscc{h}_1}) -
v_{ij}(s_i, {\mathbf0}) + {\e\over
2} &\ge&
\th_1(\textsc{h} _1, B_{\textscc{h}_1};
\textsc{h}_1, {\mathbf0}) - \th_1(s_i,
y_j; s_i, {\mathbf0}) + {\e
\over4} \\
&\ge&
{\e\over4} - \rho_1(3T\delta) \qquad\mbox{on }
E^1_{ij},
\end{eqnarray*}
where $\rho_1$ is the modulus of continuity function of $\th_1$. In
particular, $\rho_1(3T\delta) < {\e\over4} $ when $\delta$ is
small enough.
Now on $E^1_{ij}$, denoting $t_1 := \textsc{h}_1$, $x:= \omega
_{\textscc{h}_1}$, $\tilde
t:= s_i - \textsc{h}_1 + t$, by (\ref{Gmonotone}), Assumption~\ref
{assum-Guniform}(i) and (\ref{oerror}) again, we have
%
%e6.17 #&#
\begin{eqnarray}
\label{psiest2} \cL\psi(t,\omega) &\ge& \cL\psi(t,\omega) - {\mathbf
L}^{s_i,
\omega^{(s_i,
y_j)}}v_{ij} (\tilde t,x)
\nonumber
\\
&=& \rho_0(2\e) - G \bigl(t, \omega, \psi(t,\omega), D
v_{ij}(\tilde t,x), D^2 v_{ij} (\tilde t,x)
\bigr)
\nonumber
\\
&&{} + G \bigl(\tilde t\wedge T, \omega^{(s_i, y_j)}_{\cdot\wedge s_i},
v_{ij}(\tilde t, x), D v_{ij} (\tilde t,x), D^2
v_{ij} (\tilde t,x) \bigr)
\nonumber
\\[-8pt]
\\[-8pt]
\nonumber
&\ge& {\e\over4} - \rho_1(3T\delta) - G \bigl(t,
\omega^{\hat\pi_1}_{\cdot\wedge t_1},v_{ij}(\tilde t, x), D
v_{ij}(\tilde t,x),D^2 v_{ij} (\tilde t,x)
\bigr)
\nonumber
\\
&& {}+ G \bigl(\tilde t\wedge T, \omega^{(s_i, y_j)}_{\cdot\wedge s_i},
v_{ij}(\tilde t, x), D v_{ij} (\tilde t,x), D^2
v_{ij} (\tilde t,x) \bigr)
\nonumber
\\
&\ge& {\e\over4} - \rho_1(3T\delta) -
\rho_0 \bigl({\mathbf{d}}_\infty \bigl(\bigl(t,
\omega^{\hat\pi
_1}_{\cdot\wedge t_1}\bigr), \bigl(\tilde t\wedge T,
\omega^{(s_i, y_j)}_{\cdot\wedge s_i}\bigr) \bigr) \bigr).\nonumber
\end{eqnarray}
Without loss of generality, assume $\e\le T$. Then
\begin{eqnarray*}
&&{\mathbf{d}}_\infty \bigl(\bigl(t, \omega^{\hat\pi_1}_{\cdot\wedge t_1}
\bigr), \bigl(\tilde t\wedge T, \omega^{(s_i, y_j)}_{\cdot\wedge s_i}\bigr)
\bigr) \\
&&\qquad \le |t - \tilde t| + \sup_{0\le s\le T}\biggl|{s\wedge t_1\over t_1}
x - {s\wedge s_i\over s_i} y_j\biggr|
\\
&&\qquad \le \delta T + \sup_{0\le s\le T}\biggl|{s\wedge t_1\over t_1} x -
{s\wedge t_1\over t_1} y_j\biggr| + \sup_{0\le s\le T} \biggl|
{s\wedge t_1\over t_1} y_j - {s\wedge s_i\over s_i}
y_j \biggr|
\\
& &\qquad\le 2\delta T + \e\sup_{0\le s\le T}\biggl|{s\wedge t_1\over t_1} -
{s\wedge s_i\over s_i}\biggr |
\\
&&\qquad  = 2\delta T + \e\biggl[1-{t_1\over s_i}\biggr] \le 2\delta T + \e
\biggl[1-{s_{i+1}\over s_i}\biggr] =3 \delta T.
\end{eqnarray*}
Then
$
\cL\psi(t,\omega) \ge{\e\over4} - [\rho_0+\rho_1](3T\delta)$.
By choosing $\delta$ small enough, we obtain~(\ref{cLpsi}).

$\bullet$ Finally, we emphasize that the bound of $v_{ij}$ and its
derivatives depend only on the properties of $\th_1$ (and the $\eta_0$
which again depends on $\th_1$), but not on $(i,j)$. Then $\psi$
satisfies Definition~\ref{defn-barC}(iii) on $[\textsc{h}_1, \textsc{h}_2]$.
Moreover, since $\th_1$ is bounded, by comparison we see that $\psi
\ge
-C$ on $[\textsc{h}_1, \textsc{h}_2]$.

\textit{Step} 3. Repeating the arguments, we may define
$\psi$ on
$[\textsc{h}
_n, \textsc{h}_{n+1}]$ for all $n$. From the construction and
recalling Lemma~\ref{lem-chen}, we see that $\psi\in\overline C^{1,2}(\L)$ bounded
from below, $\psi(0,{\mathbf0}) = \th_0(0,{\mathbf0}) + \e+ T\rho_0(2\e)$
and $\cL\psi\ge0$ on $[0, T)$. Finally, since $\textsc{h}_n=T$ when
$n$ is
large enough, we see that $\psi(T,\omega) = \psi(\textsc
{h}_n(\omega), \omega) \ge\th_n(\omega
^\e) = \xi(\omega^\e)$.
\end{pf}

Now we are ready to prove the main result of this section.

\begin{pf*}{Proof of Proposition~\ref{prop-perron}} For any
$\e>0$, let
$\textsc{h}^\e_n$, $n\ge0$ and $\psi^\e$ be as in Lemma~\ref
{lem-phie}, and
define $\overline\psi^\e:= \psi^\e+ \rho_0(2\e)$. Then clearly
$\overline\psi^\e\in\overline C^{1,2}(\L)$, $\overline\psi^\e$ is
bounded from below, and
\begin{eqnarray*}
\overline\psi^\e(T,\omega) - \xi(\omega) = \psi^\e(T,
\omega) + \rho_0(2\e) - \xi(\omega) \ge\xi\bigl(\omega^\e
\bigr) -\xi(\omega) + \rho_0(2\e) \ge0,
\end{eqnarray*}
where the last inequality follows from (\ref{oerror}). Moreover, for
$t\in[\textsc{h}_n,\textsc{h}_{n+1})$, by~(\ref{Gmonotone}) again,
\begin{eqnarray*}
\cL\overline\psi^\e(t,\omega) &=& -\pa_t
\psi^\e(t,\omega) - G\bigl(t,\omega, \psi^\e+
\rho_0(2\e), \pa_\omega\psi^\e,
\pa^2_{\omega\omega}\psi^\e\bigr)
\\
&\ge& -\pa_t \psi^\e(t,\omega) - G\bigl(t,\omega,
\psi^\e, \pa_\omega\psi^\e,
\pa^2_{\omega\omega}\psi^\e\bigr) = \cL
\psi^\e(t,\omega) \ge0.
\end{eqnarray*}
Then by the definition of $\overline u$ we see that
\begin{eqnarray*}
\overline u(0, {\mathbf0}) \le \overline\psi^\e(0, {\mathbf0}) =
\psi^\e(0, {\mathbf0}) +\rho_0(2\e) \le
\th^\e_0(0, {\mathbf0}) + \e+ (T+1)\rho_0(2
\e).
\end{eqnarray*}
Similarly, $
\underline u(0, {\mathbf0})
\ge
\th^\e_0(0, {\mathbf0}) - \e- (T+1)\rho_0(2\e)$.
This implies that
\begin{eqnarray*}
\overline u(0, {\mathbf0}) - \underline u(0, {\mathbf0}) \le 2 \bigl( \e+ (T+1)
\rho_0(2\e) \bigr).
\end{eqnarray*}
Since $\e>0$ is arbitrary, this shows that $\overline u(0, {\mathbf0}) =
\underline u(0, {\mathbf0})$. Similarly we can show that $\overline
u(t,\omega
) = \underline u(t,\omega)$ for all $(t,\omega)\in\L$.
\end{pf*}

Fo later use, we conclude this section with a complete well-posedness
result for a special PPDE.

%co6.5 #&#
\begin{cor} \label{lem-overline-w}
Let $G(t,\omega,y,z,\g)=\overline g(y,z,\g)$ satisfy Assumptions
\ref
{assum-G} and \ref{assum-Guniform}, and assume that $\xi$ satisfies
Assumptions \ref{assum-xi0} and \ref{assum-xi2} with $N=1$. Then
$\overline u=\underline u$ and is the unique viscosity solution of the
PPDE (\ref{PPDE}).
\end{cor}

\begin{pf}
We first observe that $\overline g$ satisfies Assumption~\ref
{assum-Guniform}(i). We shall prove in Proposition~\ref{prop-Perron1}
below that it also satisfies Assumption~\ref{assum-comparison}. Then it
follows from the last proof that $\underline{u}=\overline{u}$.
Moreover, the process $\overline{w}$ introduced in (\ref
{rep-boundingPPDE}) is a viscosity solution of PPDE (\ref{PPDE}) with
terminal condition $\xi$. Then it follows from the partial comparison
of Proposition~\ref{prop-comparison} that $\underline{u}\le\overline
{w}\le\overline{u}$, hence equality.
\end{pf}

%s7 #&#
\section{Viscosity properties of $\overline{u}$ and $\underline{u}$}
\label{sect:viscosity}

This section is devoted to the proof of Proposition~\ref
{prop-viscosity}. The idea is similar to the corresponding result in
the PDE literature, and in spirit is similar to the stability of the
viscosity solutions as in \cite{ETZ1}, Theorem~5.1. However, we shall
first establish the required regularities of $\ol u$ and $\ul u$.

%le7.1 #&#
\begin{lem}
\label{lem-ubounded}
Under Assumptions \ref{assum-G} and \ref{assum-xi0}, the processes
$\overline u, \underline u$ are bounded.
\end{lem}

\begin{pf}We shall only prove the result for $\overline u$,
the proof for
$\underline u$ being similar. Fix $(t,\omega)$, and set
\[
\psi(s, \tilde\omega) := C_0(L_0+1) e^{(L_0+1)(T-s)}.
\]
Then $\psi\in C^{1,2}(\L^t)\subset\ol C^{1,2}(\L^t)$, $\psi\ge0$,
$\psi_T \ge C_0(L_0+1) \ge C_0 \ge\xi^{t,\omega}$, and we compute that
\begin{eqnarray*}
(\cL\psi)^{t,\omega}_s &=& (L_0+1)\psi_s -
G^{t,\omega}(\cdot,\psi_s,{\mathbf0},{\mathbf0}) \ge
\psi_s - G^{t,\omega}(\cdot, {\mathbf0}, {\mathbf0}, {\mathbf0})
\\
&\ge& C_0(L_0+1) - C_0 \ge0.
\end{eqnarray*}
This implies that $\psi\in\overline\cD_T^\xi(t,\omega)$, and thus
$\overline u(t,\omega) \le\psi(t,{\mathbf0})$.

On the other hand, by similar arguments one can show that $ - \psi$ is
a classical subsolution of PPDE (\ref{PPDE}) satisfying $-\psi_T \le
\xi
^{t,\omega}$. Then by partial comparison Proposition~\ref{prop-comparison},
$\overline u(t,\omega) \ge-\psi(t,{\mathbf0})$. Hence $
|\overline u(t,\omega)| \le\psi(t, {\mathbf0}) \le C_0(L_0+1) e^{(L_0+1)T}$.
\end{pf}

We next prove that $\overline{u}$ and $\underline{u}$ satisfy a partial
dynamic programming principle.

%le7.2 #&#
\begin{lem}\label{lem-uDP}
Under Assumptions \ref{assum-G} and \ref{assum-xi0}, for $0\le t_1<
t_2\le T$, we have
\[
\overline u(t_1,\omega) \ge\inf \bigl\{\psi_{t_1}\dvtx
\psi\in \overline\cD _{t_2}^{\overline u_{t_2}}(t_1, \omega)
\bigr\}, \qquad \underline u(t_1,\omega) \le\sup \bigl\{\psi_{t_1}
\dvtx \psi\in \underline\cD _{t_2}^{\underline u_{t_2}}(t_1,
\omega) \bigr\}.
\]
\end{lem}

\begin{pf} We only prove the result for $\overline{u}$. For
any arbitrary
$\psi\in\overline\cD_T^\xi(t_1,\omega)$, notice that $\psi
^{t_2, \omega'} \in
\overline C^{1,2}(\L^{t_2})$ and $\psi_{t_2}(\omega') \ge\overline
u_{t_2}^{t_1,\omega}(\omega')$ for any $\omega'\in\O^{t_1}$.
Then $\psi\in\overline
\cD_{t_2}^{\overline u_{t_2}}(t_1,\omega)$, and the result follows.
\end{pf}

The next result shows that the functions $\overline u, \underline u$
are uniformly continuous. We observe that with this regularity in hand,
and following standard techniques, we may prove that the equality holds
in Lemma~\ref{lem-uDP}, so that $\overline u, \underline u$ satisfy a
dynamic programming principle. However, this is not needed for the
present analysis. Moreover, the result is true in degenerate case
$c_0=0$ as well.

%le7.3 #&#
\begin{lem}
\label{lem-uUC}
Under Assumptions \ref{assum-G}, \ref{assum-Guniform}\textup{(ii)} and \ref
{assum-xi0}, we have $\overline u, \underline u\in \mathrm{UC}_b(\L)$.
\end{lem}

\begin{pf}
We only prove the result for $\overline u$.

(i) We first prove that $\overline u$ is uniformly continuous in
$\omega$,
uniformly in $t$.
For $t\in[0,T]$ and $\omega^1,\omega^2\in\O$, denote $\delta
:=\|\omega^1-\omega^2\|_t$.
For $\psi^1 \in\overline\cD^\xi_T(t,\omega^1)$, define
\begin{eqnarray*}
\psi^2(s, \tilde\omega) := \psi^1(s, \tilde\omega) +
\psi(s)\qquad \mbox{where } \psi(s):= e^{(L_0+1)(T-s)} \bigl[\rho_0(
\delta) +\delta \bigr].
\end{eqnarray*}
Notice that $e^{-(L_0+1)s} = e^{-(L_0+1)\textscc{h}_i}
e^{-(L_0+1)(s-\textscc{h}_i)}$,
and one can easily check that $\psi^2 \in\ol C^{1,2}(\L^t)$ with the
same $\textsc{h}_i$ as those of $\psi^1$. Moreover, $\psi^2$ is
bounded from
below, and
\begin{eqnarray*}
\psi^2_T &=& \psi^1_T +
\psi_T\ge\xi^{t,\omega^1}_T +\rho _0(
\delta) \ge\xi^{t,
\omega^2};
\\
\bigl(\cL\psi^2\bigr)^{t,\omega^2}_s &\ge& \bigl(\cL
\psi^2\bigr)^{t,\omega^2}_s-\bigl(\cL\psi^1
\bigr)^{t,\omega^1}_s
\\
&=& (L_0+1)\psi_s -G^{t,\omega^2} \bigl(s,\cdot,
\psi^2,\pa_\omega\psi^1,\pa^2_{\omega\omega}
\psi^1\bigr) \\
&&{}+G^{t,\omega^1} \bigl(s,\cdot,\psi^1,
\pa_\omega\psi^1,\pa ^2_{\omega\omega}
\psi^1\bigr)
\\
&\ge& (L_0+1)\psi_s -\rho_0(\delta) -
L_0\psi_s =\psi_s - \rho_0(
\delta)\ge\delta>0.
\end{eqnarray*}
Then $\psi^2\in\overline\cD^\xi_T(t,\omega^2)$, and therefore
$\overline
u(t,\omega^2) \le\psi^2(t, {\mathbf0})$, implying that
\begin{eqnarray*}
\overline u\bigl(t, \omega^2\bigr) - \psi^1(t, {
\mathbf0}) &\le& \psi^2(t, {\mathbf0}) -\psi^1(t, {
\mathbf0}) = e^{(L_0+1)(T-t)}\bigl[\rho_0(\delta) + \delta\bigr] \\
&\le&
C\bigl[\rho_0(\delta) + \delta\bigr].
\end{eqnarray*}
Since $\psi^1\in\overline\cD^\xi_T(t,\omega^1)$ is arbitrary, we obtain
$\overline u(t, \omega^2) - \overline u(t, \omega^1)
\le
C[\rho_0(\delta) +\delta]$. By symmetry, this shows the required uniform
continuity of $\overline u$ in $\omega$, uniformly in~$t$.

(ii) We now prove that $-\overline{u}$ satisfies (\ref{USC}). Fix $t_1<
t_2\le T$, and consider the process
%
%e7.1 #&#
\begin{eqnarray}
\label{PPDE1rep} \underline w(t, \omega):= \inf_{b\in\cB^t_{L_0}}
\underline\cE^{L_0,
c_0}_t \biggl[e^{\int_t^{t_2} b_r\,dr} \overline u
\bigl(t_2, \omega\otimes_t B^t\bigr) -
C_0 \int_t^{t_2} e^{\int_t^s b_r\,dr}\,ds
\biggr],
\nonumber
\\[-8pt]
\\[-8pt]
\eqntext{(t,\omega)\in[0, t_2]\times\O.}
\end{eqnarray}
By (\ref{rep-boundingPPDE}), $\underline w$ is a viscosity solution of
the PPDE
%
%e7.2 #&#
\begin{eqnarray}
\label{PPDE1} \underline\cL \underline w:= -\pa_t \underline w-
\underline g\bigl(\underline w,\pa_\omega\underline w, \pa
^2_{\omega\omega} \underline w\bigr) = 0,
\nonumber
\\[-8pt]
\\[-8pt]
 \eqntext{t\in[0,t_2),
\omega\in\O, \underline w(t_2, \omega) =\overline u(t_2,
\omega).}
\end{eqnarray}
Recalling (\ref{Goverlineg}) and applying partial comparison principle
Proposition~\ref{prop-comparison} on PPDE (\ref{PPDE1}), we see that
$\psi_{t_1} \ge\underline w(t_1, \omega)$ for any $\psi\in
\overline\cD
^{\overline{u}_{t_2}}_{t_2}(t_1,\omega)$. Then $\overline
u(t_1,\omega) \ge
\underline w(t_1,\omega)$, and thus
\begin{eqnarray*}
&&\overline u(t_2,\omega) - \overline u(t_1,\omega)\\
&&\qquad\le
\overline u(t_2,\omega) - \underline w(t_1,\omega)
\\
&&\qquad= \sup_{b\in\cB^{t_1}_{L_0}} \overline\cE^{L_0,c_0}_{t_1}
\biggl[ \overline u(t_2,\omega) - e^{\int_{t_1}^{t_2} b_r\,dr} \overline u
\bigl(t_2, \omega \otimes_{t_1} B^{t_1}\bigr) +
C_0 \int_{t_1}^{t_2}
e^{\int_{t_1}^s b_r\,dr}\,ds \biggr].
\end{eqnarray*}
Then it follows from (i) and Lemma~\ref{lem-ubounded} that
\begin{eqnarray*}
&&\overline u(t_2,\omega) - \overline u(t_1,\omega) \\
&&\qquad\le
C(t_2-t_1) + C \overline\cE^{L_0,c_0}_{t_1}
\bigl[ \bigl|\overline u(t_2,\omega) - \overline u\bigl(t_2,
\omega\otimes_{t_1} B^{t_1}\bigr) \bigr| \bigr]
\\
&&\qquad\le C(t_2-t_1) + C\overline\cE^{L_0,c_0}_{t_1}
\bigl[\rho \bigl({\mathbf{d}} _\infty\bigl((t_1, \omega),
(t_2,\omega)\bigr) + \bigl\|B^{t_1}\bigr\|_{t_2} \bigr)
\bigr],
\end{eqnarray*}
where $\rho$ is the modulus of continuity of $\overline u(t_2,\cdot)$.
Now it is straightforward to check that $-\overline u$ satisfies (\ref{USC}).

(iii) We finally prove that $\overline{u}$ satisfies (\ref{USC}). This,
together with Lemma~\ref{lem-ubounded} and~(ii), implies that
$\overline u \in \mathrm{UC}_b(\L)$.
For $t_1<t_2$, $\omega\in\O$ and $\psi^2 \in\overline\cD_T^\xi
(t_2,\omega)$, define
\begin{eqnarray*}
\xi_{t_2}(\tilde\omega):=\psi^2(t_2,{
\mathbf0}) +e^{L_0(T-t_2)}\rho _0 \bigl({\mathbf{d}}
_\infty \bigl((t_1, \omega), (t_2, \omega)
\bigr) + \|\tilde\omega \|_{t_2} \bigr),\qquad  \tilde\omega\in
\O^{t_2}
\end{eqnarray*}
and
%
%e7.3 #&#
\begin{eqnarray}
\label{PPDE2rep} \overline w(t,\tilde\omega):= \sup_{b\in\cB^t_{L_0}}
\overline \cE ^{L_0,c_0}_t \biggl[e^{\int_t^{t_2} b_r\,dr}
\xi_{t_2}\bigl(t_2, \tilde \omega\otimes _t
B^t\bigr) + C_0 \int_t^{t_2}
e^{\int_t^s b_r\,dr}\,ds \biggr],
\nonumber
\\[-8pt]
\\[-8pt]
\eqntext{(t,\tilde\omega)\in[t_1, t_2]
\times\O^{t_1}.}
\end{eqnarray}
By Lemma~\ref{lem-ubounded}, we may assume without loss of generality
that $| \psi^2(t_2,{\mathbf0})|\le C$. Then
%
%e7.4 #&#
\begin{eqnarray}
\label{ubarest1} &&\bigl|\overline w(t_1,{\mathbf0}) - \psi^2(t_2,{
\mathbf0})\bigr | \nonumber\\
&&\qquad\le C(t_2-t_1)+ C \overline
\cE^{L_0, c_0}_{t_1} \bigl[\rho_0 \bigl({\mathbf{d}}
_\infty \bigl((t_1, \omega), (t_2, \omega)
\bigr) + \bigl\|B^{t_1}\bigr\| _{t_2} \bigr) \bigr]
\\
&&\qquad \le C \rho \bigl({\mathbf{d}}_\infty \bigl((t_1,
\omega), (t_2, \omega) \bigr) \bigr),\nonumber
\end{eqnarray}
for some modulus of continuity $\rho$.

By (\ref{rep-boundingPPDE}), the process $\overline w$
is a viscosity solution of the PPDE
%
%e7.5 #&#
\begin{eqnarray}
\label{PPDE2} \overline\cL\overline w:=-\pa_t \overline w -
\overline g\bigl(\overline w,\pa_\omega\overline w, \pa^2_{\omega\omega}
\overline w\bigr) = 0,
\nonumber
\\[-8pt]
\\[-8pt]
\eqntext{(t, \tilde\omega)\in[t_1,t_2)
\times\O^{t_1} \mbox{ and } \overline w(t_2, \cdot) =
\xi_{t_2}.}
\end{eqnarray}
Notice that $\xi_{t_2}$ satisfies the conditions of Corollary~\ref
{lem-overline-w}, and therefore $\overline w =\overline{( \overline
w)}$, where $\overline{( \overline w)}$ is defined for PPDE (\ref
{PPDE2}) in the spirit of (\ref{baru}). Then for any $\e>0$, there exists
$\psi^0\in\ol C^{1,2}(\L^{t_1})$ bounded from below such that
%
%e7.6 #&#
\begin{eqnarray}
\label{ubarest2} \psi^0(t_1, {\mathbf0}) &\le& \overline
w(t_1, {\mathbf0}) + \e,\nonumber
\\
\psi^0(t_2, \tilde
\omega) &\ge& \overline w(t_2, \tilde\omega)\quad \mbox{and}
\\
 -
\pa_t\psi^0 - \overline g\bigl(\psi^0,
\pa_\omega\psi^0, \pa ^2_{\omega\omega}
\psi^0\bigr)& \ge&0.\nonumber
\end{eqnarray}
Therefore, for $t\in[t_1, t_2)$, by (\ref{overlineg}) and (\ref
{Goverlineg}), we have
%
%e7.7 #&#
\begin{eqnarray}
\label{ubarest3} \cL\psi^0 &=& -\pa_t\psi^0 -
G\bigl(\cdot, \psi^0,\pa_\omega\psi^0,
\pa^2_{\omega\omega
}\psi^0\bigr)
\nonumber
\\[-8pt]
\\[-8pt]
\nonumber
&\ge& \overline g_0\bigl(\psi^0, \pa_\omega
\psi^0, \pa^2_{\omega\omega} \psi^0\bigr) - G
\bigl(\cdot, \psi^0,\pa_\omega\psi^0,
\pa^2_{\omega\omega
}\psi^0\bigr) \ge 0.
\end{eqnarray}
Now define $\psi^1$ on $\L^{t_1}$ by
%
%e7.8 #&#
\begin{eqnarray}
\label{ubarest4} \psi^1(t, \tilde\omega) &:= &\psi^0(t,
\tilde\omega) \1_{[t_1, t_2)}(t)
\nonumber
\\[-8pt]
\\[-8pt]
\nonumber
&&{} + \bigl[ \psi^2\bigl(t,\tilde
\omega^{t_2}\bigr) + \bigl(\psi^0(t_2, \tilde
\omega) - \psi^2(t_2,{\mathbf0}) \bigr)e^{L_0(t_2-t)}
\bigr] \1_{[t_2, T]}(t),
\end{eqnarray}
where $\tilde\omega^{t_2}_s := \tilde\omega_s - \tilde\omega
_{t_2}$ for $\tilde\omega
\in\O^{t_1}$ and $s\in[t_2, T]$.
Since $\psi^0, \psi^2$ and $-\psi^2(t_2, {\mathbf0})$ are bounded from
below, then so is $\psi^1$.
We shall prove in (iv) below that $\psi^1\in\ol C^{1,2}(\L^{t_1})$.
Then it follows from (\ref{PPDE2}) and (\ref{ubarest2}) that $\psi^0(t_2,
\tilde\omega) \ge\overline w(t_2, \tilde\omega) \ge\psi
^2(t_2,{\mathbf0})$, and
thus $\psi^1(t, \tilde\omega) \ge\psi^2(t,\tilde\omega^{t_2})$
for $t\ge t_2$.
Then, for $t\in[t_2, T]$,
%
%e7.9 #&#
\begin{eqnarray}
\label{ubarest5} \cL\psi^1 &=& -\pa_t \psi^2
+ L_0 \bigl(\psi^1 - \psi^2\bigl(t,\tilde
\omega^{t_2}\bigr) \bigr) - G\bigl(\cdot, \psi^1,
\pa_\omega\psi^2, \pa^2_{\omega\omega
}
\psi^2\bigr)
\nonumber
\\
&\ge& L_0 \bigl(\psi^1 - \psi^2\bigl(t,
\tilde\omega^{t_2}\bigr)\bigr) + G\bigl(\cdot, \psi^2,
\pa_\omega\psi^2, \pa^2_{\omega\omega
}
\psi^2\bigr)
\nonumber
\\[-8pt]
\\[-8pt]
\nonumber
&&{} - G\bigl(\cdot, \psi^1, \pa_\omega
\psi^2, \pa^2_{\omega\omega
}\psi^2\bigr)\\
&\ge&0.\nonumber
\end{eqnarray}
Moreover, by (\ref{ubarest4}), (\ref{ubarest2}) and (\ref{PPDE2}),
\begin{eqnarray*}
\psi^1(T,\tilde\omega) &\ge& \psi^2\bigl(T, \tilde
\omega^{t_2}\bigr) + \bigl(\overline w(t_2, \tilde\omega)-
\psi^2(t_2,{\mathbf0}) \bigr) e^{L_0(t_2-T)}
\\
&\ge& \xi^{t_2, \omega}\bigl(\tilde\omega^{t_2}\bigr) +
\rho_0 \bigl({\mathbf{d}}_\infty \bigl((t_1,
\omega), (t_2, \omega) \bigr) + \|\tilde\omega\| _{t_2}
\bigr) \ge\xi^{t_1, \omega}(\tilde\omega).
\end{eqnarray*}
This, together with (\ref{ubarest3}) and (\ref{ubarest5}), implies that
$\psi^1\in\overline\cD_T^\xi(t_1,\omega)$. Then it follows from
(\ref{ubarest2}) and (\ref{ubarest1}) that
\begin{eqnarray*}
\overline u(t_1,\omega) &\le& \psi^1(t_1, {
\mathbf0}) = \psi^0(t_1, {\mathbf0}) \le\overline
w(t_1,{\mathbf0}) + \e \\
&\le& \psi^2(t_2,{
\mathbf0}) + C\rho \bigl({\mathbf{d}}_\infty \bigl((t_1,
\omega), (t_2,\omega ) \bigr) \bigr) + \e.
\end{eqnarray*}
Since $\psi^2\in\overline\cD_T^\xi(t_2,\omega)$ and $\e>0$ are
arbitrary,
this provides (\ref{USC}).

(iv) It remains to verify that $\psi^1\in\ol C^{1,2}(\L^{t_1})$. Let
$\textsc{h}^0_i, E^{0, i}_j$ correspond to $\psi^0$ and $\textsc
{h}^2_i, E^{2,i}_j$
correspond to $\psi^2$ in Definition~\ref{defn-barC}. Define a random index
\[
I := \inf\bigl\{i\dvtx \textsc{h}^0_i \ge
t_2\bigr\}.
\]
Set $\textsc{h}^1_{ i} := \textsc{h}^0_{ i}$ for $i < I$ and $\textsc
{h}^1_{i} (\omega):= \textsc{h}
^2_{i-I}(\omega^{t_2})$ for $i\ge I$. Moreover, set $E^{1,i}_{2j-1} :=
E^{0, i}_j\cap\{I>i\}$ and $E^{1,i}_{2j} := E^{2, i-I}_j\cap\{I\le i\}
$, $j\ge1$.

Noting that $\textsc{h}^1_{i+1} = \textsc{h}^0_{i+1}\wedge t_2$
whenever $\textsc{h}^0_i <
t_2$, it is clear that $\textsc{h}^1_i$ are $\dbF$-stopping times and
$(\textsc{h}
^1)^{\textscc{h}^1_i(\omega),\omega}_{i+1} \in\cH^{\textscc
{h}^1_i(\omega)}$ whenever $\textsc{h}^1_i(\omega)
< T$. From the construction of $E^{1,i}_j$ one can easily see that $\{
E^{1,i}_j, j\ge1\}\subset\cF_{\textscc{h}_i^1}$ and form a partition
of $\O
^{t_1}$. Moreover, since on each $E^{1,i}_j$, either $\textsc{h}^1_i=
\textsc{h}^0_i$
or $\textsc{h}^1_i = \textsc{h}^2_{i-I}$, Definitions \ref
{defn-barC}(ii)--(iv) are obvious.

It remains to prove
%
%e7.10 #&#
\begin{eqnarray}
\label{cLHi} \bigl\{i\dvtx \textsc{h}^1_i < T\bigr\}
\mbox{ is finite}\quad\mbox{ and}\quad \lim_{i\to
\infty} \cC
^L_{t}\bigl[\bigl(\textsc{h}^1_i
\bigr)^{t,\omega} < T\bigr] =0
\nonumber
\\[-8pt]
\\[-8pt]
\eqntext{\mbox{for any } (t,\omega) \in
\L^{t_1}.}
\end{eqnarray}
Notice that, denoting by $[{i\over2}]$ the largest integer below
${i\over2}$,
\begin{eqnarray*}
\bigl\{\textsc{h}^1_i < T\bigr\} &=& \biggl\{
\textsc{h}^1_i < T, I > \biggl[{i\over2}
\biggr]\biggr\} \cup\biggl\{\textsc{h}^1_i < T, I \le
\biggl[{i\over2}\biggr]\biggr\}
\\
&\subset& \bigl\{\textsc{h}^0_{[{i/2}]} < t_2\bigr
\} \cup\bigl\{\omega\in\O ^{t_1}\dvtx \textsc{h} ^2_{[{i/2}]}
\bigl(\omega^{t_2}\bigr) < T\bigr\}.
\end{eqnarray*}
Then $\{i\dvtx \textsc{h}_i(\omega) < T\}$ is finite for all $\omega$.
%% $\ch^1_i(\o) = T$
%when $i$ is large enough.
Furthermore, for any $L>0$ and $\dbP\in\cP^{t_1}_L$,
\begin{eqnarray*}
\dbP \bigl[\textsc{h}^1_i < T \bigr] &\le& \dbP \bigl[
\textsc{h}^0_{[{i/2}]} < t_2 \bigr] + \dbP \bigl[
\bigl\{\omega\in\O^{t_1}\dvtx \textsc{h}^2_{[{i/
2}]}
\bigl(\omega^{t_2}\bigr) < T\bigr\} \bigr]
\\
&\le& \cC^L_{t_1} \bigl[\textsc{h}^0_{[{i/2}]}<
T \bigr] + \dbE^{\dbP} \bigl[\dbP^{t_2, \omega} \bigl[
\textsc{h}^2_{[{i/
2}]} < T \bigr] \bigr] \\
&\le&
\cC^L_{t_1} \bigl[\textsc{h}^0_{[{i/2}]}<
T \bigr] +\cC^L_{t_2} \bigl[\textsc{h}^2_{[{i/2}]}
< T \bigr],
\end{eqnarray*}
and thus
\begin{eqnarray*}
\lim_{i\to\infty} \cC^L_{t_1}\bigl[
\textsc{h}^1_i < T\bigr] \le \lim_{i\to\infty}
\bigl[ \cC^L_{t_1}\bigl[\textsc{h}^0_{[{i/2}]}<
T\bigr] +\cC^L_{t_2}\bigl[\textsc {h}^2_{[{i/2}]}
< T\bigr] \bigr] = 0.
\end{eqnarray*}
Similarly one can show (\ref{cLHi}) for any $(t,\omega) \in\L^{t_1}$.
\end{pf}

\begin{pf*}{Proof of Proposition~\ref{prop-viscosity}} In view
of Lemmas
\ref{lem-ubounded} and \ref{lem-uUC}, it remains to prove that
$\overline u$ and $\underline u$ are the viscosity $L_0$-supersolution
and subsolution, respectively, of PPDE (\ref{PPDE}). Without loss of
generality, we may assume that the generator $G$ satisfies (\ref
{Gmonotone}), and we prove only that $\overline u$ is a viscosity
$L_0$-supersolution at $(0, {\mathbf0})$.

Assume to the contrary that there exists $\f\in\overline\cA
^{L_0}\overline u(0,{\mathbf0})$ such that $-c := \cL\f(0,{\mathbf0})<0$.
Following the proof of the partial dynamic programming principle of
Lemma~\ref{lem-uDP}, we observe that for any $\psi\in\overline\cD
_T^\xi(0,{\mathbf0})$ and any $(t,\omega)\in\L$, it is clear that
$\psi^{t,\omega
} \in\overline\cD_T^\xi(t,\omega)$ and then $\psi(t,\omega)
\ge\overline
u(t,\omega)$.
By the definition of $\overline u$ in~(\ref{baru}), there exist $\psi
^n\in\ol C^{1,2}(\L)$ such that
%
%e7.11 #&#
\begin{eqnarray}
\label{psin} \delta_n &:=& \psi^n(0,{\mathbf0}) -
\overline u(0,{\mathbf0}) \downarrow0
\qquad
\mbox{as } n\to\infty,
\nonumber
\\[-8pt]
\\[-8pt]
\nonumber
\bigl(\cL
\psi^n\bigr)_s &\ge&0\quad \mbox{and}\quad  \psi^n_s
\ge\overline u_s,\qquad s\in[0,T].
\end{eqnarray}
Let $\textsc{h}$ be the hitting time required in $\overline\cA
^{L_0}\overline u(0,{\mathbf0})$, and since $\varphi\in C^{1,2}(\L)$ and
$\ol u \in \mathrm{UC}_b(\L) \subset\Usup$, without loss of generality, we
may assume
%
%e7.12 #&#
\begin{eqnarray}
\label{perron-ch} \cL\f(t,\omega) \le-{c\over2}\quad \mbox{and}\quad |
\f_t - \f_0| + \ol u_t - \ol u_0
\le{c\over6L_0},
\nonumber
\\[-8pt]
\\[-8pt]
\eqntext{\mbox{for all } t\le\textsc{h}.}
\end{eqnarray}
We emphasize that the above $\textsc{h}$ is independent of $n$. Now
let $\{\textsc{h}
^n_i, i\ge1\}$ correspond to $\psi^n \in\overline C^{1,2}(\L)$.
Since $\f\in\overline\cA^{L_0}\overline u(0,{\mathbf0})$, this implies
for all $\dbP\in\cP_{L_0}$ and $n, i$ that
%
%e7.13 #&#
\begin{equation}
\label{ineq1overlineu} 0 \ge \dbE^\dbP \bigl[(\f- \overline
u)_{\textscc{h}\wedge\textscc
{h}^n_i} \bigr] \ge \dbE^\dbP \bigl[\bigl(\f-
\psi^n\bigr)_{\textscc{h}\wedge\textscc{h}^n_i} \bigr].
\end{equation}
Recall the processes $\a^\dbP$, $\b^\dbP$ in the definition of
$\dbP\in
\cP_L$ [see (\ref{cPL})], and denote $\cG^\dbP\phi:=\a^\dbP\cdot
\pa_\omega
\phi+\frac{1}2(\b^\dbP)^2\dvtx \pa^2_{\omega\omega}\phi$. Then,
applying functional
It\^o formula in (\ref{ineq1overlineu}) and recalling that $\psi^n$
is a
semi-martingale on $[0, \textsc{h}^n_i]$, it follows from (\ref
{psin}) that
\begin{eqnarray*}
\delta_n &\ge& \dbE^\dbP \bigl[ \psi^n_0
- \psi^n_{\textscc{h}\wedge\textscc
{h}^n_i} + \f_{\textscc{h}\wedge
\textscc{h}^n_i} - \f_0
\bigr] \\
&= &\dbE^\dbP \biggl[\int_0^{\textscc
{h}\wedge\textscc{h}^n_i}
\bigl(\pa_t+\cG^\dbP\bigr) \bigl(\f-\psi^n
\bigr)\,ds \biggr]
\\
&\ge& \dbE^\dbP \biggl[\int_0^{\textscc{h}\wedge\textscc{h}^n_i}
\biggl(\frac{c}{2} - G\bigl(\cdot,\f,\pa_\omega\f,
\pa^2_{\omega\omega}\f\bigr) +G\bigl(\cdot,\psi^n,
\pa_\omega\psi^n,\pa^2_{\omega\omega}
\psi^n\bigr)\\
&&\hspace*{213pt}{} + \cG^\dbP\bigl(\f-\psi^n\bigr)
\biggr) \,ds \biggr]
\\
&\ge& \dbE^\dbP \biggl[\int_0^{\textscc{h}\wedge\textscc{h}^n_i}
\biggl(\frac{c}{2} - G\bigl(\cdot,\f,\pa_\omega\f,
\pa^2_{\omega\omega}\f\bigr) +G\bigl(\cdot,\overline u,
\pa_\omega\psi^n,\pa^2_{\omega\omega}
\psi^n\bigr) \\
&&\hspace*{203pt}{}+ \cG^\dbP\bigl(\f-\psi^n\bigr)
\biggr) \,ds \biggr],
\end{eqnarray*}
where the last inequality follows from (\ref{Gmonotone}) and the fact
that $\overline u\le\psi^n$ by (\ref{psin}). Since $\f_0 =
\overline
u_0$, by (\ref{perron-ch}) and (\ref{Gmonotone}), we get
%$\f,\overline u\in C^0(\L)$, and $G$ is uniformly continuous in $y$,
%for possibly smaller ${\ch\wedge\ch^n_i}$:
%
\begin{eqnarray*}
&&\delta_n \ge \dbE^\dbP \biggl[\int
_0^{\textscc{h}\wedge\textscc{h}^n_i} \biggl(\frac{c}{3} - G\bigl(\cdot,
\overline u_0,\pa_\omega\f,\pa^2_{\omega\omega}
\f\bigr) +G\bigl(\cdot,\overline u_0,\pa_\omega
\psi^n,\pa^2_{\omega\omega}\psi^n\bigr) \\
&&\hspace*{235pt}{}+
\cG^\dbP\bigl(\f-\psi^n\bigr) \biggr) \,ds \biggr].
\end{eqnarray*}
Now let $\eta>0$ be a small number. For each $n$, define $\t^n_0 :=
0$, and
\begin{eqnarray*}
\t^n_{j+1} &:=& \textsc{h}\wedge\inf \bigl\{ t\ge
\t^n_j\dvtx \rho _0 \bigl({\mathbf{d}}
_\infty\bigl((t, \omega), \bigl(\t^n_j, \omega
\bigr)\bigr) \bigr) + \bigl|\pa_\omega\f (t,\omega) - \pa_\omega\f
\bigl(\t ^n_j, \omega\bigr)\bigr|
\\
&&\hspace*{34pt}{}+ \bigl|\pa_{\omega\omega}^2\f(t,\omega) - \pa_{\omega\omega
}^2
\f\bigl(\t^n_j, \omega\bigr)\bigr| + \bigl|\pa_\omega
\psi ^n(t,\omega) - \pa_\omega\psi^n\bigl(
\t^n_j, \omega\bigr)\bigr|
\\
&&\hspace*{140pt}{} + \bigl|\pa_{\omega\omega}^2\psi^n(t,\omega) -
\pa_{\omega
\omega}^2 \psi^n\bigl(\t^n_j,
\omega\bigr)\bigr| \ge \eta \bigr\}.
\end{eqnarray*}
Recalling Definitions \ref{defn-barC}(iii)--(iv), we see the uniform
regularity of $\psi^n$ on $[0, \textsc{h}^n_i]$ for each $i$. Then, together
with the smoothness of $G$ and $\f$, one can easily check that $\t^n_j
\uparrow\textsc{h}$ as $j\to\infty$. Thus
\begin{eqnarray*}
\delta_n & \ge& \biggl[\frac{c}{3} - C\eta\biggr]
\dbE^\dbP\bigl[{\textsc{h}\wedge\textsc{h}^n_i}
\bigr] \\
&&{}+ \sum_{j\ge0}\dbE^\dbP \bigl[ \bigl(
\t^n_{j+1}\wedge\textsc{h}^n_i-
\t^n_j \wedge\textsc {h}^n_i
\bigr)
\\
&&\hspace*{45pt}{}\times \bigl(G\bigl(\cdot,\overline u_0,\pa_\omega
\psi^n,\pa^2_{\omega
\omega
}\psi^n\bigr) - G
\bigl(\cdot,\overline u_0,\pa_\omega\f,\pa^2_{\omega\omega}
\f\bigr) \\
&&\hspace*{198pt}{}+ \cG^\dbP\bigl(\f-\psi^n\bigr)
\bigr)_{\tau^n_j} \bigr]
\\
&=& \biggl[\frac{c}{3} - C\eta\biggr]\dbE^\dbP\bigl[{\textsc{h}
\wedge\textsc{h}^n_i}\bigr]\\
&&{} + \sum
_{j\ge0}\dbE^\dbP \biggl[ \bigl(\t^n_{j+1}
\wedge\textsc{h}^n_i-\t^n_j
\wedge\textsc {h}^n_i \bigr)
\\
&& \hspace*{45pt}{}\times\biggl(\a_{\t^n_j}\cdot\pa_\omega\bigl(\psi^n-\f
\bigr) + {1\over2} \b _{\t^n_i}^2 \dvtx
\pa^2_{\omega\omega}\bigl(\psi^n-\f\bigr) +
\cG^\dbP\bigl(\f-\psi^n\bigr)_{\tau
^n_j} \biggr)
\biggr]
\end{eqnarray*}
for some appropriate $\a_{\t^n_j}, \b_{\t^n_j}$. Now
%Since $\xi\in\mbox\mathrm{ UC}_b(\cF_T,\dbR)$, we may
choose $\dbP_n\in\cP_{L_0}$ such that $\a^{\dbP_n}_t = \a_{\t^n_j}$,
$\b^{\dbP_n}_t = \b_{\t^n_j}$ for all $\t^n_j \le t < \t^n_{j+1}$.
% \beaa
% \big\{G(\cdot,\overline u_0,\pa_\o\psi^n,\pa^2_{\o\o}\psi^n)
% - G(\cdot,\overline u_0,\pa_\o\f,\pa^2_{\o\o}\f)
% + \cG^{\dbP_n}(\f-\psi^n)
% \big\}_{\tau^n_j}
% &=& 0.
% \eeaa
Then $\delta_n\ge[\frac{c}{3} - C\eta]\dbE^{\dbP_n}[{\textsc
{h}\wedge\textsc{h}
^n_i}]$. Set $\eta:= {c\over6C}$, send $i\to\infty$ and recall from
Definition~\ref{defn-barC} that $\lim_{i\to\infty}\cC
^{L_0}_0(\textsc{h}^n_i <
T) = 0$. This leads to $\delta_n \ge{c\over6}\dbE^{\dbP
_n}[\textsc{h}]\ge
\underline\cE^{L_0}_0[\textsc{h}]$, and by sending $n\to\infty$,
we obtain
$\underline\cE^{L_0}_0[\textsc{h}] =0$. However, since $\textsc
{h}\in\cH$, by \cite
{ETZ1}, Lemma~2.4, we have $\underline\cE^{L_0}[\textsc{h}] >0$.
This is a
contradiction.
\end{pf*}

%s8 #&#
\section{On Assumptions \texorpdfstring{\protect\ref{assum-comparison}}{3.8} and \texorpdfstring{\protect\ref{assum-Guniform}}{3.2}(i)}
\label{sect-Perron}
%\setcounter{equation}{0}

%s8.1 #&#
\subsection{Sufficient conditions for Assumption \texorpdfstring{\protect\ref{assum-comparison}}{3.8}}
In this subsection we discuss the validity of our Assumption~\ref
{assum-comparison} which is clearly related to the classical Perron
approach, the key argument for the existence in the theory of viscosity
solutions, as shown by Ishii \cite{Ishii}. However, our definition of
$\overline v$ and $\underline v$ involves classical supersolutions and
subsolutions, while the classical definition in \cite{Ishii} involves
viscosity solutions. We remark that Fleming and Vermes \cite{FV1,FV2}
have some studies in this respect. %, see Remark~\ref{rem-Perron1} (ii)
%below.
The main issue here is to approximate viscosity solutions by classical
supersolutions or subsolutions. This is a difficult problem which
requires some restrictions on the nonlinearity. In this section, we
provide some sufficient conditions, and we hope to address this issue
in a more systematic way in future research.
% which require, in particular, that $g$ is convex in $\g$. We
%emphasize however that we may also provide sufficient conditions for
%nonlinearities which have no convexity structure in $\g$, see e.g.
%Jakobsen \cite{Jakobsen}, Pham and Zhang \cite{PZ} in the context of
%the Isaacs equation.

For ease of presentation, we first simplify the notation in Assumption~\ref{assum-comparison}. Let
%
%e8.1 #&#
\begin{eqnarray}
\label{Perron-Oet}
O&:=&\bigl\{x\in
\dbR^d\dvtx |x|<1\bigr\},\qquad \overline O :=\bigl\{x\in\dbR^d
\dvtx |x|\le1\bigr\}, \nonumber\\
 \pa O &:=&\bigl\{x\in\dbR^d\dvtx |x|=1\bigr\};
\nonumber
\\[-8pt]
\\[-8pt]
\nonumber
Q &:=& [0, T)\times O,\qquad \overline Q := [0, T]\times\overline O,\\
 \pa Q &:=&
\bigl([0,T]\times\partial O \bigr) \cup \bigl(\{T\}\times O \bigr).\nonumber
\end{eqnarray}
We shall consider the following (deterministic) PDE on $Q$:
%
%e8.2 #&#
\begin{eqnarray}
\label{Perron-PDEe} \mathbf{L}v &:=& -\pa_t v - g \bigl(s, x, v, D v,
D^2 v\bigr) = 0 \qquad\mbox{in } Q \quad\mbox{and}
\nonumber
\\[-8pt]
\\[-8pt]
\nonumber
 v &= &h\qquad \mbox{on } \pa Q.
\end{eqnarray}
We remark that in (\ref{PDEe}) the generator $g$ is independent of $x$.

%as8.1 #&#
\begin{assum}\label{assum-Perron1}
{(i)} $g$ and $h$ are continuous in $(t,x)$;

{(ii)} $g$ is uniformly Lipschitz continuous in $(y,z,\g)$ and
uniformly elliptic in $\g$.
%\\
%{(iii)} The PDE \reff{Perron-PDEe} satisfies existence and
%comparison in the sense of viscosity solutions within the class of
%bounded functions.
\end{assum}
As in Lemma~\ref{lem-PDEcomparison}, under the above assumption, we see
that PDE (\ref{Perron-PDEe}) has a unique viscosity solution $v$, and
the comparison principle holds in the sense of viscosity solutions
within the class of bounded functions.
%More precisely, the last item (iii) states that for any bounded
%functions $v^1,v^2$ satisfying $\mathbf{L}v^1\le0\le\mathbf{L}v^2$ on
%$Q$, in the sense of viscosity solutions, and $v^1\le h\le v^2$ on $
%\pa Q$, we have $v^1\le v^2$.
Define
\begin{eqnarray*}
\overline v(t,x) &:=& \inf \bigl\{w(t,x)
\dvtx w \mbox{ classical supersolution of PDE (\ref{Perron-PDEe})} \bigr\},
\\
\underline v(t,x) &:=& \sup \bigl\{w(t,x)\dvtx w \mbox{ classical subsolution of
PDE (\ref{Perron-PDEe})} \bigr\}.
\end{eqnarray*}
By the comparison principle we have $\underline v \le v \le\overline
v$. %, where $v$ denotes the unique viscosity solution of PDE
%\reff{Perron-PDEe}.

Denote $\dbS^d_+ := \{\g\in\dbS^d\dvtx \g\ge{\mathbf0}\}$. The following
proposition is the main result of this section:
%
%pr8.2 #&#
\begin{prop}
\label{prop-Perron1}
Under Assumption~\ref{assum-Perron1}, we have $\overline v =
\underline
v$ if $g$ is either convex in $\g$ or the dimension $d\le2$.
%in the following three cases:
%\\
%{(i)} $g$ is convex in $(y, z, \g)$, $ g_\d(\cdot,\g):=\inf_{A\in
%\dbS^d_+}\big\{ g(\cdot,\g+A)- \d I_d\dvtx A\big\}>-\infty$ for $0\le\d\le
%c_0$, for some $c_0>0$, and $g_\d\longrightarrow g$ as $\d\searrow0$,
%\\
%{(ii)} $g$ is convex in $\g$ and uniformly elliptic: for some
%constant $c_0>0$,
%\beaa
%g(\cd, \g_1) - g(\cd, \g_2) \ge c_0 I_d\dvtx (\g^1-\g^2) &\mbox{for any}&
%\g^1\ge\g^2.
%\eeaa
%{(iii)} $g$ is uniformly elliptic and $d\le2$.
\end{prop}
\begin{pf}
For the case $d\le2$ we refer to Pham and Zhang
\cite{PZ}.
Below, we prove the result only for the case when $g$ is convex in $\g
$. As in (\ref{Gmonotone}), we assume without loss of generality that
%
%e8.3 #&#
\begin{eqnarray}
\label{gstrict-decr} g(\cdot,y_1,\cdot)-g(\cdot,y_2,\cdot) \le
y_2-y_1 \qquad\mbox{for all } y_1 \ge
y_2.
\end{eqnarray}

For any $\a>0$, we define $O^\a:=\{x\in\dbR^d\dvtx |x|< 1+\a\}$,
$Q^\delta:=
[0,(1+\a)T)\times O^\a$, and similar to (\ref{Perron-Oet}), define
their closures and boundaries. Let $\mu,\eta$ be smooth mollifiers on
$Q$ and $Q^1 \times\dbR\times\dbR^d\times\dbS^d$, and define for
any $\a'>0$,
\begin{eqnarray*}
h_\a(t,x) &:=& (h * \mu_\a) \biggl(\frac{t}{1+\a},
\frac{x}{(1+\a)} \biggr), \qquad(t,x)\in\overline Q^\a,
\\
g_{0}(t,x,y,z,\gamma) &:=& \min_{(t',x')\in Q} {\bigl\{g
\bigl(t',x',y,z,\gamma\bigr) +2 \rho_0
\bigl(\bigl|t-t'\bigr|+\bigl|x-x'\bigr|\bigr)\bigr\}},
\\
g_{\a'}&:=&(g_{0}*\eta_{\a'}),\qquad  (t,x,y,z,\g)\in
Q^1\times\dbR\times\dbR^d\times\dbS^d.
\end{eqnarray*}
By the uniform continuity of $g$, we have $c(\a'):=\|g-g_{\a'}\|
_{\infty
}\rightarrow0$ as $\a'\searrow0$. Set
\[
\underline{g}_{\a'} := g_{\a'}-c\bigl(\a'
\bigr)\quad \mbox{and}\quad \overline{g}_{\a'} := g_{\a'}+c\bigl(
\a'\bigr).
\]
By our assumptions on $g$ and $h$, it follows from Theorem~14.15 of
Lieberman~\cite{Lieberman} that there exist $\underline{v}_{\a,\a'},
\overline{v}_{\a,\a'}\in C^{1,2}(Q^\a) \cap C(\overline Q^\a)$
solutions of the equations
\begin{eqnarray*}
(\underline{E}_{\a,\a'} )\dvtx -\pa_t v -
\underline{g}_{\a'}\bigl(\cdot,v, D v,D^2 v\bigr)&=& 0\qquad
\mbox{in } Q^\a\quad \mbox{and}\quad v=h_\a\qquad \mbox{on } \partial
Q^\a,
\\
(\overline{E}_{\a,\a'} )\dvtx -\pa_t v -\overline{g}_{\a'}
\bigl(\cdot,v, D v,D^2 v\bigr)&=& 0\qquad \mbox{in } Q^\a\quad
\mbox{and}\quad v=h_\a \qquad\mbox{on } \partial Q^\a,
\end{eqnarray*}
respectively. In particular, their restriction to $\overline Q$ are in
$C^{1,2}(\overline Q)$. By the comparison principle, $\underline
{v}_{\a
,\a'} \le\overline{v}_{\a,\a'}$. Moreover, it follows from (\ref
{gstrict-decr}) that
\begin{eqnarray*}
\overline{g}_{\a'}\bigl(\cdot,y+2 c\bigl(\a'\bigr),\cdot
\bigr) \le \overline{g}_{\a'}(\cdot,y,\cdot)-2 c\bigl(
\a'\bigr) = \underline{g}_{\a'}(\cdot,y,\cdot).
\end{eqnarray*}
This shows that $v_{\a,\a'}+2c(\alpha')$ is a classical supersolution
of $(\overline{E}_{\a,\a'})$, and therefore
\[
\underline{v}_{\a,\a'}+2c\bigl(\a'\bigr) \ge
\overline{v}_{\a,\a'} \ge \underline{v}_{\a,\a'}.
\]
Additionally, notice that the solutions $ \underline{v}_{\a,\a'},
\overline{v}_{\a,\a'}$ are bounded uniformly in ${\a,\a'}$ for
${\a,\a
'}$ small enough. The generators $\underline{g}_{\a'}, \overline
{g}_{\a
'}$ have the same uniform ellipticity constants as $g$, and they verify
the hypothesis of Theorem~14.13 of Liebermann \cite{Lieberman}
uniformly in $\a'$. Therefore $\underline{v}_{\a,\a'}, \overline
{v}_{\a
,\a'}$ are Lipschitz continuous with the same Lipshitz constant for all
$\a,\a'$. Then, denoting $\overline{h}_{\a,\a'}:=\overline{v}_{\a
,\a
'} |_{\partial Q}$ and $\underline{h}_{\a,\a'}:=\underline
{v}_{\a,\a
'} |_{\partial Q}$, this implies that
\begin{eqnarray}
c\bigl(\a,\a'\bigr) := \max \bigl\{\|\overline{h}_{\a,\a'}-h
\|_\infty, \|\underline{h}_{\a,\a'}-h\|_\infty \bigr\}
\longrightarrow 0
\nonumber
\\
\eqntext{\mbox{as } \a\rightarrow0, \mbox{ uniformly in }
\a'.}
\end{eqnarray}
Now for fixed $\e>0$, choose $\a_0,\a'_0>0$ so that $c(\a_0,\a
')<\e/4$
for all $\a'>0$, and $c(\a'_0)\leq\e/4$. Then $\overline{w}_{\a
_0,\a
'_0}:= \overline{v}_{\a_0,\a'_0}+c(\a_0,\a'_0)$ and $\underline
{w}_{\a
_0,\a'_0}:= \underline{v}_{\a_0,\a'_0}-c(\a_0,\a'_0)$
are respectively the classical supersolution and subsolution of (\ref
{Perron-PDEe}) on $\overline Q$. Thus $\underline{w}_{\a_0,\a'_0}\le
\underline v$ and $\overline{w}_{\a_0,\a'_0} \ge\overline v$. Therefore,
\begin{eqnarray*}
\overline v - \underline v& \le&\overline{w}_{\a_0,\a'_0} - \underline
{w}_{\a_0,\a'_0} = \overline{v}_{\a_0,\a'_0} - \underline{v}_{\a
_0,\a
'_0}
+ 2c\bigl(\a_0,\a'_0\bigr) \le2 c\bigl(
\a_0'\bigr) + 2c\bigl(\a_0,
\a'_0\bigr) \\
&\le&\e.
\end{eqnarray*}
Then it follows from the arbitrariness of $\e$ that $\overline v =
\underline v$.
\end{pf}

%s8.2 #&#
\subsection{A weaker version of Assumption
\texorpdfstring{\protect\ref{assum-Guniform}\normalfont{(i)}}{3.2(i)}}
\label{sect:Guniform}
We remark that, while seemingly reasonable, the uniform continuity of
$G$ in $(t,\omega)$ is violated even for semilinear PPDEs when the
diffusion coefficient $\si$ depends on $(t,\omega)$. In this
subsection we
weaken the uniform regularity in Assumption~\ref{assum-Guniform}
slightly so as to fit into the framework of Pham and Zhang \cite{PZ},
which deals with path-dependent Bellman--Isaacs equations associated to
stochastic differential games.
%We also note that, in a more recent paper Ren, Touzi, and Zhang
%\cite{RTZ}, by using a quite different approach we established the
%wellposedness for semilinear PPDEs with random and possibly degenerate
%diffusion $\si$. However, it is still unclear to us whether or not
%that approach can be extended to fully nonlinear case.

%as8.3 #&#
\begin{assum}
\label{assum-Guniform2}
There exist a modulus of continuity functions $\rho_0$, $\tilde\rho_0$
such that, for any $(t, \omega), (\tilde t, \tilde\omega) \in\L$
and any $(y,
z, \g)$,
\begin{eqnarray*}
&&\bigl|G(t,\omega, y,z,\g) - G(\tilde t, \tilde\omega, y, z, \g)\bigr|
\\
&&\qquad\le \tilde\rho
_0\bigl(|t-\tilde t|\bigr) \bigl[|z|+|\g|\bigr] + \rho_0 \bigl({
\mathbf{d}}_\infty \bigl((t,\omega), (\tilde t,\tilde\omega) \bigr)
\bigr).
\end{eqnarray*}
\end{assum}

Recall the parameters $\e,\delta,\eta_0$ and the functions $v_{ij}$
introduced in the proof of Lemma~\ref{lem-phie}.
Notice that Assumption~\ref{assum-Guniform} is used only in the proof
of Lemma~\ref{lem-phie}, more precisely in (\ref{psiest1}) and (\ref
{psiest2}). We also note that the smooth functions $v_{ij}$ are
typically constructed as the classical solution to some PDE, as in
Section~\ref{sect-Perron} and in~\cite{PZ}, and thus satisfy certain
estimates. Assume the following:
%
%e8.4 #&#
\begin{equation}\label{vbound}
\begin{tabular}{p{300pt}@{}}
There exists a constant $C_{\eta_0}>0$, which may
depend on $\eta _0$ (and $\e$), but is independent
of $\delta$, such that $|D v_{ij}(t,x)| \le
C_{\eta_0}, |D^2 v_{ij}(t,x)|\le C_{\eta
_0}$
for all $(t, x) \in\overline Q_0^\e.$
\end{tabular}
\end{equation}
We claim that Lemma~\ref{lem-phie}, hence our main result, Theorem~\ref
{thmm-wellposedness}, still holds true if we replace Assumption~\ref
{assum-Guniform} by (\ref{vbound}) and Assumption~\ref{assum-Guniform2}.

Indeed, in (\ref{psiest1}), note that
\begin{eqnarray*}
&& G \bigl(t, \omega, \bigl(v_0, D v_0,
D^2v_0\bigr) (t,\omega_t) \bigr) -
g^{0,{\mathbf0}} \bigl(t, \bigl(v_0, D v_0,
D^2v_0\bigr) (t,\omega_t) \bigr)
\\
&&\quad= G \bigl(t, \omega, \bigl(v_0, D v_0,
D^2v_0\bigr) (t,\omega _t) \bigr) - G
\bigl(t, {\mathbf0}, \bigl(v_0, D v_0,
D^2v_0\bigr) (t,\omega_t) \bigr) \le
\rho_0(\e),
\end{eqnarray*}
thanks to Assumption~\ref{assum-Guniform2}. Thus we still have (\ref
{psiest1}).

To see (\ref{psiest2}) under our new assumption, we first note that as
in (\ref{psiest2}) and by~(\ref{Gmonotone}),
\begin{eqnarray*}
\cL\psi(t,\omega) &\ge& \rho_0(2\e) +{\e\over4} - \rho
_1(3T\delta) - G \bigl(t, \omega, v_{ij}(\tilde t, x), D
v_{ij}(\tilde t,x), D^2 v_{ij} (\tilde t,x)
\bigr)
\\
&&{} + G \bigl(\tilde t\wedge T, \omega^{(s_i, y_j)}_{\cdot\wedge s_i},
v_{ij}(\tilde t, x), D v_{ij} (\tilde t,x), D^2
v_{ij} (\tilde t,x) \bigr).
\end{eqnarray*}
Now by Assumption~\ref{assum-Guniform2} and (\ref{vbound}) we have, at
$(\tilde t, x) \in\overline Q_0^\e$,
\begin{eqnarray*}
&&G \bigl(t, \omega, v_{ij}, D v_{ij}, D^2
v_{ij} \bigr) - G \bigl(\tilde t\wedge T, \omega^{(s_i, y_j)}_{\cdot\wedge s_i},
v_{ij}, D v_{ij}, D^2 v_{ij} \bigr)
\\
&&\qquad=G \bigl(t, \omega, v_{ij}, D v_{ij}, D^2
v_{ij} \bigr) - G \bigl(t, \omega^{\hat
\pi_1}_{\cdot\wedge t_1},
v_{ij}, D v_{ij}, D^2 v_{ij} \bigr)
\\
&&\qquad\quad{}+ G \bigl(t, \omega^{\hat\pi_1}_{\cdot\wedge t_1}, v_{ij}, D
v_{ij}, D^2 v_{ij} \bigr) - G \bigl(\tilde t
\wedge T, \omega^{(s_i, y_j)}_{\cdot
\wedge s_i}, v_{ij}, D
v_{ij}, D^2 v_{ij} \bigr)
\\
&&\qquad\le \rho_0 \bigl(\bigl\|\omega- \omega^{\hat\pi_1}_{\cdot\wedge
t_1}
\bigr\|_t \bigr) + \tilde \rho_0\bigl(|t-\tilde t\wedge T|\bigr) \bigl[|D
v_{ij}| +\bigl |D^2 v_{ij}\bigr| \bigr]
\\
&&\qquad\quad{}+ \rho_0 \bigl({\mathbf{d}}_\infty \bigl(\bigl(t,
\omega^{\hat\pi
_1}_{\cdot\wedge t_1}\bigr), \bigl(\tilde t\wedge T,
\omega^{(s_i, y_j)}_{\cdot\wedge s_i}\bigr) \bigr) \bigr)
\\
&&\qquad\le \rho_0(2\e) + C_{\eta_0} \tilde\rho_0(T
\delta) + \rho _0 \bigl({\mathbf{d}} _\infty \bigl(\bigl(t,
\omega^{\hat\pi_1}_{\cdot\wedge t_1}\bigr), \bigl(\tilde t\wedge T, \omega
^{(s_i, y_j)}_{\cdot\wedge s_i}\bigr) \bigr) \bigr).
\end{eqnarray*}
Thus
\begin{eqnarray*}
\cL\psi(t,\omega) &\ge& {\e\over4} - \rho_1(3T\delta)-
C_{\eta_0} \tilde\rho _0(T\delta) - \rho_0 \bigl({
\mathbf{d}}_\infty \bigl(\bigl(t, \omega ^{\hat\pi_1}_{\cdot\wedge
t_1}
\bigr), \bigl(\tilde t\wedge T, \omega^{(s_i, y_j)}_{\cdot\wedge
s_i}\bigr)
\bigr) \bigr).
\end{eqnarray*}
Substituting this inequality to (\ref{psiest2}), we see that the rest of
the proof of Lemma~\ref{lem-phie} remains the same.

%s8.3 #&#
\subsection{Concluding remarks}

We now summarize the conditions under which we have the complete
wellposedness result.

%th8.4 #&#
\begin{thmm}
\label{thmm-summary} Assume the following hold true:
\begin{itemize}
\item Assumptions \ref{assum-G} and \ref{assum-Guniform}\textup{(ii)};

\item Assumptions \ref{assum-xi0} and \ref{assum-xi2} or, more
specifically, the sufficient conditions of Lemma~\ref{lem-xi};

\item $G$ is either convex in $\g$ or the dimension $d\le2$;

\item Assumption~\ref{assum-Guniform}\textup{(i)}, or more generally,
Assumption~\ref{assum-Guniform2} and (\ref{vbound}).
\end{itemize}

Then the results of Theorem~\ref{thmm-wellposedness} hold true.
\end{thmm}

We conclude with some final remarks on our assumptions. We first note
that the highly technical requirements of the space $\ol C^{1,2}(\L)$
are needed only in the proofs, and are not part of our assumptions.
Assumptions \ref{assum-G} and \ref{assum-xi0} are more or less
standard, and are in fact the conditions used in \cite{ETZ1}. In
particular, due to the failure of the dominated convergence theorem
under $\ol\cE^{\cP_L}$, the regularity of the involved processes
become crucial, and some assumptions on regularity of data are more or
less necessary.

Assumption~\ref{assum-xi2} on the additional structure of $\xi$ is
purely technical, due to our current approach. Indeed, in situations
where we have a representation for the viscosity solution, for example,
in the semilinear case, as in \cite{ETZ1}, Section~7, this assumption
is not needed. We believe this assumption can also be removed if we
consider path-dependent HJB equations where the function $\th^\e_n$ in
Lemma~\ref{lem-the} can be constructed directly via second-order BSDEs.

The uniform continuity of $G$ in $(t,\omega)$ in Assumption~\ref
{assum-Guniform}(i) excludes the dependence of the diffusion
coefficient $\si$ on $(t,\omega)$ for stochastic control or stochastic
differential game problems (see \cite{ETZ1}, Section~4 and \cite{PZ})
and thus is not desirable. This is due to our approach of approximating
PPDEs by path-frozen PDEs. This assumption may not be needed if we do
not use this approximation.
%, for example we may extend the the arguments in \cite{EKTZ} to cover
%semilinear PPDEs with random diffusion $\si$.

The uniform nondegeneracy of $G$ in Assumption~\ref{assum-Guniform}(ii) is of course serious, as in PDE literature.
%In the semilinear case of \cite{RTZ}, the PPDE can be degenerate.
%However, in fully nonlinear case we have to study the regularity,
%which can be subtle when $G$ is degenerate as we see in Example~\ref{eg-rep}.

Finally, Assumption~\ref{assum-comparison} is crucial in our current
approach. For path-dependent HJB equations, namely when $G$ is convex
in $\g$, we have, more or less, complete results in the uniformly
nondegenerate case. However, in the present paper we verify this
assumption by the existence of classical solutions of the mollified
path-frozen PDE. Unfortunately, for Bellman--Isaacs equations, we are
able to obtain classical solutions only when $d \le2$; see \cite{PZ}.
It will be very interesting to explore more PDE estimates to see if we
can verify Assumption~\ref{assum-comparison} directly without getting
into classical solutions of high-dimensional Bellman--Isaacs equations.

We note that the essential point of our whole argument is to find
approximations $\ol u^\e, \ul u^\e\in\ol C^{1,2}(\L)$ such that
$\cL
\ol u^\e\ge0 \ge\cL\ul u^\e$. Assumptions \ref{assum-Guniform},
\ref
{assum-xi2} and \ref{assum-comparison} all serve this purpose.
There is potentially an alternative way to prove the comparison
principle directly. Let $u^1$ be a viscosity subsolution and $u^2$ a
viscosity supersolution such that $u^1_T \le u^2_T$. Instead of
mollifying the PDE to obtain classical solutions, we may try to mollify
$u^i$ directly so that the corresponding $u^{i,\e}$ will be
automatically smooth (in some appropriate sense). In fact, in the PDE
literature, the convex/concave convolution exactly serves this purpose.
However, in this case, the main challenge is that we need to check that
$u^{1,\e}$ is a classical subsolution, and $u^{2,\e}$ a classical
supersolution, which, if true, will imply the comparison immediately.
It will be interesting to explore this approach as well in future research.

% imsref loaded by akundreckaite, 2015-05-21 12:07:57
%

\printaddresses
\end{document}